%% file: main.tex
\documentclass[11pt,a4paper,reqno]{amsart}

\usepackage{times}
\usepackage{mathrsfs,mathtools}
\usepackage{enumerate}
\usepackage[inline]{enumitem}
\usepackage{stmaryrd} 
\usepackage{amsmath,amsthm,amssymb,amsfonts,accents} 
\usepackage[usenames,dvipsnames]{xcolor}
\usepackage{subfiles}
\usepackage{dsfont}
\usepackage[mathscr]{euscript}
\usepackage{graphicx}
\usepackage[colorlinks=true, linkcolor=tub-red ,citecolor=my-blue]{hyperref}
\usepackage[T1]{fontenc}
\usepackage[english]{babel}
\usepackage{todonotes}
\usepackage{geometry}
\usepackage{xspace,cancel}
\usepackage[capitalize]{cleveref}
\crefname{subsection}{Subsection}{Subsections}
\usepackage[numbers,sort&compress,longnamesfirst]{natbib}

\allowdisplaybreaks

\def \E{\mathbb{E}}

\def \N{\mathbb{N}}

\def \R{\mathbb{R}}

\def\d{\mathrm{d}}

\definecolor{my-blue}{RGB}{0,85,150}
\definecolor{tub-red}{RGB}{187,14,31}

\def\mage{\color{magenta}}

\setlist[enumerate]{label=\textup{(\roman*)},itemindent=1cm,leftmargin=0cm}
\newcommand{\oX}{\vphantom{X}\smash{\overline X}}

\makeatletter \@addtoreset{equation}{section}
\def\theequation{\thesection.\arabic{equation}}
\newtheorem{theorem}{Theorem}[section]
\newtheorem*{assumption*}{\assumptionName}
    \providecommand{\assumptionName}{}
    \makeatletter
    \newenvironment{assumption}[1]
    {
        \renewcommand{\assumptionName}{Assumption #1}
        \begin{assumption*}\protected@edef\@currentlabel{#1}%
    }
    {%
        \end{assumption*}
    }
    \makeatother

\newtheorem*{condition*}{\conditionName}
    \providecommand{\conditionName}{}
    
\newtheorem{corollary}[theorem]{Corollary}
\newtheorem{lemma}[theorem]{Lemma}
\newtheorem{proposition}[theorem]{Proposition}
\newtheorem{definition}[theorem]{Definition}
\newtheorem{remark}[theorem]{Remark}

\newtheorem*{iassumption}{Intermediate Assumption}

\theoremstyle{definition}

\newtheorem{notation}[theorem]{Notation}
\makeatletter
\def\namedlabel#1#2{\begingroup
    #2%
    \def\@currentlabel{#2}%
    \phantomsection\label{#1}\endgroup
}

\numberwithin{equation}{section}

\newcommand{\numberthis}{\addtocounter{equation}{1}\tag{\theequation}}
\newcommand{\tnorm}[1]{\left\vert\kern-0.25ex\left\vert\kern-0.25ex\left\vert #1 \right\vert\kern-0.25ex\right\vert\kern-0.25ex\right\vert}
\newcommand{\tnormbig}[1]{\bigl\vert\kern-0.25ex\bigl\vert\kern-0.25ex\bigl\vert #1 \bigr\vert\kern-0.25ex\bigr\vert\kern-0.25ex\bigr\vert}

\newcommand{\ud}{\ensuremath{\mathrm{d}}}

\newcommand{\ds}{\ud s}
\newcommand{\dx}{\ud x}
\newcommand{\abs}[1]{\left\vert#1\right\vert}
\newcommand{\norm}[1]{\left\Vert#1\right\Vert}

\newcommand{\Pred}{\mathcal{P}}
\newcommand{\borel}[1]{\mathcal{B}\left({#1}\right)}
\newcommand{\Rp}{\mathbb{R}_+}
\newcommand{\Pm}{\mathbb{P}}
\newcommand{\cadlag}{c\`adl\`ag\xspace}
\newcommand{\dC}[1]{\ud C_{#1}}
\newcommand{\e}[1]{\mathrm{e}^{#1}}
\newcommand{\FilG}{\mathbb{G}}
\newcommand{\weakFil}{\xrightarrow{\hspace*{0.2cm}\w\hspace*{0.2cm}}}

\newcommand{\LGconvmulti}[2]{\xrightarrow{\hspace{0.2cm}\mathbb{L}^{#1}(\Omega,\mathcal{G},\Pm;\RR{#2})\hspace{0.2cm}}}
\newcommand{\RR}[1]{\mathbb{R}^{#1}}

\DeclareMathOperator{\w}{w}

\DeclareFontFamily{U}{mathx}{\hyphenchar\font45}
\DeclareFontShape{U}{mathx}{m}{n}{
<5><6><7><8><9><10>
<10.95><12><14.4><17.28><20.74><24.88>
mathx10
}{}
\DeclareSymbolFont{mathx}{U}{mathx}{m}{n}
\DeclareFontSubstitution{U}{mathx}{m}{n}
\DeclareMathAccent{\widebar}{0}{mathx}{"73}

\makeatletter
\newcommand \Dotfill {\leavevmode \leaders \hb@xt@ 6pt{\hss \hss }\hfill \kern \z@}
\makeatother

\makeatletter
\def\@tocline#1#2#3#4#5#6#7{\relax
  \ifnum #1>\c@tocdepth 
  \else
    \par \addpenalty\@secpenalty\addvspace{#2}%
    \begingroup \hyphenpenalty\@M
    \@ifempty{#4}{%
      \@tempdima\csname r@tocindent\number#1\endcsname\relax
    }{%
      \@tempdima#4\relax
    }%
    \parindent\z@ \leftskip#3\relax \advance\leftskip\@tempdima\relax
    \rightskip\@pnumwidth plus4em \parfillskip-\@pnumwidth
    #5\leavevmode\hskip-\@tempdima
      \ifcase #1
       \or\or \hskip 1.65em \or \hskip 3.3em \else \hskip 4.95em \fi%
      #6\nobreak\relax
    \Dotfill
    \hbox to\@pnumwidth{\@tocpagenum{#7}}\par
    \nobreak
    \endgroup
  \fi}
\makeatother

\makeatletter
\def\l@section{\@tocline{1}{0pt}{1pc}{}{\scshape}}
\renewcommand{\tocsection}[3]{%
\indentlabel{\@ifnotempty{#2}{\ignorespaces#1 #2.\hskip 0.7em}}#3}
\def\l@subsection{\@tocline{2}{0pt}{1pc}{5pc}{}}

\def\l@subsubsection{\@tocline{3}{0pt}{1pc}{7pc}{}}

\makeatother

\begin{document}

\title[Stability of Backward Stochastic Differential Equations]{Stability of Backward Stochastic Differential Equations:\\ the general  Lipschitz case}

\author[A. Papapantoleon]{Antonis Papapantoleon}
\author[D. Possama\"{i}]{Dylan Possama\"{i}}
\author[A. Saplaouras]{Alexandros Saplaouras}

\address{Delft Institute of Applied Mathematics, TU Delft, 2628 Delft, The Netherlands \& Department of Mathematics, NTUA, 15780 Athens, Greece \& Institute of Applied and Computational Mathematics, FORTH, 70013 Heraklion, Greece}
\email{papapan@mail.ntua.gr, a.papapantoleon@tudelft.nl}

\address{Department of Mathematics, ETH Z\"urich, R\"amistrasse 101, 8092 Z\"urich, Switzerland}
\email{dylan.possamai@math.ethz.ch}

\address{Department of Mathematics, NTUA, Zografou Campus, 15780 Athens, Greece}
\email{alsapl@mail.ntua.gr}

\thanks{Antonis Papapantoleon gratefully acknowledges the financial support from the Hellenic Foundation for Research and Innovation Grant No. HFRI-FM17-2152. 
Dylan Possama\"i gratefully acknowledges the financial support from the ANR project {\sc PACMAN} (ANR-16-CE05-0027). 
Alexandros Saplaouras gratefully acknowledges the financial support from the DFG Research Training Group 1845 ``Stochastic Analysis with Applications in Biology, Finance and Physics''. 
Moreover, all authors gratefully acknowledge the financial support from the {\sc Procope} project ``Financial markets in transition: mathematical models and challenges''.}


\subjclass[2010]{}

\date{}

\begin{abstract}
In this paper, we obtain stability results for backward stochastic differential equations with jumps (BSDEs) in a very general framework.
More specifically, we consider a convergent sequence of standard data, each associated to their own filtration, and we prove that the associated sequence of (unique) solutions is also convergent.
The current result extends earlier contributions in the literature of stability of BSDEs and unifies several frameworks for numerical approximations of BSDEs and their implementations. 
\end{abstract}


\maketitle

\frenchspacing

\tableofcontents

\section{Introduction}
\subfile{Introduction}

\section{Framework}\label{sec:Framework}
\subfile{Framework}

\section{Stability of backward stochastic differential equations}\label{sec:MainTheorem}
\subfile{Stability_BSDEJ_Theorem}

\subsection{Outline of the proof}\label{subsec:BodyProof}
\subfile{Stability_BSDEJ_Proof_Main_Body}
\subsection{The first step of the induction is valid}\label{subsec:FirstStepInduction}
\subfile{Stability_BSDEJ_Induction_FirstStep}

\subsection{The \texorpdfstring{$p$-th}{p-th} step of the induction is valid}\label{subsec:p-step}
\subfile{Stability_BSDEJ_Induction_pStep_Remarks}

\subfile{Stability_BSDEJ_Induction_pStep_Convergence}

\subfile{Stability_BSDEJ_Induction_pStep_UI}

\subsection{On the nature of the conditions}\label{subsec:CommentsonConditions}
\subfile{Stability_BSDEJ_Comments_on_Framework}

\section{Examples and applications}\label{subsec:Examples}
\subfile{Examples}

\appendix

\section{Auxiliary results}\label{AppendixAux}

\subfile{Appendix_Notation}
\subfile{Appendix_SBSDE_2}
\subfile{Appendix_SBSDE_1}

\bigskip
\bibliographystyle{plainnat}
{\small
\bibliography{bibliographyDylan}}

\end{document}

%% file: Introduction.tex

The goal of this paper is to provide a suitable framework under which the stability {\color{black}-- sometimes also called robustness---}property of well-posed backward stochastic differential equations with jumps (BSDEs for short) is valid.
More precisely, {\color{black} consider} the \emph{standard data} $\mathscr{D}:=(\mathbb{G}, \tau,(X^{\circ},X^{\natural}),C,\xi,f)$ {\color{black} of a BSDE---\textit{i.e.} assuming that a probability space is given, $\mathbb{G}$ is a filtration, $\tau$ is the stopping time serving as terminal horizon for the BSDE, $X^{\circ}$ and $X^{\natural}$ are square-integrable martingales which are the integrators of the corresponding stochastic integrals, $C$ is an increasing predictable process playing the role of Lebesgue--Stieltjes integrator, $\xi$ is the terminal condition of the BSDE, and $f$ its generator---see }\citeauthor*{papapantoleon2016existence} \cite[Definitions 3.2 and 3.16]{papapantoleon2016existence}, then we know from \cite[Theorems 3.5 and 3.23]{papapantoleon2016existence} that the BSDE   
\begin{align}\label{BSDE-Integral}
Y_t = \xi + \int_{(t,{\tau}]} f\big(s, Y_{s-},Z_s,U(s,\cdot)\big)\ud C_s
- \int_t^{\tau} Z_s \ud X^{\circ}_s - \int_t^{\tau}\int_{\mathbb{R}^n} U(s,x)\widetilde{\mu}^{X^{\natural}}(\ds,\dx) - \int_t^{\tau} \ud N_s,
\end{align}
is well-posed, \emph{i.e.} it has a unique solution $\mathscr{S}:=(Y,Z,U,N)$ in appropriate spaces.
Assuming now that we are given a sequence of standard data $(\mathscr{D}^{k})_{k\in\mathbb{N}\cup\{\infty\} }$, whose associated solutions form the sequence $(\mathscr{S}^{k})_{k\in\mathbb{N}\cup\{\infty\}}$, we address the following question:
\[\text{ `How should the sequence $(\mathscr{D}^{k})_{k\in\mathbb{N}\cup\{\infty\} }$ converge, so that the sequence $(\mathscr{S}^{k})_{k\in\mathbb{N}\cup\{\infty\} }$ converges as well?'}\]
Our aim is to make precise the way(s) convergence should be understood in the previous statement, and to provide a framework which is as general as possible, while the conditions are as weak as possible.
Naturally, the existing stability results in the literature should then be recovered as special cases.
In particular, numerical schemes---and in a loose sense their implementations---for BSDEs with stochastic Lipschitz generators will be considered as sub-cases of our main results (see \cref{BSDERobMainTheorem}).
We would like to point out the fact that the techniques we use are purely probabilistic, so that the numerical schemes that can potentially be derived from our results will not rely on the connection of BSDEs to partial integro--differential equations (PIDEs for short). 
An important consequence of this fact is that our results potentially provide stochastic numerical methods that could be an alternative to finite difference or finite element schemes for the solution of semi-linear PIDEs.
This could be particularly interesting for high-dimensional non-linear PIDEs arising, for example, in the presence of valuation adjustments in option pricing; see \textit{e.g.} \citeauthor*{crepey2015bilateral} \cite{crepey2015bilateral,crepey2015bilateral2} or \citeauthor*{bichuch2018arbitrage} \cite{bichuch2018arbitrage}.

\medskip
Let us start by providing some historical references for the problem we are interested in. 
When the limit-BSDE---that is, the one that corresponds to the standard data $\mathscr{D}^{\infty}$---is solely driven by a Brownian motion, the articles of \citeauthor*{briand2002robustness} \cite{briand2001donsker,briand2002robustness} provide a suitable framework for the stability property to hold. 
It is noteworthy that in these articles, the filtration $\mathbb{G}^k$ is neither required to coincide with $\mathbb{G}^{\infty}$, nor to be a discretisation of $\mathbb{G}^{\infty}$, which, as far as we know, makes \cite{briand2002robustness} the most general result for the stability of Brownian-driven BSDEs. 
More precisely, in \cite{briand2001donsker}, respectively in \cite{briand2002robustness}, the authors approximate the Brownian motion driving the BSDE by a sequence of scaled random walks, respectively by square-integrable martingales.  
In both articles, the sequence of filtrations associated to the It\^o integrators weakly converges to the Brownian filtration, thus allowing for the aforementioned generality.
The earlier work of \citeauthor*{hu1997stability} \cite{hu1997stability} falls in the category where $\mathbb{G}^{k}=\mathbb{G}^{\infty}$, for every $k\in\mathbb{N}$, and the stability is investigated with respect to the pair $(\xi,f)$ only.
The articles of \citeauthor*{ma2002numerical} \cite{ma2002numerical}, \citeauthor*{toldo2006stability} \cite{toldo2006stability}, \citeauthor*{geiss2020random} \cite{geiss2020random}, and \citeauthor*{briand2021donsker} \cite{briand2021donsker}, follow the spirit of \cite{briand2001donsker}, \emph{i.e.}, the Brownian motion is approximated by a sequence of random walks, and respectively study numerical schemes, BSDEs with random horizon, H\"older-continuous terminal conditions, and rates of convergence. {Remaining in the Brownian framework for the limit-BSDE, \citeauthor*{jusselin2019scaling} \cite{jusselin2019scaling} provide stability results when the limit-BSDE is driven by a Brownian martingale and is approximated by scaled point processes.
}

\medskip
When the sequence $(\mathbb{G}^{k})_{k\in\mathbb{N}}$ corresponds to refined discretizations of $\mathbb{G}^{\infty}$, then we are essentially dealing with an Euler scheme for the BSDE itself.  
It is then noteworthy that the corresponding sequence of discrete filtrations weakly converges to the original filtration; see \citeauthor*{coquet2001weak} \cite[Proposition 3 or Proposition 4.A]{coquet2001weak}. In this line of literature, the majority of the articles consider the case where the BSDE is driven by a Brownian motion, and only a few deal with a more general case, \textit{e.g.} BSDEs driven by L\'evy processes. 
Notably, only the articles of \citeauthor*{bouchard2008discrete} \cite{bouchard2008discrete}, \citeauthor*{aazizi2013discrete} \cite{aazizi2013discrete} (in the pure jump case), \citeauthor*{lejay2014numerical} \cite{lejay2014numerical}, \citeauthor*{geiss2016simulation} \cite{geiss2016simulation} (in \cite{lejay2014numerical,geiss2016simulation} the jump part of the driving martingale is a Poisson process), \citeauthor*{kharroubi2015decomposition} \cite{kharroubi2015decomposition}, with a jump process depending on the Brownian motion itself, \citeauthor*{madan2015convergence} \cite{madan2015convergence} (which follows in spirit the approach of \cite{briand2001donsker}) and \citeauthor*{dumitrescu2016reflected} \cite{dumitrescu2016reflected} (where the jump part of the driving martingale is a Poisson process and the authors actually consider reflected BSDEs) consider BSDEs which include stochastic integrals with respect to an integer-valued measure, associated to the jumps of a L\'evy process. See also \citeauthor*{khedher2016discretisation} \cite{khedher2016discretisation} for BSDEs driven by c\`adl\`ag martingales.

\medskip
When the martingale driving the BSDE is only an It\^o integral, the literature is very rich and there are mainly two approaches: one that relies on purely probabilistic techniques, and another that relies on the connection of the BSDEs driven by Brownian motion with second-order semi-linear parabolic PDEs. 
In the former approach, the first work that dealt with the stability of BSDEs is due to \citeauthor*{bally1997approximation} \cite{bally1997approximation}, who provides an approximation scheme for a BSDE (driven by a Brownian motion) based on a time-discretisation scheme constructed on a Poisson net. Contemporary contributions are due to \citeauthor*{chevance1997numerical} \cite{chevance1997numerical} and \citeauthor*{coquet1998stability} \cite{coquet1998stability}.
With the exception of some of the aforementioned articles, the majority of numerical approximations for BSDEs are implementations of the standard Euler scheme. 
Indicatively, and with no claim of comprehensiveness, we can mention
\citeauthor*{bender2007forward} \cite{bender2007forward}, 
\citeauthor*{bender2010importance} \cite{bender2010importance}, 
\citeauthor*{bender2008time} \cite{bender2008time}, 
\citeauthor*{bouchard2009strong} \cite{bouchard2009strong}, 
\citeauthor*{bouchard2004discrete} \cite{bouchard2004discrete}, 
\citeauthor*{bouchard2009discrete} \cite{bouchard2009discrete}, 
\citeauthor*{briand2014simulation} \cite{briand2014simulation}, 
\citeauthor*{chassagneux2014runge} \cite{chassagneux2014runge}, 
\citeauthor*{crisan2010solving} \cite{crisan2010solving,crisan2012solving,crisan2014second}, 
\citeauthor*{crisan2010monte} \cite{crisan2010monte}, 
\citeauthor*{gobet2010solving} \cite{gobet2010solving}, 
\citeauthor*{gobet2008numerical} \cite{gobet2008numerical}, 
\citeauthor*{gobet2010l2} \cite{gobet2010l2}, 
\citeauthor*{gobet2015analytical} \cite{gobet2015analytical}, 
\citeauthor*{gobet2016approximation} \cite{gobet2016approximation,gobet2017adaptive}, 
\citeauthor*{gobet2005regression} \cite{gobet2005regression}, 
\citeauthor*{hu2011malliavin} \cite{hu2011malliavin},
\citeauthor*{pages2018improved} \cite{pages2018improved}, 
\citeauthor*{zhang2004numerical} \cite{zhang2001some,zhang2004numerical}.
Numerical solutions of BSDEs via their connection to branching processes have also been recently explored, and correspond to an alternative incarnation of the probabilistic approach, see \textit{e.g.} 
\citeauthor{bouchard2019numerical} \cite{bouchard2019numerical},
\citeauthor{bouchard2016numerical} \cite{bouchard2016numerical},
\citeauthor{henry2012counterparty} \cite{henry2012counterparty},
\citeauthor{henry2016branching} \cite{henry2016branching},
\citeauthor{henry2014numerical} \cite{henry2014numerical}.

\medskip
In the latter approach, the work of \citeauthor*{ma1994solving} \cite{ma1994solving}, which deals with the existence and uniqueness of the solution of a (coupled) forward--backward stochastic differential equation (FBSDE) as well as the stability of a parametrised family, with its proposed `four-step scheme' constitutes the cornerstone of the PDE approach. 
The papers by  
\citeauthor*{bouchard2004discrete} \cite{bouchard2004discrete}, 
\citeauthor*{douglas1996numerical} \cite{douglas1996numerical}, 
\citeauthor*{gobet2007error} \cite{gobet2007sequential,gobet2007error},
\citeauthor*{milstein2007discretization} \cite{milstein2007discretization}, 
are implementations of numerical schemes based either on the four-step scheme or PDE arguments.   

\medskip
Let us also mention that recently, new implementations have been proposed based on modern techniques, such as the use of deep neutral networks, \emph{e.g.}, see 
\citeauthor*{e2017deep} \cite{e2017deep}, 
\citeauthor*{e2019multilevel} \cite{e2019multilevel},
\citeauthor*{beck2019machine} \cite{beck2019machine},
\citeauthor*{hure2020deep} \cite{hure2019some,hure2020deep},
\citeauthor*{germain2020deep} \cite{germain2020deep,germain2021neural},
and of parallel programming, \emph{e.g.}, see \citeauthor*{gobet2016stratified} \cite{gobet2016stratified} or \citeauthor*{gobet2020quasi} \cite{gobet2020quasi} and also \citeauthor*{abbas2018xva} \cite{abbas2018xva} for an application to XVA computations.

\medskip
A closer examination of the existing literature indicates that there exists no general (at least in the spirit of \citet{briand2002robustness}) stability result for processes with jumps, even in the L\'evy case. 
Our paper aspires to fill this gap, and to describe and study a framework for the stability property for a sequence of BSDEs driven by square-integrable martingales, covering thus also the subclass of L\'evy drivers. 
To this end, we mainly have two mathematical tools at hand: the Moore--Osgood theorem, and a general result for the stability of martingale representations. 
The former was introduced and successfully used by \citeauthor*{briand2001donsker} \cite{briand2001donsker,briand2002robustness} in order to efficiently control the doubly indexed sequence obtained by the sequence of Picard schemes; we will follow the same approach.
The latter is presented in \citeauthor*{papapantoleon2019stability} \cite{papapantoleon2019stability}. 
It comes as no surprise that the stability of martingale representations comes into play, since the existence and uniqueness of the solution of the BSDE \eqref{BSDE-Integral} with stochastically Lipschitz generator is obtained by means of a martingale representation argument.
After all, the BSDE \eqref{BSDE-Integral} is frequently referred to as the $f$--\emph{nonlinear martingale representation} of $\xi$.
Additionally to these main tools, our toolbox includes \cref{CharactWeakConv} and its offsprings. 
\cref{CharactWeakConv} provides a characterisation of weak convergence of measures to an atomless measure on the real (half-)line and its corollaries are going to guarantee the convergence of the Lebesgue--Stieltjes integrals when the integrators are regular enough.
In view of the above brief comments, it is natural to combine the frameworks of \cite{papapantoleon2016existence,papapantoleon2019stability}\footnote{To the best of our knowledge, the most general well-posedness results for {\color{black} Lipschitz} BSDEs are given in \cite{papapantoleon2016existence}.} and to properly enrich them in order to obtain the framework suitable for our purpose.
We postpone the detailed description and the required technicalities until \Cref{subsec:BodyProof}. 
However, we would like to underline that if we were to derive the unique solution of BSDE \eqref{BSDE-Integral} by means of a different method compared to the two referred above, we can, \emph{mutatis mutandis}, adapt the arguments and derive the stability property as soon as the alternative method is based on a fixed-point argument, \emph{i.e.} a Picard scheme.
This would require to introduce appropriate (uniform) integrability conditions related to the existence and uniqueness of the solution of BSDE \eqref{BSDE-Integral}.

\medskip


Backward SDEs have been used for modelling a vast array of random phenomena in the areas of economics, finance, game theory, and others, and their numerical solutions is a topic of significant importance both in theory and in practice. {\color{black}In particular, allowing for random, and possibly infinite horizon, has found applications in optimal control and optimal stopping theory, see \citeauthor*{hu2011some} \cite{hu2011some}, or even more recently in contract theory, see \citeauthor*{sannikov2008continuous} \cite{sannikov2008continuous}, \citeauthor*{pages2014mathematical} \cite{pages2014mathematical}, \citeauthor*{lin2020random} \cite{lin2020random}. Notice that it is typical with unbounded horizon to allow for stochastic Lipschitz conditions for integrability purposes, but we also emphasise that such BSDEs also naturally arises when one considers Malliavin derivatives of the solutions (see for instance \citeauthor*{imkeller2012differentiability} \cite{imkeller2012differentiability}), the latter being useful to study for instance existence of densities, see \citeauthor*{mastrolia2014malliavin} \cite{mastrolia2014density,mastrolia2014malliavin}.}
Another important application of the results of this paper is in the area of numerical schemes for BSDEs, where Theorem \ref{BSDERobMainTheorem} provides a flexible and general framework for automatically deducing the convergence of such numerical schemes. 
These can be viewed as generalisations of the celebrated Donsker's theorem.
In that case, the martingales of the sequence that approximates the limit-BSDE have to be defined with respect to their own filtrations, hence the requirement of the convergence of filtrations becomes necessary.
This creates additional technical difficulties in the presence of jumps in the limit-BSDE, \textit{e.g.} for the convergence of Lebesgue--Stieltjes integrals as well as for the required uniform integrability, that have to be carefully handled.


\medskip
The structure of the paper is as follows.
In \cref{sec:Framework} we provide the set of conditions a sequence of standard data should satisfy 
for existence and uniqueness of the solution
alongside notation and some helpful comments.
In \cref{sec:MainTheorem} the theorem for the stability of BSDEs is stated, while in \Cref{subsec:BodyProof} we provide the sketch of the proof by describing the main arguments.
There, it will become evident that the stability of martingale representations plays an important role in obtaining the stability property of BSDEs, as we have already stated.
In \Cref{subsec:FirstStepInduction} and \Cref{subsec:p-step} the required technical lemmata, which verify our claims in \Cref{subsec:BodyProof}, are presented. 
In \Cref{subsec:CommentsonConditions} we briefly discuss the nature of the imposed conditions.
Finally, in \Cref{subsec:Examples} we present some examples that demonstrate the generality and applicability of \Cref{BSDERobMainTheorem}. 
In \Cref{AppendixAux} several auxiliary results are proved, and we expand upon the notation used in the paper.

\subsection*{General notation} 

Let $\mathbb R_+$ denote the set of non-negative real numbers, $\overline{\mathbb{R}}_+:=\mathbb{R}_+\cup\{\infty\}$ and $\vert x\vert$ denote the absolute value of a real number $x\in\R$. 
We denote the set of positive natural numbers by $\mathbb{N}$ and we define $\overline{\mathbb{N}}:=\mathbb{N}\cup\{\infty\}.$ 
For two arbitrary $(p,q)\in\mathbb{N}\times\mathbb{N}$, we identify $\R^{p\times q}$ with the set of matrices with $p$ rows and $q$ columns and real entries. The transpose of $v\in\mathbb{R}^{p\times q}$ will be denoted by $v^\top\in\mathbb{R}^{q\times p}$.
The element at the $i$-th row and $j$-th column of $v$ will be denoted by $v^{ij}$, for $ i\in\{1,\dots,p\}$ and $j\in\{1,\dots,q\}$, and it will be called the $(i,j)$-element of $v$.
However, the notation needs some care when we deal with sequences of elements of $\mathbb{R}^{p\times q}$, \emph{e.g.} if $(v^k)_{k\in\overline{\mathbb{N}}}\subset \mathbb{R}^{p\times q}$, then we will denote by $v^{k,ij}$ the $(i,j)$-element of $v^k$, for $ i\in\{1,\dots,p\}$, $j\in\{1,\dots,q\}$, and $k\in\overline{\mathbb{N}}$.
The trace of a square matrix $w\in\mathbb{R}^{p\times p}$ is given by $\textrm{Tr}[w]:=\sum_{i=1}^pw^{ii}.$
We endow $\mathbb{R}^{p \times q}$ with the norm defined for any $z\in\mathbb{R}^{p \times q}$ by $\Vert z\Vert_2^2:=\textrm{Tr} [z^\top z]$ and remind the reader that this norm derives from the inner product defined for any $(z,u)\in\mathbb{R}^{p \times q}\times\mathbb{R}^{p \times q}$ by $\textrm{Tr}[z^{\top}u]$. 
We will additionally endow the space $\mathbb{R}^{p\times q}$ with the norm $\Vert \cdot\Vert_1$, which is defined by $\Vert v\Vert_1:=\sum_{i=1}^p\sum_{j=1}^q|v^{ij}|$.
We identify $\mathbb{R}^p$ with $\mathbb{R}^{p\times 1}$, \emph{i.e.} an arbitrary $x\in\mathbb{R}^p$ will be identified as a column vector of length $p$.
If $(x^k)_{k\in\overline{\mathbb{N}}}\subset\mathbb{R}^p$, then $x^{k,i}$ denotes the $i$-th element of $x^k$, for $i\in\{1,\dots,p\}$ and $k\in\overline{\mathbb{N}}.$

\medskip
{We abuse notation} and denote by $0$ the neutral element in the group $(\mathbb R^{p\times q},+)$. 
Furthermore, for any finite dimensional topological space $E$, $\mathcal B(E)$ will denote the associated Borel $\sigma$-algebra. 
In addition, for any other finite dimensional space $F$, and for any non-negative measure $\nu$ on $\big(\mathbb R_+\times E,\mathcal B(\mathbb R_+)\otimes\mathcal B(E)\big)$, we will denote the Lebesgue--Stieltjes integrals of any measurable map $f:\big(\mathbb R_+\times E,\mathcal B(\mathbb R_+)\otimes\mathcal B(E)\big)\longrightarrow (F,\mathcal B(F))$, by
\begin{gather*}
\int_{(u,t]\times A}f(s,x)\nu(\ud s,\ud x), \text{ for any }  (t,A)\in\mathbb R_+\times\mathcal B(E),\; \text{and}\; \int_{(u,\infty)\times A}f(s,x)\nu(\ud s,\ud x),  \text{ for any }   A\in\mathcal B(E),
\end{gather*}
where the integrals are to be understood in a component-wise sense. 
Finally, the letters $p$, $q$, $i$, $j$, $k$, $l$, $\ell$, $m$, and $n$ are reserved to denote arbitrary positive integers.
Specifically, the calligraphic letter $\ell$ will denote the dimension of the state space of a solution of a BSDE,
$m$ will denote the dimension of the state space of an It\^o integrator, and $n$ will denote the dimension of the state space of a process associated to an integer-valued random measure.

\medskip
Let us define the maps $\mathbb R^{n}\ni x\overset{\textrm{Id}}{\longmapsto}x\in\mathbb R^{n}$ and $\mathbb R^{n}\ni x\overset{\textrm{q}}{\longmapsto}x x^\top\in\mathbb R^{n\times n}$, where we suppress the dependence on the dimension $n$ for notational simplicity.
The space of functions defined on $\mathbb{R}_+$ and with state space $\mathbb{R}^{p\times q}$, which are right-continuous with left-limits (c\`adl\`ag) will be denoted by $\mathbb{D}(\mathbb{R}_+;\mathbb{R}^{p\times q})$.
The metric induced by the $\textup{J}_1$-topology will be denoted by $\delta_{\textup{J}_1}$, while the supremum norm will be denoted by $\Vert \cdot \Vert_\infty$. 
Finally, we will denote the locally uniform convergence of the sequence $(\alpha^{k})_{k\in\overline{\mathbb{N}}}\subset \mathbb{D}(\mathbb{R}_+;\mathbb{R}^{p\times q}) $ by $\alpha^{k}\xrightarrow[k\to\infty]{\textup{lu}}\alpha^{\infty}$. 

\medskip
In order to simplify the presentation and minimise the introductory remarks, we will adopt the notation and definitions from \cite{papapantoleon2019stability}. 
Moreover, definitions related to BSDEs from \cite{papapantoleon2016existence} are also adopted.
Nonetheless, for ease of reference, we present in \cref{AppendixAux} the main ones that will be used throughout this work.
For example, given the probability space $(\Omega,\mathcal{F},\mathbb{P})$, the expectation under $\mathbb{P}$ will be denoted by $\mathbb{E}[\cdot]$. 
Another example is the following: given a filtration $\mathbb{F}$, a stopping time $\tau$ and $X\in \mathcal{H}^{2}(\mathbb{F},\tau;\mathbb{R})$, the space of It\^o integrands is 
\begin{align*}
\mathbb H^2(X,\mathbb{F},\tau;\mathbb{R}^{\ell \times m}) 
  :=\bigg\{
    Z:(\Omega\times\mathbb{R}_+,\mathcal{P}^\mathbb{F}) \longrightarrow (\mathbb{R}^{\ell \times m},\mathcal B(\mathbb{R}^{\ell \times m})):
    \mathbb{E}\bigg[
    \int_{(0,\tau]} \textrm{Tr}\Big[Z_t \frac{\ud \langle {X} \rangle^{\mathbb{F}}_t}{\ud C^{X}_t} Z_t^{\top}\Big]\ud C^{X}_t
    \bigg]<\infty
    \bigg\}.
\end{align*}

In a nutshell, we have adopted the notation
from \cite{papapantoleon2019stability} and adapted the dimension of the state space; the results of \cite{papapantoleon2019stability} remain valid for multi-dimensional martingale representations.

%% file: Framework.tex

Let us fix the probability space $(\Omega,\mathcal G,\mathbb P)$ for the remainder of this article, as well as a sequence $(\mathscr{D}^k)_{k\in\overline{\mathbb{N}}}$\label{fix-SD} where, for any $k\in\overline\N$, $\mathscr{D}^k:=(\oX^k,\mathbb{G}^k, \tau^k,C^k, \xi^k,f^k)$ are standard data under (a sufficiently large) $\hat{\beta}>0$\footnote{\label{foot:beta-hat}$\hat{\beta}$ is independent of $k\in\overline{\mathbb{N}}$. We will refer to it as \emph{the common value}.} in the sense of \citeauthor*{papapantoleon2016existence} \cite[Section 3.1]{papapantoleon2016existence}.
{\color{black} More precisely, for every $k\in\overline{\mathbb{N}}$, they satisfy the following conditions:

\begin{enumerate}[label={(F\arabic*)},leftmargin=*,itemindent=0.0cm]
	\item\label{Fone} The martingale $\oX^{k}:=(X^{k,\circ}, X^{k,\natural})$ belongs to $\mathcal H^{2,k}_{\hat{\beta}}(\mathbb R^m)\times \mathcal H^{2,k}_{\hat{\beta}}(\mathbb R^n)$, with $X^{k,\natural}$ being purely discontinuous, and $(\oX^{k},C^{k})$ satisfies Assumption \ref{assumptionC} (see \cref{appendix-A1}).
	\item\label{Ftwo} The terminal condition satisfies $\xi^{k}\in \mathbb{L}^{2,k}_{\hat{\beta}}$. 
	\item\label{Fthree} The generator $f^{k}:{\Omega\times\mathbb{R}_{+}\times\R^d\times\R^{d\times m}\times\mathfrak H^{k}\longrightarrow\mathbb{R}^d}$ is such that for any 
				$(y,z,u)\in\R^d\times \R^{d\times m}\times \mathfrak H^{k}$, 
				the map 
				$$(t,\omega)\longmapsto f^{k}(t,\omega,y,z,u_t(\omega;\cdot))\ \text{is $\mathcal{G}^{k}_t\otimes\mathcal B([0,t])-$ measurable.}$$ 
				Moreover, $f^{k}$ satisfies a stochastic Lipschitz condition, \textit{i.e.} there exist 
				\[
					r^{k}:\left(\Omega\times\mathbb{R}_+,\Pred\right)\longrightarrow(\R_+,\borel{\R_+})\ \textrm{ and }\
					\vartheta^{k} =(\theta^{k,\circ},\theta^{k,\natural}):\left(\Omega\times\mathbb{R}_+,\Pred\right)\longrightarrow(\Rp^2,\borel{\Rp^2}),
				\] such that, for $\ud C^{k} \otimes \ud \Pm-a.e.$ $(t,\omega)\in\mathbb R_+\times\Omega$
				\begin{align}\label{fstochLip}
				\begin{split}
					& \quad \left|f^{k}(t,\omega,y,z,u_t(\omega;\cdot))-f^{k}(t,\omega,y',z',u_t'(\omega;\cdot))\right|^2\\
							&\hspace{2em}\leq 	r^{k}_t(\omega) |y-y'|^2 	
							+ \theta^{k,\circ}_t(\omega)\Vert c_t(\omega) (z-z')\Vert^2 
							+ \theta^{k,\natural}_t(\omega) \left(\tnorm{u_t(\omega;\cdot)-u'_t(\omega;\cdot)}_{k,t}(\omega)\right)^2.		
				\end{split}
				\end{align}
	\item\label{Ffour} Let 
				$(\alpha^k_\cdot)^{2}:= \max\{\sqrt{r_\cdot}^{k}, \theta^{k,\circ}_\cdot,\theta^{k,\natural}_\cdot\}$
				and define the increasing, $\mathbb G^{k}-$predictable and \cadlag \ process
				\begin{equation}\label{eq:A}
					A_\cdot^{k} :=\int_0^{\cdot} (\alpha_s^{k})^2 \dC{s}^{k}.
				\end{equation}
				\noindent Then there exists $\Phi^{k}>0$ such that
				\begin{equation}\label{Ffouruniformjumps}
					\Delta A_t^{k}(\omega) \leq \Phi^{k}, \ \text{for $\ud C^{k} \otimes \ud \Pm-a.e.$ $(t,\omega)\in\mathbb R_+\times\Omega$}.
				\end{equation}
	\item\label{Ffive} It holds
 		 	\[
 				\E\left[\int_0^T\e{\hat{\beta} A^{k}_t} \frac{\abs{f^{k}(t,0,0,\mathbf 0)}^2}{(\alpha^{k})^2_t} \ \dC{t}^{k}\right] <\infty,
 			\]
 	where $\mathbf 0$ denotes the null application from $\mathbb R^n$ to $\mathbb R$.		
\end{enumerate}
These conditions guarantee that the BSDE \eqref{k-BSDEintegralform} has a unique solution for every $k\in\overline{\mathbb N}$; see also Remark \ref{rem:BSDE-solution} below for more detailed statements. }


\medskip
At this point let us clarify that the processes associated to the stochastic basis $(\Omega,\mathcal{G},\mathbb{G}^k,\mathbb{P})$ will be stopped at $\tau^k$, for each fixed $k\in\overline{\mathbb{N}}$.
If, for example, $L^k$ is a $\mathbb{G}^k$-adapted process, for some $k\in\overline{\mathbb{N}}$, then it will be assumed to be stopped at time $\tau^k$. 
In particular, the process $\oX^k$ is stopped at $\tau^k$, for every $k\in\overline{\mathbb{N}}.$

\smallskip
\begin{notation}
Some notational simplifications and rules are in order:
\begin{enumerate}[itemindent=0.cm,leftmargin=*]
\item[\textbullet] we simplify the notation of \cite{papapantoleon2019stability} associated to integer-valued measures. The integer-valued measure $\mu^{X^{k,\natural}}$ will be denoted by $\mu^{k,\natural}$, its $\mathbb{G}^k$-compensator $\nu^{(X^{k,\natural},\mathbb{G}^k)}$ will be denoted by $\nu^{k,\natural}$ and the compensated integer-valued measure 
$\widetilde{\mu}^{(X^{k,\natural},\mathbb{G}^k)}$ will be denoted by $\widetilde{\mu}^{k,\natural}$, for every $k\in\overline{\mathbb{N}}$;
\item[\textbullet] let $k\in\overline{\mathbb{N}}$ and $Y$ be a $\mathbb{G}^{k}$-semimartingale. 
In the notation of its quadratic covariation $[Y]^{\mathbb{G}^{k}}$, resp. its predictable quadratic covariation $\langle Y \rangle^{\mathbb{G}^{k}}$, we will suppress the dependence on the filtration, and we will simply write $[Y]$, resp. $\langle Y \rangle$;
%
%
\item[\textbullet] for every $k\in\overline{\mathbb{N}}$, the notation for the spaces of \cite[Section 2.3]{papapantoleon2016existence} associated to the standard data $\mathscr{D}^k$ has been the simplified notations introduced in \cite[page 16]{papapantoleon2016existence} extended as follows: the index $k$ will be affixed, succeeding the number $2$ (if any) and preceding any other symbol (if any), \emph{e.g.}, $\mathbb H^{2,k,\natural}_{\beta}$, $\mathfrak{H}^{k}_{t,\omega}$, and so on.
The rule for the norms makes the index $k$ a subscript which precedes the value $\beta$ or $t$, \emph{e.g.}, $ \Vert \cdot\Vert_{\star,k,\beta}$, $\tnorm{\cdot}_{k,t}$. 
\end{enumerate}
\end{notation}

\begin{notation}\label{notation:tildeE-Cc}
We introduce some auxiliary notation for two subsets of $E:=\mathbb{R}_+\times\mathbb{R}^n$ as well as for sets of continuous functions with compact support.
Let $(F, \Vert \cdot\Vert_2), (G, \Vert \cdot\Vert_2)$ be Euclidean finite-dimensional spaces.
\begin{enumerate}[itemindent=0.cm,leftmargin=*]
\item[\textbullet] Set $E_0:=(\mathbb{R}_+\times\{0\})\cup(\{0\}\times\mathbb{R}^n)$ and $\widetilde{E}:=(\mathbb{R}_+\times\mathbb{R}^n)\setminus E_0$.
\smallskip
\item[\textbullet] Let $f:(F,\Vert \cdot\Vert_2)\rightarrow (G, \Vert \cdot\Vert_2)$. 
	The support of the function $f$ is the set $\textup{supp}(f):=\overline{\big\{s\in F: f(s)\neq 0 \big\}}^{\Vert \cdot\Vert_2}$, where $\overline{A}^{\Vert \cdot\Vert_2}$ denotes the closure of the set $A$ under the metric associated to the norm $\Vert \cdot\Vert_2$.
\smallskip%
\item[\textbullet] $C_{c}(F;G):=
	\big\{f:(F, \Vert \cdot \Vert_2) \longrightarrow (G, \Vert \cdot \Vert_2):	f \text{ is continuous with compact support}\big\}$.
\smallskip%
\item[\textbullet] $C_{c|\widetilde{E}}(E;\mathbb{R}^{\ell}):=
	\big\{f \in C_c(E;\mathbb{R}^{\ell}):\textup{supp}(f) \subset \widetilde{E}\big\}$.
\smallskip%
\item[\textbullet] $D^{\circ,\ell\times m}\subset C_{c}(\mathbb{R}_+;\mathbb{R}^{\ell\times m})$ is a fixed, countable set, which is dense in $\mathbb{L}^{2}\big(\mathbb{R}_+,\mathcal{B}(\mathbb{R}_+),\langle X^{\infty,\circ}\rangle(\omega)\big)$ for $\mathbb{P}$--a.e. $\omega \in \Omega$; for the existence of such a set the reader may consult \cref{lem:Construction_DCirc}.
\smallskip%
\item[\textbullet] $D^{\natural}\subset C_{c|\widetilde{E}}(E;\mathbb{R}^\ell)$ is a fixed, countable set, which is dense in 
	$\mathbb{L}^2\big(E,\mathcal{B}(E), \nu^{\infty,\natural}(\omega); \mathbb{R}^{\ell}\big)$ for $\mathbb{P}$--a.e. $\omega \in \Omega$ (see \cref{lem:Construction_DNatural}), keeping in mind that $\nu^{\infty,\natural}(E_0)=0$, $\mathbb{P}\text{\rm--a.s.}$, as well as 
	$\int_E \Vert x \Vert_2^{2}\, \nu^{\infty,\natural}(\textup{d}s,\textup{d}x)<\infty$, $\mathbb{P}\text{\rm--a.s.}$
\end{enumerate}
\end{notation}
We proceed now with the description of the conditions for the sequence $(\mathscr{D}^{k})_{k\in\overline{\mathbb{N}}}$ to be convergent, \emph{i.e.}, we make precise in which sense the convergence mentioned in the introduction should be understood. 
We remind the reader that we have already mentioned that the stability for martingale representations as given in \cite[Corollary 3.10]{papapantoleon2019stability} will be central to our approach.
Therefore, it is natural to complement the conditions of that corollary with conditions that ensure the convergence of the Lebesgue--Stieltjes integrals associated to the generators of the BSDEs. This is exactly the role of Conditions \ref{BIConditionAk}--\ref{Bintegrator}.

\begin{enumerate}[label=\textup{{(S\arabic*)}},itemindent=0.cm,leftmargin=*]
%
%
{\color{black}
\item\label{MFilqlc} 	
	The filtration $\FilG^{\infty}$ is quasi--left--continuous and the process $X^{\infty,\circ}$ is continuous; the process $X^{\infty,\natural}$ is $\mathbb G^\infty-$quasi--left--continuous, because of its martingale property. 
\item\label{MSIprime}
	The pair $(X^{k,\circ}, X^{k,\natural})\in\mathcal H^{2,k}(\mathbb R^m)\times\mathcal H^{2,d}(\mathbb G^k,\infty;\mathbb R^\ell)$ with $M_{\mu^{X^{k,\natural}}}\big[\Delta X^{k,\circ}\big| \widetilde{\mathcal P}^{\mathbb G^k}\big] = 0$ for every $k\in\overline{\mathbb N}$, such that in addition $X^k = X^{k,\circ} + X^{k,\natural}$, and
	\begin{gather}
		\big(X^{k,\circ}_{\infty}, X^{k,\natural}_{\infty}\big) \LGconvmulti{2}{\ell\times 2} \big(X^{\infty, c}_{\infty},X^{\infty,d}_{\infty}\big).
	\label{MSISeparateFinalRV}
	\end{gather}
\item\label{MXWPRP}
	{The martingale $X^{\infty}$ possesses the ${\FilG^{\infty}}-$predictable representation property.} \phantom{$\xrightarrow{L^2}$}
\item\label{MFilweak}
	The filtrations converge weakly, \textit{i.e.}
	$\FilG^k\weakFil \FilG^{\infty}$. \phantom{$\xrightarrow{L^2}$}
\item\label{Mfinalrv} 	
	The random variable $\xi^k\in \mathbb{L}^2(\Omega,\mathcal{G}^k_{\infty},\mathbb P;\mathbb R)$, for every $k\in\overline{\N}$, and 
			$\xi^k \xrightarrow{\hspace*{0.2cm}\mathbb{L}^2(\Omega,\mathcal{G},\mathbb P;\mathbb R)\hspace*{0.2cm}} \xi^{\infty}.$
}

\item\label{BgeneratorUI} 
The sequence
 	\begin{align*}
 	\bigg(\int_{(0,\tau^k)} \textup{e}^{\hat{\beta} A^k_t} \frac{\big\Vert f^k(t,0,0,\mathbf{0})\big\Vert^2}{(\alpha_t^{k})^2} \textup{d}C_t^k\bigg)_{k\in\overline{\mathbb N}},
 	\end{align*}
 	is uniformly integrable, where $\hat{\beta}$ is the common value from \cref{foot:beta-hat}, while $\mathbf0$ denotes the null application from $\R^n$ to $\R^\ell$.
\smallskip
\item\label{BIConditionAk} 
The sequence $(A^k_\infty)_{k\in\overline{\mathbb N}}$ is bounded by a constant $\overline{A}>0$.\footnote{This allows to write Condition \ref{BgeneratorUI} without the exponential functions.}
\smallskip
\item\label{Bgenerator}	
	The generators $(f^k)_{k\in\overline{\mathbb{N}}}$ possess additionally the following properties:
 \begin{enumerate}[label=\textup{{(\roman*)}},itemindent=0.cm,leftmargin=0.5em]
\medskip
\item \label{BgeneratorPath}
 	for every $k\in\overline{\mathbb{N}}$, $W\in\mathbb{D}(\mathbb{R}_+;\mathbb{R}^\ell)$, $Z\in D^{\circ,\ell\times m}$ and $U\in D^\natural$, it holds that\footnote{In \cref{rem:support-well-posed-tnorm} we verify that we are allowed to use elements of $C_{c|\widetilde{E}}(E;\mathbb{R}^\ell)$ as elements of $\mathfrak{H}^{k}$.
 	Clearly, the elements of $D^{\circ,\ell\times m}$  also possess the suitable measurability properties as deterministic processes which are $\mathcal{B}(\mathbb{R}_+)\backslash \mathcal{B}(\mathbb{R}^{\ell \times m})$-measurable.}%
 	\begin{align*}
 	\big( f^k(t,W_t,Z_t,U(t,\cdot)\big)_{t\in\mathbb{R}_+} \in\mathbb{D}(\mathbb{R}_+;\mathbb{R}^\ell),\; \mathbb{P}\text{\rm--a.s.;}
 	\end{align*}%
 	%
\item\label{Bgenerator:respects-J1-conv}
	for every $Z\in D^{\circ,\ell\times m}$, $U\in D^\natural$, if $(W^{k})_{k\in\overline{\mathbb{N}}}$ is a sequence of c\`adl\`ag maps such that $W^k\xrightarrow{\hspace{0.5em}\textup{J}_1\hspace{0.5em}}W^\infty$
then
\begin{align*}
 \big(f^k(t, W^k_t,Z_t,U(t,\cdot))\big)_{t\in\mathbb{R}_+} 
 \xrightarrow[k\to\infty]{\hspace{1em}\textup{J}_1(\mathbb{R}^\ell)\hspace{1em}}
 \big(f^\infty(t, W^\infty_t,Z_t,U(t,\cdot))\big)_{t\in\mathbb{R}_+},\; \mathbb{P}\text{\rm--a.s.} 
 \end{align*}
Besides, if $\sup_{k\in\overline{\mathbb{N}}}\Vert W^k(\omega)\Vert_\infty <\infty,$ $\mathbb{P}\text{\rm--a.s.}$, then  
$\sup_{k\in\overline{\mathbb{N}}}\big\Vert \big(f^k(t,W_t^k,Z_t,U(t,\cdot))\big)_{t\in\mathbb{R}_+}\big\Vert_\infty<\infty$, $\mathbb{P}\text{\rm--a.s.}$
 \end{enumerate}
%
\medskip
\item\label{Bintegrator} The sequences $(C^k)_{k\in\overline{\mathbb N}}$ and $(\Phi^k)_{k\in\overline{\mathbb N}}$ satisfy  
\begin{enumerate}[label=\textup{{(\roman*)}},itemindent=0.cm,leftmargin=0.5em]
	\item\label{BintegratorJ1P}
		$C^k \xrightarrow{\hspace{0.2em}\left(\textup{J}_1(\mathbb{R}),\mathbb P\right)\hspace{0.2em}} C^\infty$;
	\item\label{BintegratorInfinityProb}
		$C^k_\infty:=\lim_{t\to\infty}C^k_t\in\mathbb{R}_+$, $\mathbb P$--a.s., with $C^k_\infty\xrightarrow{\hspace{0.2em}\mathbb P\hspace{0.2em}}C^\infty_\infty$;
	\item\label{BintegratorPhiBounds} 	
		$\Phi^k\xrightarrow{\hspace{0.2em}|\cdot|\hspace{0.2em}} \Phi^\infty :=0$;
\end{enumerate}
%

\medskip
\item\label{BFiniteStoppingTime} The stopping time $\tau^{\infty}$ is finite and $\tau^{k}\xrightarrow[k\to\infty]{\mathbb{P}}\tau^{\infty}$.
\end{enumerate}

\begin{remark}\label{rem:BSDE-solution}
We include below some brief remarks on the nature of the conditions assumed above. 
A more detailed discussion will be provided in {\rm\Cref{subsec:CommentsonConditions}}.
\begin{enumerate}[itemindent=0.cm,leftmargin=*]
\item  {We have assumed that $(\mathscr{D}^k)_{k\in\overline\N}$ are standard data under $($the common value$)$ $\hat{\beta}$, and that $\hat{\beta}$ is sufficiently large such that the quantities $M^{\Phi^k}(\hat{\beta})$, and $M_{\star}^{\Phi^k}(\hat{\beta})$ as defined in {\rm\cite[Lemma 3.4]{papapantoleon2016existence}} are sufficiently close to $18\textup{e} \Phi^k$.\label{foot:StandardData}
The above in conjunction with \ref{Bintegrator}.\ref{BintegratorPhiBounds} and {\rm\cite[Corollary 3.6]{papapantoleon2016existence}} ensure that for \emph{all but finitely many} $k\in\overline{\mathbb{N}}$ there exists a unique quadruplet 
\begin{align*}
\begin{multlined}[0.8\textwidth]
(Y^k,Z^k,U^k,N^k)\in\mathbb{H}^{2,k}_{\hat{\beta}}\times\mathbb{H}^{2,k,\circ}_{\hat{\beta}}\times\mathbb{H}^{2,k,\natural}_{\hat{\beta}}\times\mathbb{H}^{2,k,\perp}_{\hat{\beta}}\\ 
\text{\rm or}\;
(Y^k,Z^k,U^k,N^k)\in\mathcal{S}^{2,k}_{\hat{\beta}}\times\mathbb{H}^{2,k,\circ}_{\hat{\beta}}\times\mathbb{H}^{2,k,\natural}_{\hat{\beta}}\times\mathbb{H}^{2,k,\perp}_{\hat{\beta}},
\end{multlined}
\end{align*}
such that for any $t\in\llbracket 0,\tau^k\rrbracket$
\begin{align}\label{k-BSDEintegralform}
\begin{multlined}[0.8\textwidth]
Y_t^k = \xi^k + \int_t^{\tau^k} f^k\big(s, Y_s^k,Z_s^k,U^k(s,\cdot)\big)\ud C^k_s \\- \int_t^{\tau^k} Z_s^k \ud X^{k,\circ}_s - \int_t^{\tau^k}\int_{\mathbb{R}^n} U^k(s,x)\widetilde{\mu}^{k,\natural}(\ds,\dx) - \int_t^{\tau^k} \ud N_s^k,\; \mathbb{P}\text{\rm--a.s.}
\end{multlined}
\end{align}
We will assume, without loss of generality, that {\rm BSDE} \eqref{k-BSDEintegralform} admits a unique solution, for every $k\in\overline{\mathbb{N}}$.
The unique solution associated to the standard data $\mathscr{D}^k$ will be denoted by $\mathscr{S}^k$, for every $k\in\overline{\mathbb{N}}$. }

\smallskip
\item The assumption $\Phi^\infty=0$ is compatible with the assumption that $\oX^{\infty}$ is a quasi--left-continuous martingale; recall {\rm Condition \ref{MFilqlc}}.
Indeed, in this case $X^{\infty,\circ}$ $($which has been assumed continuous$;$ see \ref{MFilqlc}$)$ has a continuous angle bracket process, and $\mu^{\infty,\natural}$ has an atomless compensator$;$ see {\rm\citeauthor*{jacod2003limit} \cite[Corollary II.1.19]{jacod2003limit}}. 
\medskip
\item {\rm Condition \ref{MXWPRP}} states that the pair $(X^{\infty,\circ},X^{\infty,\natural})$ possesses the $\mathbb{G}^\infty$--martingale {\color{black} predictable?} representation property. 
Therefore, the element $N\in\mathcal{H}^{2,\infty,\perp}_{\beta}$ of the solution $Y^\infty$ has to be the zero process. 
For the sake of a unified notation, even if $k=\infty$, we will write the orthogonal martingale $N$ in {\rm BSDE} \eqref{k-BSDEintegralform}.
\medskip
\item\label{rem:OnBoundednessA-1} 
{\rm Condition} \ref{BIConditionAk} states that the sequence $(A^k_\infty)_{k\in\overline{\mathbb N}}$ has to be bounded.
On the one hand, this condition implies the equivalence of the norms in the weighted spaces indexed by $\beta$ with the respective non-weighted norms, \emph{i.e.} when $\beta=0$.
On the other hand, {\rm Condition} \ref{BIConditionAk} does not necessarily imply that the sequences $(\alpha^k)_{k\in\overline{\mathbb{N}}}$ and $(C^k)_{k\in\overline{\mathbb N}}$ are bounded.
In \ref{Bintegrator}.\ref{BintegratorInfinityProb} we do assume, however, that $C^k$ is $\mathbb P$--almost surely finite, for every $k\in\overline{\mathbb N}$.
\medskip
\item\label{rem:OnBoundednessA-ii}
In view of {\rm Condition} \ref{BIConditionAk} again, {\rm Condition} \ref{BgeneratorUI} is equivalent to the sequence
\begin{equation*}
 \bigg(\int_{(0,\infty)} \frac{\big\Vert f^k(t,0,0,\mathbf{0})\big\Vert_2^2}{(\alpha_t^{k})^2} \,\textup{d}C_t^k\bigg)_{k\in\overline{\mathbb N}},
\end{equation*}
being uniformly integrable.
Moreover, {\rm Condition \ref{Mfinalrv}} states that $\xi^k\xrightarrow[k\to\infty]{\mathbb{L}^2(\Omega, \mathcal{G},\mathbb{P};\mathbb{R}^{\ell})} \xi^\infty $.
In particular $\big(\Vert \xi^k\Vert_2^2\big)_{k\in\overline{\mathbb{N}}}$ is uniformly integrable, which implies that the sequence $(\xi^k)_{k\in\overline{\mathbb{N}}}$ is $\mathbb{L}^2(\Omega, \mathcal{G},\mathbb{P};\mathbb{R}^{\ell})$-bounded
\begin{gather*}
\sup_{k\in\overline{\mathbb N}}\big\Vert \xi^k\big\Vert_{\mathbb{L}^2(\Omega, \mathcal{G},\mathbb{P};\mathbb{R}^{\ell})}<\infty,\; 
\text{\rm or, equivalently by {\rm Condition} \ref{BIConditionAk},}\; 
\sup_{k\in\overline{\mathbb N}}\; \bigl\Vert \xi^k\big\Vert_{\mathbb{L}^{2,k}_{\hat{\beta}}}<\infty,
\end{gather*}
 where $\hat{\beta}$ is the common value $($see {\rm\cref{foot:beta-hat}}$)$.
\end{enumerate}
\end{remark}

\begin{remark}
\label{rem:support-well-posed-tnorm}
Let $U\in  C_{c|\widetilde{E}}(E;\mathbb{R}^{\ell})$ and $k\in\mathbb{N}$. 
We have that $\textup{supp}(U)\cap E_0=\emptyset$, which in conjunction with the compactness of $\textup{supp}(U)$ allows%
\footnote{We use that every Euclidean space with its usual topology is a normal Hausdorff space.}
us to write $\inf\big\{ \Vert x \Vert_2:x\in\textup{supp}(U) \big\}>0.$ Then, it follows
\begin{align*}
\mathbb{E}\bigg[ \int_{(0,\infty)\times \mathbb{R}^n} \Vert U(t, x)\Vert_2^2 \ \mu^{k,\natural}(\textup{d}t,\textup{d}x) \bigg]
=\mathbb{E}\bigg[ \sum_{t>0} \Vert U(t, \Delta X^{k,\natural}_t)\Vert_2^2 \bigg]
<\infty,
\end{align*}
where we also used the boundedness of the function $U$ and the fact that 
$\mathbb{E}\Big[ \sum_{t>0} \Vert \Delta X^{k,\natural}_t\Vert_2^2 \Big]<\infty;$ recall that $X^{k,\natural}\in \mathcal{H}^{2,d}(\mathbb{G}^{k},\tau^{k};\mathbb{R}^{n})$.
Additionally, by {\rm\citeauthor*{jacod2003limit} \cite[Theorem II.1.8]{jacod2003limit}} we have
\begin{align}\label{prop:sq-int-comp}
\mathbb{E}\bigg[ \int_{(0,\infty)\times\mathbb{R}^n} \Vert U(t, x)\Vert_2^2 \ \nu^{k,\natural}(\textup{d}t,\textup{d}x) \bigg]
<\infty.
\end{align}
Therefore, the stochastic integral $U\star \widetilde{\mu}^{k,\natural}$ is well-defined.
Moreover, by {\rm\citeauthor*{he1992semimartingale} \cite[Theorem 11.21.3]{he1992semimartingale}} and by means of {\rm Property \eqref{prop:sq-int-comp}}, we deduce that
\begin{align*}
\langle U\star \widetilde{\mu}^{k,\natural}\rangle_\cdot  
= U^2*\nu^{k,\natural}_\cdot - \sum_{s\le \cdot}
\int_{\mathbb{R}^n} U(t,x)\nu^{k,\natural}(\{t\}\times \textup{d}x)
\int_{\mathbb{R}^n} U^\top(t,x)\nu^{k,\natural}(\{t\}\times \textup{d}x),\; \mathbb{P}\text{\rm--a.s.}
\end{align*}
The last identity further allows us to write for every $t\in\mathbb{R}_+$
\begin{align}\label{rem:tnorm-U-compact}
\big(\tnorm{U(t,\cdot)}_{k}(\omega)\big)^2 
= \widehat{K}^k_t \big(\Vert U(t,\cdot) \Vert_2^2\big)(\omega) 
- \Delta C^k_t(\omega) \big\Vert \widehat{K}^k_t(U(t,\cdot))(\omega)\big\Vert_2^2,\; \mathbb{P}\text{\rm--a.s.}
\end{align}
The finiteness of $\tnorm{U(t,\cdot)}(\omega)$ is a result of the compact support of $U$ and the $\sigma$-finiteness of the transition kernel $\widehat{K}^{k}_t(\omega)$ $($see {\rm\cite[Proposition II.2.9.(iv)]{jacod2003limit}}$)$.
In particular, observe that $U(t,\cdot)\in \mathfrak{H}^{k}_{\omega,t}$ for $C^{k}\otimes\mathbb{P}\text{\rm--a.e.}$ $(t,\omega)\in\mathbb{R}_+\times \Omega$.
Therefore, seeing $U$ as a deterministic process, we can conclude that $D^{\natural}\subset \mathfrak{H}^{k}$, for every $k\in\overline{\mathbb{N}}$.
Consequently we are allowed to evaluate elements of $D^{\natural}$ in the last argument of the generator $f^k$, for every $k\in\overline{\mathbb{N}}$.\footnote{The reader may recall the domain of the generator; see Condition \ref{Fthree}.}
\end{remark}

%% file: Stability_BSDEJ_Theorem.tex

We start this section with the statement of the main theorem of this article. 
In Subsection \ref{subsec:BodyProof} we outline the strategy we will follow in order to prove Theorem \ref{BSDERobMainTheorem}, which is based on Moore--Osgood's theorem, see \cref{MooreOsgoodA} and the references therein, while the remaining technical parts will be presented in \Cref{subsec:FirstStepInduction} and \Cref{subsec:p-step}. 
%
%
%
%
%
%
%
\begin{theorem}\label{BSDERobMainTheorem}
Let {\rm Conditions \ref{MFilqlc}--\ref{BFiniteStoppingTime}} hold, and for every $k\in\overline{\mathbb N}$, denote by $\mathscr{S}^k:=(Y^k,Z^k,U^k,N^k)$ the unique solution of the {\rm BSDE} \eqref{k-BSDEintegralform} associated to the standard data $\mathscr{D}^k$.
Then
\begin{gather}
\big(Y^k, Z^k\cdot X^{k,\circ} +  U^k\star\widetilde{\mu}^{k,\natural}, N^k\big) 		
\xrightarrow{\hspace{0.2em}\left(\textup{J}_1(\mathbb R^{\ell \times 3}),\mathbb L^2\right)\hspace{0.2em}} 	
\big(Y^\infty, Z^\infty\cdot X^{\infty,\circ} + U^\infty\star \widetilde{\mu}^{\infty,\natural},0 \big),
\label{BSDERobi}\\
\begin{multlined}[c][0.9\displaywidth]
\big([ Y^k],[ Z^k\cdot X^{k,\circ} + U^k\star\widetilde{\mu}^{k,\natural}], [N^k], [ Y^k,X^{k,\circ}],[ Y^k,X^{k,\natural}], [Y^k,N^k] \big)		
\xrightarrow{\hspace{0.2em}\left(\textup{J}_1(\mathbb R^{\ell \times 6\ell}),\mathbb L^1\right)\hspace{0.2em}}\\
\big([ Y^\infty], [ Z^\infty\cdot X^{\infty,\circ} + U^\infty\star\widetilde{\mu}^{\infty,\natural}], 0, [ Y^\infty,X^{\infty,\circ}],  [ Y^\infty,X^{\infty,\natural}], 0\big),
\end{multlined}
\label{BSDERob-square-bracket}
\shortintertext{}
\begin{multlined}[c][0.9\displaywidth]
\big(\langle Y^k\rangle,\langle Z^k\cdot X^{k,\circ} \rangle, \langle U^k\star\widetilde{\mu}^{k,\natural}\rangle, \langle N^k\rangle, \langle Y^k,X^{k,\circ}\rangle,\langle Y^k,X^{k,\natural}\rangle,\langle Y^k,N^k\rangle \big)		
\xrightarrow{\hspace{0.2em}\left(\textup{J}_1(\mathbb R^{\ell \times 7\ell}),\mathbb L^1\right)\hspace{0.2em}}\\
\big(\langle Y^\infty\rangle, \langle Z^\infty\cdot X^{\infty,\circ} \rangle, \langle U^\infty\star\widetilde{\mu}^{\infty,\natural}\rangle, 0,\langle Y^\infty,X^{\infty,\circ}\rangle,  \langle Y^\infty,X^{\infty,\natural}\rangle, 0\big),
\end{multlined}
\label{BSDERobii}
\end{gather}
 where $0$ denotes the zero process whose state space is a finite-dimensional Euclidean space.
\end{theorem}

%% file: Stability_BSDEJ_Proof_Main_Body.tex
%
%
%

The main strategy for the proof can be visualised in \cref{table:Picardscheme}. 
Conditions \ref{MFilqlc}--\ref{BFiniteStoppingTime} ensure the convergence of the standard data; this corresponds to the first column of the table, and the respective convergence is denoted with a solid arrow.
Using \citeauthor*{papapantoleon2016existence} \cite[Corollary 3.15]{papapantoleon2016existence}, we can associate to the standard data $\mathscr{D}^k$, for each $k\in\overline{\mathbb{N}}$, the sequence of Picard iterations $(\mathscr{S}^{k,(p)})_{p\in\mathbb{N}\cup\{0\}}$ (where $\mathscr{S}^{k,(0)}$ is the zero element of the respective product space), which converges to the unique solution $\mathscr{S}^k$; this corresponds to the $k$-row of the table.    
Our aim is to prove the convergence $\mathscr{S}^k\xrightarrow[k\to\infty]{} \mathscr{S}^{\infty}$, which corresponds to the convergence in the last column of the table, and is denoted by a wiggly arrow.

\begin{table}[ht]
\begin{tabular}{c|ccccccc}
$\mathscr{D}^1$ 	&{$\mathscr{S}^{1,(0)}$} 			&{$\mathscr{S}^{1,(1)}$} 			&{$\mathscr{S}^{1,(2)}$}		
 					&{$\cdots$}							&{$\mathscr{S}^{1,(p)}$}			&{$\xrightarrow{\ p\to\infty\ }$} 	&{$\mathscr{S}^{1}$}\\[0.3cm]
$\mathscr{D}^2$ 	&{$\mathscr{S}^{2,(0)}$} 			&{$\mathscr{S}^{2,(1)}$} 			&{$\mathscr{S}^{2,(2)}$}
 					&{$\cdots$}							&{$\mathscr{S}^{2,(p)}$}			&{$\xrightarrow{\ p\to\infty\ }$}	&{$\mathscr{S}^{2}$}\\[0.3cm]
$\mathscr{D}^3$ 	&{$\mathscr{S}^{3,(0)}$} 			&{$\mathscr{S}^{3,(1)}$} 			&{$\mathscr{S}^{3,(2)}$}
					&{$\cdots$}							&{$\mathscr{S}^{3,(p)}$}			&{$\xrightarrow{\ p\to\infty\ }$}	&{$\mathscr{S}^{3}$}\\[0.1cm]
$\vdots$			&{$\vdots$} 						&{$\vdots$} 						&{$\vdots$} 						& 
					&{$\vdots$}							& 												&{$\vdots$} \\[0.2cm]
$\mathscr{D}^k$ 	&{$\mathscr{S}^{k,(0)}$} 			&{$\mathscr{S}^{k,(1)}$} 			&{$\mathscr{S}^{k,(2)}$}			&{$\cdots$}
					&{$\mathscr{S}^{k,(p)}$}			&{$\xrightarrow{\ p\to\infty\ }$}	&{$\mathscr{S}^{k}$}\\
$\big\downarrow$	&\rotatebox[origin=c]{270}{$\dashrightarrow{}$} 	&\rotatebox[origin=c]{270}{$\dashrightarrow{}$}  	&\rotatebox[origin=c]{270}{$\dashrightarrow{}$}  	&  		
					&\rotatebox[origin=c]{270}{$\dashrightarrow{}$}  	&									& \rotatebox[origin=c]{270}{\Large${\rightsquigarrow}$} \\[0.2cm]
$\mathscr{D}^\infty$&{$\mathscr{S}^{\infty,(0)}$} 		&{$\mathscr{S}^{\infty,(1)}$} 		&{$\mathscr{S}^{\infty,(2)}$}		&{$\cdots$}
					&{$\mathscr{S}^{\infty,(p)}$}		&{$\xrightarrow{\ p\to\infty\ }$}	&{$\mathscr{S}^\infty$}\\[0.2cm]
\end{tabular}\\[.3em]
\caption{The doubly-indexed Picard scheme.}
\label{table:Picardscheme}
\end{table}

Strictly speaking, one can generally construct the elements of the doubly-indexed sequence $(\mathscr{S}^{k,(p)})_{k\in\overline{\mathbb{N}}, p\in\mathbb{N}}$ and not the elements of $(\mathscr{S}^{k})_{k\in\overline{\mathbb{N}}}$.
Consequently, in order to achieve our aim, we will apply Moore--Osgood's theorem on $(\mathscr{S}^{k,(p)})_{k\in\overline{\mathbb{N}}, p\in\mathbb{N}}$.
The aforementioned theorem provides a sufficient framework for the existence of the (unconditional) limit of a doubly-indexed sequence. In our case we will obtain 
\begin{align*}
	\lim_{(k,p)\to (\infty,\infty)} \Vert \mathscr{S}^{k,(p)} - \mathscr{S}^{\infty}\Vert =
	\lim_{k\to \infty}\lim_{p\to \infty} \Vert \mathscr{S}^{k,(p)} - \mathscr{S}^{\infty}\Vert =
	\lim_{p\to \infty}\lim_{k\to \infty} \Vert \mathscr{S}^{k,(p)} - \mathscr{S}^{\infty}\Vert =0.
\end{align*}
 Moore--Osgood's theorem requires the uniform convergence in one direction (here, say horizontally) and the pointwise convergence in the other direction (here, say vertically).
\cref{UniformPicardApproximation} guarantees the finally\footnote{	We are going to use the following convention: whenever we write that a sequence $(\alpha^{k,p})_{k,p\in\overline{\mathbb{N}}}$ `converges \emph{finally} uniformly in $k$' (under the metric $\delta$), we mean that there exists $k_0\in\mathbb{N}$ such that
$\sup_{k\ge k_0} \delta(\alpha^{k,p}, \alpha^{k,\infty})\xrightarrow[p\to \infty]{} 0.$}
 uniform in $k$ convergence of the sequence of Picard approximations $\{(\mathscr{S}^{k,(p)})_{p\in\mathbb{N}}\}_{k\in\overline{\mathbb{N}}}$.
In other words, the first condition of Moore--Osgood's theorem is relatively effortlessly satisfied. 
In \cref{table:Picardscheme}, we have denoted these convergences with a solid arrow. 
The second condition of Moore--Osgood's theorem amounts to proving the convergence 
\begin{align}\label{conv:Skp-induction}
	\lim_{k\to\infty} \Vert \mathscr{S}^{k,(p)} - \mathscr{S}^{\infty,(p)}\Vert = 0, \;  \text{for every}\; p\in\mathbb{N}.
\end{align}
In \cref{table:Picardscheme}, we have denoted these convergences with a dashed arrow. 
Naturally, we will prove the required pointwise convergence by means of induction. However, a series of helpful comments will reduce the complexity of the proof of Convergence \eqref{conv:Skp-induction}; 
the details are postponed until \cref{subsec:auxiliaryForInduction}.
The first step of the induction will be proved in \cref{subsec:FirstStepInduction} and the $p$-th step of the induction in \cref{subsec:p-step}.

\subsection{Uniform \emph{a priori} BSDE estimates}\label{subsec:UniformAPrioriEstimates}

The next result provides uniform \emph{a priori} estimates for the tail of the Picard approximations, see \citeauthor*{papapantoleon2016existence} \cite[Corollary 3.15]{papapantoleon2016existence}. 
In particular, \cref{UniformPicardApproximation} generalises \citeauthor*{briand2002robustness} \cite[Corollary 10]{briand2002robustness}.
\begin{proposition}\label{UniformPicardApproximation}
For every $k\in\overline{\mathbb{N}}$, we associate to the standard data $\mathscr{D}^k$ the sequence of Picard iterations $(\mathscr{S}^{k,(p)})_{p\in\mathbb{N}\cup\{0\}}$, where $\mathscr{S}^{k,(0)}$ is the zero element of 
$\mathcal{S}^{2}_{k,\hat{\beta}}\times  \mathbb{H}^{2,k,\circ}_{\hat{\beta}}\times  \mathbb{H}^{2,k,\natural}_{\hat{\beta}}\times  \mathbb{H}^{2,k,\perp}_{\hat{\beta}}$.
There exists 
$k_{\star,0}\in{\mathbb{N}}$ s.t.
\begin{align*}
\lim_{p\to\infty}\sup_{k\geq k_{\star,0}} \big\Vert\mathscr{S}^{k} -\mathscr{S}^{k,(p)}\big\Vert_{\star,k,\hat{\beta}}^2=0.
\end{align*}
In particular, $\sup_{k\ge k_{\star,0}} \Vert \mathscr{S}^{k}\Vert_{\star,k,\hat{\beta}}^{2}<\infty.$
\end{proposition}
\begin{proof} 
We choose $k_{\star,0}$ as the one determined by Lemma \ref{lem:Uniform-Picard-Constamts}.
Essentially, we have constructed contractions associated to the standard data $(\mathscr{D}^{k})_{k\ge k_{\star,0}}$ whose constant is smaller than $1/\sqrt{2}$. 
Then, for any integer $k\ge k_{\star,0}$, we have
\begin{align*}
	\|\mathscr{S}^k -\mathscr{S}^{k,(p)}\|^2_{\star,k,\hat{\beta}} 
		&\leq \sum_{n=0}^{\infty} 2^{n+1} \|\mathscr{S}^{k,(p+n)} - \mathscr{S}^{k,(p+n+1)}\|^2_{\star,k,\hat{\beta}}\leq \sum_{n=0}^{\infty} \frac{2^{n+1}}{ 2^{p+n}} \|\mathscr{S}^{k,(1)}\|^2_{\star,k,\hat{\beta}} 
		= 4^{1-p} \|\mathscr{S}^{k,(1)}\|^2_{\star,k,\hat{\beta}},
		\numberthis\label{UBstarNorm}
\end{align*}
where in the second inequality we have used \cite[Inequality (3.42)]{papapantoleon2016existence} and that in 
Lemma \ref{lem:Uniform-Picard-Constamts} $M_{\star}^{\Phi^k}(\hat{\beta})<\frac{1}{4}$ for $k\ge k_{\star,0}$.
Since $\mathscr{S}^{k,(0)}$ is the zero element, we obtain by \cite[Lemma 3.8]{papapantoleon2016existence}
\begin{equation*}
\big\|\mathscr{S}^{k,(1)}\big\|_{\star,k,\hat{\beta}}^2
	\le \widetilde{\Pi}_\star^{\hat{\beta},\Phi^k}\big\|\xi^k\big\|_{\mathbb{L}^2_{k,\hat{\beta}}}^2  +  \Pi_\star^{\Phi^k}\big(\overline{\gamma}_\star^{\Phi^k}, \hat{\beta}\big) \bigg\|\frac{f^k\bigl(\cdot,0,0,\mathbf 0)}{\alpha^k}\bigg\|_{\mathbb H^2_{k,\hat{\beta}}}^2.
\end{equation*}
Again by \Cref{lem:Uniform-Picard-Constamts}, we derive the uniform bound
\begin{align*}
		\sup_{k\geq k_0^{\star}}\|\mathscr{S}^{k,(1)}\|^2_{\star,k,\hat{\beta}}
&\leq	
		\widetilde{\Pi}^{\hat{\beta},\Phi^{k_{\star,0}}}
		\sup_{k\geq k_0^{\star}} \|\xi^k\|^2_{\mathbb{L}^2_{\hat{\beta}}(\mathcal{G}_{\tau^k}^k)}
	+ 	\Pi_\star^{\Phi^{k}}\big(\overline{\gamma}_\star^{\Phi^{k_{\star,0}}},\hat{\beta}\big)
		\sup_{k\geq k_0^{\star}} \bigg\|\frac{f^k\bigl(\cdot,0,0,\textbf{0})}{\alpha^k}\bigg\|^2_{\mathbb H^2_{\star,k,\hat{\beta}}}<\infty,
		\numberthis\label{UBstarNorm1stPicard}
\end{align*}
which implies the desired result in conjunction with Inequality \eqref{UBstarNorm} and Condition \ref{BgeneratorUI}.

\vspace{0.5em}
For the second statement, we combine Inequalities \eqref{UBstarNorm} and \eqref{UBstarNorm1stPicard} and we obtain
\begin{align*}
\sup_{k\ge k_{\star,0}} \Vert \mathscr{S}^{k}\Vert_{\star,k,\hat{\beta}}^{2}
		\le  2\big(1+\widetilde{\Pi}^{\hat{\beta},\Phi^{k_{\star,0}}}_{\star}\big)
		\sup_{k\geq k_0^{\star}} \|\xi^k\|^2_{\mathbb{L}^{2,k}_{\hat{\beta}}}
	+ 	2\big(1+\Pi_\star^{\Phi^{k}}\big(\overline{\gamma}_\star^{\Phi^{k_{\star,0}}},\hat{\beta}\big)\big)
		\sup_{k\geq k_0^{\star}} \bigg\|\frac{f^k\bigl(\cdot,0,0,\textbf{0})}{\alpha^k}\bigg\|^2_{\mathbb H^{2,k}_{\hat{\beta}}}
	<\infty,
\end{align*}
where $\widetilde{\Pi}^{\hat{\beta},\Phi^{k_{\star,0}}}_{\star}$ is defined in \cite[Lemma 3.8]{papapantoleon2016existence} as follows
\begin{equation*}
	\mathbb{R}_+\times \mathbb{R}_+ \ni (\delta,\Phi) \longmapsto 17+9\textup{e}^{\delta \Phi} \in \mathbb{R}_+.\qedhere
\end{equation*}
\end{proof}

\begin{corollary}\label{cor:UniformConvTerminal}
	Let $M^{k,(p)}_{\infty}:= Z^{k,(p)}\cdot X^{k,\circ}_\infty + U^{k,(p)}\star \widetilde{\mu}^{k,\natural}_\infty + N^{k,(p)}_\infty$, for $k\in\overline{\mathbb{N}}$, $p\in\mathbb{N}\cup\{0\}$ and $M^{k}_{\infty}:= Z^{k}\cdot X^{k,\circ}_\infty + U^{k}\star \widetilde{\mu}^{k,\natural}_\infty + N^{k}_\infty$, for $k\in\overline{\mathbb{N}}$.
	Then, there exists $\bar{k}_{\star,0}\in \mathbb{N}$ such that  
	\begin{align*}
		\lim_{p\to\infty}\sup_{k\ge \bar{k}_{\star,0}} \big\Vert M^{k,(p)}_\infty - M^{k}_{\infty} \big\Vert_{\mathbb{L}^{2}} = 0,
	\end{align*}
	and
	\begin{align*}
		\lim_{p\to\infty}\sup_{k\ge \bar{k}_{\star,0}} \mathbb{E} \Bigg[ \bigg \Vert \textup{Var} \bigg[ \int_{(0,\infty)}
		f^{k}\bigl(s,Y_s^{k,(p)}, Z_s^{k,(p)}, U^{k,(p)}(s,\cdot)\bigr) - 
		f^{k}\bigl(s,Y_s^{k}, Z_s^{k}, U^{k}(s,\cdot)\bigr) \,
		 \textup{d}C_{s}^k\bigg] \bigg \Vert_2^2 \Bigg]= 0.
	\end{align*}
\end{corollary}
\begin{proof}
	The first limit is immediate from \cref{UniformPicardApproximation}, It\^o's isometry, the orthogonality of the respective spaces (see \cite[Corollary 2.7]{papapantoleon2016existence}) and the inequalities
	\begin{align*}
		\big\Vert M^{k,(p)}_{\infty} - M _{\infty}^{k}\big\Vert_{\mathbb{L}^2} \le  
		\big\Vert M^{k,(p)}_{\infty} - M _{\infty}^{k}\big\Vert_{\mathbb{L}^{2,k}_{\hat{\beta}}} \le 
			\big\Vert \mathscr{S}^{k,(p)} - \mathscr{S}^{k}\big\Vert_{\star,k,\hat{\beta}}.
	\end{align*}
	The second limit is again immediate from \cref{UniformPicardApproximation}, Cauchy--Schwarz's inequality (as applied in \cite[Inequality 3.16]{papapantoleon2016existence}), the Lipschitz property of the generator $f^{k}$ for every $k\in\overline{\mathbb{N}}$, the inequalities 
	\begin{align*}
		\frac{r^k}{(\alpha^{k})^2} \le(\alpha^{k})^2,\;  
		\frac{\theta^{k,\circ}}{(\alpha^{k})^{2}},\; \frac{\theta^{k,\natural}}{(\alpha^{k})^{2}}\le 1,\; \text{for every}\;k\in\overline{\mathbb{N}},
	\end{align*}
	and Condition \ref{BIConditionAk}, \emph{i.e.}, the boundedness of $(A^k)_{k\in\overline{\mathbb{N}}}$.
	More precisely, we have
	\begin{align*}
		&\sup_{k\ge \bar{k}_{\star,0}} \mathbb{E} \Bigg[
		\bigg \Vert \textup{Var} 
		\bigg[ \int_{(0,\infty)}
		f^{k}\bigl(s,Y_s^{k,(p)}, Z_s^{k,(p)}, U^{k,(p)}(s,\cdot)\bigr) - 
		f^{k}\bigl(s,Y_s^{k}, Z_s^{k}, U^{k}(s,\cdot)\bigr) \,
		 \textup{d}C_{s}^k
		\bigg] \bigg \Vert_2^2 \Bigg]\\
		&\le 
		\sup_{k\ge \bar{k}_{\star,0}} \mathbb{E} \Bigg[
		\sum_{i=1}^{\ell}\bigg(\int_{(0,\infty)}
		\big| f^{k,i}\bigl(s,Y_s^{k,(p)}, Z_s^{k,(p)}, U^{k,(p)}(s,\cdot)\bigr) - 
		f^{k,i}\bigl(s,Y_s^{k}, Z_s^{k}, U^{k}(s,\cdot)\bigr)\big| \,
		 \textup{d}C_{s}^k\bigg)^2\Bigg]\\
		&\le \frac{1}{\hat{\beta}}
		\sup_{k\ge \bar{k}_{\star,0}} \mathbb{E} \bigg[
		\int_{(0,\infty)} \textup{e}^{\hat{\beta} A^k_s} 
		\frac{\Vert f^{k}\bigl(s,Y_s^{k,(p)}, Z_s^{k,(p)}, U^{k,(p)}(s,\cdot)\bigr) - 
		f^{k}\bigl(s,Y_s^{k}, Z_s^{k}, U^{k}(s,\cdot)\bigr)\Vert^2_2}{(\alpha^k)^2}
		\textup{d}C_{s}^k \bigg]\\
		&\le \frac{1}{\hat{\beta}}
		\sup_{k\ge \bar{k}_{\star,0}} \bigg\{ \overline{A} \Vert Y^{k,(p)} - Y^{k}\Vert^2_{\mathscr{S}^{2,k}_{\hat{\beta}}} 
		+ \Vert (Z^{k,(p)}\cdot X^{k,\circ} + U^{k,(p)}\star \widetilde{\mu}^{k,\natural}) - (Z^{k}\cdot X^{k,\circ} + U^{k}\star \widetilde{\mu}^{k,\natural})\Vert_{\mathcal{H}^{2,k}_{\hat{\beta}}}^{2}\bigg\}\\
		&\le \frac{\overline{A} \vee 1}{\hat{\beta}}
		\sup_{k\ge \bar{k}_{\star,0}} \big\Vert \mathscr{S}^{k,(p)} - \mathscr{S}^{k}\big\Vert_{\star,k,\hat{\beta}} \xrightarrow[p\to\infty]{}0.\qedhere
	\end{align*}
\end{proof}
\subsection{Reducing the complexity of the induction steps and proving Theorem \ref{BSDERobMainTheorem} }\label{subsec:auxiliaryForInduction}

The purpose of the current subsection, is to explain the strategy that will allow us to reduce the complexity of Convergence \eqref{conv:Skp-induction} and to prove Theorem \ref{BSDERobMainTheorem}.

\begin{notation}\label{notation:LS-mpn}
In order to use as compact notation as possible, we introduce the following
\begin{enumerate}[itemindent=0.cm,leftmargin=*]
\item[\textbullet] $\displaystyle L^{k,(p)}_t:= \int_{(0,t]} f^k\bigl(s,Y_s^{k,(p)}, Z_s^{k,(p)}, U^{k,(p)}(s,\cdot)\bigr) \textup{d}C_{s}^k$,  for $t\in[0,\infty]$,\footnote{For $t=\infty$, we have abused notation and we understand the interval $(0,\infty]$ as $(0,\infty)$. Actually, the processes $C^k$, for $k\in\overline{\mathbb{N}}$, are (left) continuous at $t=\infty$ and therefore there is no difference on which integral we consider.} $k\in\overline{\mathbb{N}}$ and $p\in\mathbb{N}\cup\{0\}$.
	For $i\in\{1,\dots,\ell\}$, we denote the $i$-element of $L^{k,(p)}_t$ by $L^{k,(p),i}_t$. 

\medskip
\item[\textbullet] $\displaystyle L^{k}_t:= \int_{(0,t]} f^k\bigl(s,Y_s^{k}, Z_s^{k}, U^{k}(s,\cdot)\bigr) \textup{d}C_{s}^k$,  for $t\in[0,\infty]$, $k\in\overline{\mathbb{N}}$ and $p\in\mathbb{N}\cup\{0\}$.
	For $i\in\{1,\dots,\ell\}$, we denote the $i$-element of $L^{k}_t$ by $L^{k,i}_t$. 

\medskip

	\item[\textbullet] $\displaystyle \Gamma^{k,(p)}:= \int_{(0,\infty)}{}{} \frac{\big\Vert f^k\big(s, Y^{k,(p)}_s, Z^{k,(p)}_s, U^{k,(p)}(s,\cdot)\big)\big\Vert_2^2}{(\alpha^k_s)^2}\textup{d}C^k_s$, for $k\in\overline{\mathbb{N}}$ and $p\in\mathbb{N}\cup\{0\}$.

\medskip
	\item[\textbullet] $\Delta^{k,(p)} := \textup{Tr}\bigl[\langle Z^{k,(p)}\cdot X^{k,\circ}+U^{k,(p)}\star\widetilde{\mu}^{k,\natural}\rangle_{\infty}\big]$,  for $k\in\overline{\mathbb{N}}$ and $p\in\mathbb{N}$.
	%

\end{enumerate}
\end{notation}

After the introduction of these helpful notations, we focus on the aim of the current subsection.
We claim that Convergence \eqref{conv:Skp-induction} is equivalent to proving the validity of the following two convergences
\begin{equation*}
L^{k,(p)}_{\infty}\xrightarrow[k\to\infty]{\hspace{0.2em}\mathbb{L}^2(\Omega,\mathcal{G},\mathbb{P};\mathbb{R}^\ell)\hspace{0.2em}}
L^{\infty,(p)}_{\infty}
\label{LSpInfinity}
\tag{LS$^{(p)}_\infty$},
\end{equation*}
\begin{equation*}
L^{k,(p)}_{\cdot}
\xrightarrow{\hspace{0.2em}\left(\textup{J}_1(\mathbb{R}^{\ell}),\mathbb{L}^2\right)}
L^{\infty,(p)}_{\cdot},
\tag{LS$^{(p)}$}\label{LSp}
\end{equation*}
for every $p\in\mathbb{N}\cup \{0\}$, in conjunction with   
\begin{align*}
\text{the sequence }	\big( \Gamma^{k,(p)})_{k\in\overline{\mathbb{N}}}\text{ is uniformly integrable,} 
\tag{$\textup{UI}^{(p)}$}
\label{UIp}
\end{align*}
for every $p\in\mathbb{N}\cup \{0\}$.
As should be expected, they will be proven by induction.

\vspace{0.5em}
Next, we will assume within the current subsection that the aforementioned convergences are true for every $p\in\mathbb{N}\cup\{0\}$, together with the uniform \textit{a priori} estimates of subsection \ref{subsec:UniformAPrioriEstimates}, and we will prove \cref{BSDERobMainTheorem} under these assumptions.
Then, we will prove in the upcoming subsections the validity of these assumptions, \textit{i.e.} the validity of convergences \eqref{LSpInfinity}, \eqref{LSp} and \eqref{UIp}, for every $p\in\mathbb{N}\cup\{0\}$.

\subsubsection{Convergence \eqref{BSDERobi} is true}
Unsurprisingly, we are going to transform the BSDEs associated to the Picard schemes into martingale representations and then use the stability already proved for the latter; see \cite[Corollary 3.10]{papapantoleon2019stability}. 
Recall that for every $k\in\overline{\mathbb{N}}$, we have stopped the processes indexed by $k$ at time $\tau^k$.
Therefore, we can substitute for every $k\in\overline{\mathbb{N}}$ the terminal time $\tau^k$ by $\infty$ and we will do so for notational convenience.
The reader may also keep in mind that $N^{\infty,(p)}=0$, for every $p\in\mathbb{N}$.

\vspace{0.5em}
For every $k\in\overline{\mathbb{N}}$ and every $p\in\mathbb{N}\cup\{0\}$ it is true (by construction) that
\begin{align}
Y^{k,(p+1)}_t 	
\nonumber	&= \xi^k + \int_{(t,\tau^k]} f^k\bigl(s,Y_s^{k,(p)}, Z_s^{k,(p)}, U^{k,(p)}(s,\cdot)\bigr)\textup{d}C_{s}^k \\
\nonumber	&\quad	-\int_t^{\tau^k}  Z_s^{k,(p+1)} \textup{d} X^{k,\circ}_s
		-\int_t^{\tau^k}\int_{\mathbb{R}^{n}} U^{k,(p+1)}(s,x)\, \widetilde{\mu}^{k,\natural}(\textup{d}s,\textup{d}x) - \int_t^{\tau^k} \textup{d}N^{k,(p+1)}_s\\
\begin{split}
	&= \xi^k + \int_{(t,\infty)} f^k\bigl(s,Y_s^{k,(p)}, Z_s^{k,(p)}, U^{k,(p)}(s,\cdot)\bigr)\textup{d}C_{s}^k \\
	&\quad	-\int_t^{\infty}  Z_s^{k,(p+1)} \textup{d} X^{k,\circ}_s
		-\int_t^{\infty}\int_{\mathbb{R}^{n}} U^{k,(p+1)}(s,x) \widetilde{\mu}^{k,\natural}(\textup{d}s,\textup{d}x) - \int_t^{\infty} \textup{d}N^{k,(p+1)}_s .
	\label{PicardBSDEorthogonal}
\end{split}
\end{align}
By \Cref{UniformPicardApproximation}, we get that 
$\mathscr{S}^{k,(p)} \xrightarrow[p\to\infty]{\hspace{0.2em}\|\cdot\|_{\star,k,\hat{\beta}}\hspace{0.2em}} \mathscr{S}^{k}$
finally uniformly in $k$. 
In particular, this convergence implies
$Y^{k,(p)} \underset{p\to\infty}{\xrightarrow{\hspace{0.2em}\mathcal{S}_{k,\hat{\beta}}^2\hspace{0.2em}}} Y^{k}$ finally  uniformly in $k$, 
which in turn implies   
$Y^{k,(p)} \xrightarrow[p\to\infty]{\hspace{0.2em}\left(\textup{J}_1(\mathbb{R}^{\ell}),\mathbb{L}^2\right)\hspace{0.2em}} Y^{k}$ finally  uniformly in $k$. 
Hence, by orthogonality of the respective parts, It\^o's isometry and Doob's inequality, see \cite[Section 3.5]{papapantoleon2016existence}, we obtain
\begin{align*}
\big(Y^{k,(p)}, Z^{k,(p)} \cdot X^{k,\circ},  U^{k,(p)}\star \widetilde{\mu}^{k,\natural},	N^{k,(p)} \big) 
\xrightarrow[p\to\infty]{\hspace*{0.2em}(\|\cdot\|_{\infty}, \mathbb{L}^{2})\hspace*{0.2em}}
\big(Y^{k},  Z^{k}\cdot X^{k,\circ}, U^{k}\star \widetilde{\mu}^{k,\natural}, N^k\big),
\end{align*}
finally uniformly in $k.$
At this point, we will combine two facts in order to rewrite the above convergence under the $\textup{J}_1$-topology.
The first one is that the convergence under the $\|\cdot\|_\infty$-norm allows us to conclude the convergence of the sum of two convergent sequences.
The second is that every $\|\cdot\|_\infty$-convergent sequence is also $\textup{J}_1$-convergent; see \cite[Proposition VI.1.17]{jacod2003limit}.
The latter argument was also used a few lines above.
Therefore
\begin{align}\label{PicardUnifOnk-J1}
\big(Y^{k,(p)}, Z^{k,(p)} \cdot X^{k,\circ} +  U^{k,(p)}\star \widetilde{\mu}^{k,\natural}, N^{k,(p)} \big) 
\xrightarrow[p\to\infty]{\hspace*{0.2em}(\textup{J}_1(\mathbb{R}^{\ell\times 3}), \mathbb{L}^{2})\hspace*{0.2em}}
\big(Y^k,  Z^k\cdot X^{k,\circ}+ U^k\star \widetilde{\mu}^{k,\natural}, N^k\big),
\end{align}
finally uniformly in $k.$
Consequently, in order to apply Moore--Osgood's theorem, see \cref{MooreOsgoodA}, it is sufficient to prove that for every $p\in \mathbb{N}\cup\{0\}$, we have
\begin{equation}\label{PicardFixedp}		
\big(	Y^{k,(p)},  Z^{k,(p)}\cdot X^{k,\circ} + U^{k,(p)}\star \widetilde{\mu}^{k,\natural}, 	N^{k,(p)}\big) 
\xrightarrow[k\to\infty]{\hspace{0.2em}(\textup{J}_1(\mathbb{R}^{\ell\times 3}),\mathbb{L}^2)\hspace{0.2em}}
\big(	Y^{\infty,(p)}, Z^{\infty,(p)}\cdot X^{\infty,\circ} +	U^{\infty,(p)}\star \widetilde{\mu}^{\infty,\natural}, N^{\infty,(p)} \big).
\end{equation}
To this end, let us relate BSDE \eqref{PicardBSDEorthogonal} to appropriate martingales.\footnote{The same technique was used in the proof of \cite[Theorem 3.5]{papapantoleon2016existence} in order to use finally the orthogonal decomposition of square-integrable martingales.}
The aforementioned transformation will allow us to use the stability of martingale representations; see \cite[Corollary 3.10]{papapantoleon2019stability}.
For fixed $k\in\overline{\mathbb{N}}$, $p\in\mathbb{N}\cup\{0\}$ we define for $t\in[0,\infty]$
\begin{align*}
	M^{k,(p)}_t\ 	&\overset{\phantom{\eqref{PicardBSDEorthogonal}}}{:=}
						\ Y^{k,(p+1)}_t + \int_{(0,t]} f^k\bigl(s,Y_s^{k,(p)}, Z_s^{k,(p)}, U^{k,(p)}(s,\cdot)\bigr)\textup{d}C_{s}^k \numberthis\label{def:Mk}\\
					&\overset{\eqref{PicardBSDEorthogonal}}{=}  \xi^k + \int_{(0,\infty)} f^k\bigl(s,Y_s^{k,(p)}, Z_s^{k,(p)}, U^{k,(p)}(s,\cdot)\bigr)\textup{d}C_{s}^k \\
					&\qquad -\int_t^\infty  Z_s^{k,(p+1)} \textup{d} X^{k,\circ}_s
					  	-\int_t^\infty\int_{\mathbb{R}^{n}} U^{k,(p+1)}(s,x)\, \widetilde{\mu}^{k,\natural}(\textup{d}s,\textup{d}x) 
					  	- \int_t^\infty \textup{d}N^{k,(p+1)}_s\\
					&\overset{\eqref{PicardBSDEorthogonal}}{=} Y^{k,(p+1)}_0 +\int_0^\infty  Z_s^{k,(p+1)} \textup{d} X^{k,\circ}_s
						+\int_0^\infty\int_{\mathbb{R}^{n}} U^{k,(p+1)}(s,x)\, \widetilde{\mu}^{k,\natural}(\textup{d}s,\textup{d}x) 
						+ \int_0^\infty \textup{d}N^{k,(p+1)}_s\\
					&\qquad -\int_t^\infty  Z_s^{k,(p+1)} \textup{d} X^{k,\circ}_s
						-\int_t^\infty\int_{\mathbb{R}^{n}} U^{k,(p+1)}(s,x)\, \widetilde{\mu}^{k,\natural}(\textup{d}s,\textup{d}x) 
						- \int_t^\infty \textup{d}N^{k,(p+1)}_s\\
					&\overset{\phantom{\eqref{PicardBSDEorthogonal}}}{=} Y^{k,(p+1)}_0 + \int_0^t  Z_s^{k,(p+1)} \textup{d} X^{k,\circ}_s + 
						\int_0^t\int_{\mathbb{R}^{n}} U^{k,(p+1)}(s,x)\, \widetilde{\mu}^{k,\natural}(\textup{d}s,\textup{d}x) 
						+ \int_0^t \textup{d}N^{k,(p+1)}_s.
	\numberthis\label{Mkpmartingale} 
\end{align*}
Hence $M^{k,(p)}\in\mathcal{H}^2(\mathbb{G}^{k},\tau^{k};\mathbb{R}^{\ell})$ for every $k\in\overline{\mathbb{N}}$ and every $p\in\mathbb{N}\cup\{0\}$.
At this point, for fixed $p\in\mathbb{N}\cup\{0\}$ we can obtain the convergence
\begin{equation}\label{conv:Picardfixedp-partial}
\big( Z^{k,(p+1)}\cdot X^{k,\circ} + U^{k,(p+1)}\star \widetilde{\mu}^{k,\natural}, 	N^{k,(p+1)}\big) 
\xrightarrow[k\to\infty]{\hspace{0.2em}\left(\textup{J}_1(\mathbb{R}^{\ell\times 2}),\mathbb{L}^2\right)\hspace{0.2em}}
\big( Z^{\infty,(p+1)}\cdot X^{\infty,\circ} +U^{\infty,(p+1)}\star \widetilde{\mu}^{\infty,\natural}, 0 \big),
\end{equation}
if we apply \cite[Corollary 3.10]{papapantoleon2019stability} to the sequence $(M^{k,(p)})_{k\in\overline{\mathbb{N}}}$.
In view of Conditions \ref{MFilqlc}, \ref{MSIprime}, \ref{MXWPRP} and \ref{MFilweak},  we need only to prove the convergence
\begin{equation}\label{conv:TerminalM}
M^{k,(p)}_\infty
\xrightarrow[k\to\infty]{\hspace{0.2em}\mathbb{L}^2(\Omega,\mathcal{G},\mathbb{P};\mathbb{R}^\ell)\hspace{0.2em}}
M^{\infty,(p)}_\infty,
\end{equation}
for every $p\in\mathbb{N}\cup\{0\}$, in order to apply \cite[Corollary 3.10]{papapantoleon2019stability}.

\medskip
However, in view of Condition \ref{Mfinalrv}, which states that the sequence $(\xi^k)_{k\in\overline{\mathbb{N}}}$ is $\mathbb{L}^2$-convergent, and recalling that for the Picard schemes holds $Y^{k,(p+1)}_\infty=\xi^k,$ for every $k\in\overline{\mathbb{N}}$ and $p\in\mathbb{N}\cup\{0\},$
we immediately obtain from Identity \eqref{def:Mk} that (for the same $p$) 
\begin{equation*}
L^{k,(p)}_{\infty}
\xrightarrow[k\to\infty]{\hspace{0.2em}\mathbb{L}^2(\Omega,\mathcal{G},\mathbb{P};\mathbb{R}^\ell)\hspace{0.2em}}
L^{\infty,(p)}_{\infty},
\tag{\ref{LSpInfinity}}
\end{equation*}
is equivalent to Convergence \eqref{conv:TerminalM}.
\medskip
Assume for the following that \eqref{LSpInfinity} is valid for every $p\in\mathbb{N}\cup\{0\},$ \emph{i.e.} we can apply \cite[Corollary 3.10]{papapantoleon2019stability} for the martingale sequence $(M^{k,(p)})_{k\in\overline{\mathbb{N}}}$.
Then, we obtain the convergence
\begin{equation}\label{conv:MZXUXN}
\begin{multlined}[c][0.87\displaywidth]
\big(M^{k,(p)}, Z^{k,(p+1)}\cdot X^{k,\circ} + U^{k,(p+1)}\star \widetilde{\mu}^{k,\natural}, N^{k,(p+1)}\big)
\xrightarrow{\hspace{0.2em}\left(\textup{J}_1(\mathbb{R}^{\ell\times 3}),\mathbb{L}^2\right)}\\
\big(M^{\infty,(p)}, Z^{\infty,(p+1)}\cdot X^{\infty,\circ} + U^{\infty,(p+1)}\star \widetilde{\mu}^\infty, 0\big),
\end{multlined}
\end{equation}
 for every $p\in\mathbb{N}\cup\{0\}$.
Comparing now Convergence \eqref{conv:MZXUXN} with Convergence \eqref{PicardFixedp} associated to the $(p+1)$--Picard step, we realise that they differ only in the first element.
Recall the definition of $M^{k,(p)}$, for $k\in\overline{\mathbb{N}}$; see \eqref{def:Mk}.
It is immediate\footnote{Since we are using the $\textup{J}_1$-topology, we have to be careful with arguments like this.
However, the  continuity of the process corresponding to the Lebesgue--Stieltjes integrals with respect to $C^\infty$ (recall that $\Phi^\infty=0$) allows us to proceed.} that we can obtain the convergence 
\begin{equation}\label{conv:MY}
\big(M^{k,(p)}, Y^{k,(p+1)}\big)
\xrightarrow{\hspace{0.2em}\left(\textup{J}_1(\mathbb{R}^2),\mathbb{L}^2\right)}
\big(M^{\infty,(p)}, Y^{\infty,(p+1)}\big), \text{ for every }p\in\mathbb{N}\cup\{0\},
\end{equation}
if the convergence 
\begin{equation*}
\begin{multlined}[c][0.87\displaywidth]
L^{k,(p)}_{\cdot}
\xrightarrow{\hspace{0.2em}\left(\textup{J}_1(\mathbb{R}^{\ell}),\mathbb{L}^2\right)},
L^{\infty,(p)}_{\cdot}
\end{multlined}
\tag{\ref{LSp}}
\end{equation*}
holds for every $p\in\mathbb{N}\cup\{0\}$.
Hence, our claim as stated in the title of the current section is valid.
\begin{corollary}
	If the convergences \eqref{LSpInfinity} and \eqref{LSp} are true for every $p\in\mathbb{N}\cup \{0\}$, then 
	\begin{align}\label{cor:MartTermValMooreOsgood}
		M^{k,(p)}_\infty
\xrightarrow[(k,p)\to(\infty,\infty)]{\hspace{0.2em}\mathbb{L}^2(\Omega,\mathcal{G},\mathbb{P};\mathbb{R}^\ell)\hspace{0.2em}}
		 M^{\infty}_{\infty},
	\end{align}
	\begin{align}\label{cor:LSIntegTermValMooreOsgood}
		L^{k,(p)}_{\infty}
\xrightarrow[(k,p)\to(\infty,\infty)]{\hspace{0.2em}\mathbb{L}^2(\Omega,\mathcal{G},\mathbb{P};\mathbb{R}^\ell)\hspace{0.2em}}
		L^{\infty}_{\infty}.
	\end{align}
	In both cases, the iterated limits exist and are equal to the respective right-hand side.

\end{corollary}

\begin{proof}
	In view of \cref{cor:UniformConvTerminal} and \Cref{conv:TerminalM}, Moore--Osgood's theorem (see \cref{MooreOsgoodA}) ensures the convergence in \eqref{cor:MartTermValMooreOsgood}. 
	For the convergence in \eqref{cor:LSIntegTermValMooreOsgood}, one uses \eqref{def:Mk} and \eqref{cor:MartTermValMooreOsgood} in conjunction with the convergence 
	\begin{align*}
		\xi^k \xrightarrow[k\to\infty]{\hspace{0.2em}\mathbb{L}^2(\mathcal{G};\mathbb{R}^\ell)\hspace{0.2em}} \xi^\infty;
	\end{align*}
	see Condition \ref{Mfinalrv}.
\end{proof}%
\subsubsection{Convergence \eqref{BSDERob-square-bracket} is true}
The reader should recall that the convergences \eqref{LSpInfinity} and \eqref{LSp} are assumed true for every $p\in\mathbb{N}\cup \{0\}$.
We will apply Moore--Osgood's theorem for the doubly-indexed sequence $([\mathscr{S}^{k,(p)}])_{k\in\overline{\mathbb{N}},p\in\mathbb{N}\cup\{0\}}$.
It is well understood that the convergence of a sequence of (special) semimartingales does not guarantee the convergence of the associated square brackets. 
However, if the sequence is P--UT (see \citeauthor{jacod2003limit} \cite[Definition VI.6.1]{jacod2003limit}), then \cite[Theorem VI.6.26]{jacod2003limit} ensures that we have the desired convergence.
To this end, we will prove that the P--UT property holds in a proper sense; the details will be provided below, but one could describe it as `the sequence $(\mathscr{S}^{k,(p)})_{k\in\overline{\mathbb{N}},p\in\mathbb{N}\cup \{0\}}$ is P--UT finally uniformly in $k$'.

%
%
%


\vspace{0.5em}
In this paragraph we focus initially on the sequences corresponding to the martingales of the Picard schemes, where we have also adjoined the sequence $(M^{k,(p)})_{k\in\overline{\mathbb{N}}, p\in\mathbb{N}\cup \{0\} }$.
In this case, an integrability condition, namely \cite[Corollary VI.6.30]{jacod2003limit}, is sufficient for the P--UT property.\footnote{The reader can immediately verify from Convergence \eqref{cor:MartTermValMooreOsgood} that the martingale sequences are $\mathbb{L}^2$-bounded, a property that implies Condition 6.31 of \cite[Corollary 6.30]{jacod2003limit}.\label{foot:PUT-Martingales}}
Moreover, from Convergence \eqref{cor:MartTermValMooreOsgood} we have, in particular, that the sequence $(\Vert M^{k,(p)}_\infty \Vert_2^2)_{k\in\overline{\mathbb{N}}, p\in\mathbb{N}\cup \{0\}}$ is uniformly integrable.
Therefore, in view of \cite[Corollary VI.6.30]{jacod2003limit}, we can conclude that the sequence of martingales is P-UT and the (joint) convergence of the square brackets is obtained.

\vspace{0.5em}
In this paragraph we turn our attention to the sequence $(Y^{k,(p)})_{k\in\overline{\mathbb{N}}, p\in\mathbb{N}\cup \{0\}}$.
The reader should recall that $Y^{k,(p)} = M^{k,(p)} + L^{k,(p)}$ for every $k\in\overline{\mathbb{N}},$ $p\in\mathbb{N}\cup \{0\}$.
We have already argued about the P--UT property of the martingale sequence $(M^{k,(p)})_{k\in\overline{\mathbb{N}}, p\in\mathbb{N}\cup \{0\} }$ and the convergence of the associated square bracket sequence. 
The continuity of $L^{\infty}$, hence the continuity of $[L^{\infty}]=0$, allows us to simply derive the convergence of $([L^{k,(p)}])_{k\in\overline{\mathbb{N}},p\in\mathbb{N}\cup \{0\}}$ to the zero process. 
Then, we can have the joint convergence of the square bracket sequences, and the polarisation identity will allow us to conclude the desired Convergence \eqref{BSDERob-square-bracket}. 
We leave these details to the reader.
Hence, it is left to prove the convergence 
\begin{align*}
	[L^{k,(p)}]\xrightarrow[(k,p)\to(\infty,\infty)]{\hspace{0.2em}(\textup{J}_1(\mathbb{R}^{\ell\times\ell}), \mathbb{L}^{1})\hspace{0.2em} } 0.
\end{align*}
We will use Moore--Osgood's theorem once again, this time for the sequence $([L^{k,(p)}])_{k\in\overline{\mathbb{N}},p\in\mathbb{N}\cup \{0\}}$. On the one hand, from \cite[Comment VI.6.6]{jacod2003limit}, the sequence $(L^{k,(p)})_{p\in\mathbb{N}\cup \{0\}}$ is P--UT, for every $k\in\overline{\mathbb{N}}$, if the sequence
\begin{align*}
\Big( \textup{Tr}\Big[\textup{Var}\big[ L^{k,(p)}_\cdot\big]_t \Big]\Big)_{p\in\mathbb{N}\cup \{0\}}, 
\end{align*}
is tight in $\mathbb{R}$, for every $t\in\mathbb{R}_+$, where the total variation is calculated element-wise.
Using \cref{cor:UniformConvTerminal} and \cref{UniformPicardApproximation}, we have that for the arbitrary $\varepsilon>0$, there exists $p_0\in\mathbb{N}$, which depends only on $\varepsilon$ and not on $k$, such that
\begin{align*}
 \sup_{p\ge p_0} \sup_{k\ge k_{\star,0}}\mathbb{E} \Big[\big \Vert \textup{Var}\big[ L^{k,(p)}_\cdot\big]_t \big \Vert_2^{2}\Big]
 &\le 2\sup_{p\ge p_0} \sup_{k\ge k_{\star,0}}\mathbb{E} \Big[\big \Vert\textup{Var}\big[ L^{k,(p)}_\cdot - L^{k}\big]_t \big \Vert_2^{2}\Big] + 2
 \sup_{k\ge k_{\star,0}}\mathbb{E} \Big[\big \Vert\textup{Var}\big[ L^{k}_\cdot\big]_t \big \Vert_2^{2}\Big]\\
 &\le 2\varepsilon +
 2\sup_{k\ge k_{\star,0}}\mathbb{E} \Big[\big \Vert\textup{Var}\big[ L^{k}_\cdot\big]_t \big \Vert_2^{2}\Big]\\
 &\le \varepsilon + \sup_{k\ge k_{\star,0}} \Vert \mathscr{S}^{k} \Vert_{\star,k,\hat{\beta}} +
 \sup_{k\ge k_{\star,0}}\Big\Vert \frac{f^{k}(\cdot,0,0,\mathbf{0})}{\alpha^{k}} \Big\Vert_{\mathbb{H}^{2,k}_{\hat{\beta}}}  <\infty.
\end{align*}
For the derivation of the upper bound of the total variation of $L^{k}$ we have used arguments analogous to \cref{cor:UniformConvTerminal}.  
Now, Markov's inequality implies the boundedness in probability (uniformly in $k$) of the desired sequence. 
Hence, for every $k\ge k_{\star,0}$, the sequence $(L^{k,(p)})_{p\in\mathbb{N}\cup \{0\}}$ is P--UT and consequently
\begin{align*}
	[L^{k,(p)}]\xrightarrow[p\to\infty]{\hspace{0.2em}(\Vert\cdot \Vert_{\infty},\mathbb{L}^{1})\hspace{0.2em} } [L^{k}],\text{ for every $k\ge k_{\star,0}$},
\end{align*}
by means of \cite[Theorem VI.6.26]{jacod2003limit}.
We can derive the convergence of the square brackets finally uniform in $k$ from the following inequality
\begin{align*}
 	\textup{Tr}\big[[L^{k,(p)} - L^{k}]_\infty\big] 
 	= \sum_{i=1}^{\ell} \sum_{t\ge 0} \big( \Delta (L^{k,(p),i} - L^{k,i})_t\big)^2
 	&\overset{\footnotemark}{\le} \sum_{i=1}^{\ell} \bigg(\sum_{t\ge 0} \big| \Delta (L^{k,(p),i} - L^{k,i})_t\big| \bigg)^2 \\
	&\le \big \Vert \textup{Var}\big[ L^{k,(p)} - L^{k}\big]\big\Vert_2^{2},
\end{align*}%
\footnotetext{For the sequence spaces $\ell^1(\mathbb{N})$ and $\ell^2(\mathbb{N})$ it is true that $\ell^{1}(\mathbb{N})\subset \ell^2(\mathbb{N})$, 
\emph{i.e.}, for a sequence of real numbers $x:=(x_p)_{p\in\mathbb{N}}$ it holds 
$\Vert x\Vert_{\ell^2(\mathbb{N})} := (\sum_{p\in\mathbb{N}} |x_p|^{2})^{\frac{1}{2}} \le \sum_{p\in\mathbb{N}} |x_p| = \Vert x \Vert_{\ell^1(\mathbb{N})}$.}which is true for every $k$. On the other hand, the induction hypothesis for every $p\in\mathbb{N}\cup \{0\}$ allows us to apply \cite[Corollary 3.10]{papapantoleon2019stability}, which in particular provides
\begin{align*}
	[L^{k,(p)}] \xrightarrow[k\to\infty]{\hspace{0.2em}(\textup{J}_1(\mathbb{R}^{\ell\times\ell}), \mathbb{L}^{1})\hspace{0.2em} }  [L^{\infty,(p)}], \text{ for every } p\in\mathbb{N}\cup \{0\}.
\end{align*}
Overall, the conditions of Moore--Osgood's theorem are satisfied, hence we have 
\begin{align*}
	[L^{k,(p)}] \xrightarrow[(k,p)\to(\infty,\infty)]{\hspace{0.2em}(\textup{J}_1(\mathbb{R}^{\ell\times\ell}), \mathbb{L}^{1})\hspace{0.2em} }  [L^{\infty}]=0.
\end{align*}
Both iterated limits of the sequence $([L^{k,(p)}])_{k\in\overline{\mathbb{N}}, p \in \mathbb{N}\cup \{0\}}$ exist and are equal to $0$.
The continuity of the limit allows us to derive the (joint) convergence \eqref{BSDERob-square-bracket}.

\subsubsection{Convergence \eqref{BSDERobii} is true}
The convergence of its martingale parts can be justified because of \cite[Corollary 3.10]{papapantoleon2019stability}, whose validity implied the convergence of the martingale parts of \eqref{BSDERobi}.
For the parts associated to the square-integrable $\mathbb{G}^{k}$--special semimartingales $(Y^{k})_{k\in\overline{\mathbb{N}}}$, we observe that for every $k\in\overline{\mathbb{N}}$ the process $[Y^{k}]$ is a special $\mathbb G^k$-semimartingale with canonical decomposition $[Y^{k}]=([Y^{k}]-\langle Y^{k}\rangle) + \langle Y^{k}\rangle$. 
Moreover, the arguments presented above allow us to conclude that the sequence $\big(\text{Tr}\big[[Y^{k}]_\infty\big]\big)_{k\in\overline{\mathbb{N}}}$ is uniformly integrable.
This further implies\footnote{See \citet[Th\'eor\`eme 3.2.1)]{lenglart1980presentation} and use the de La Vall\`ee-Poussin criterion} the uniform integrability of $\big(\text{Tr}\big[\langle Y^{k}\rangle_\infty\big]\big)_{k\in\overline{\mathbb{N}}}$. 
Hence, the sequence associated to the angle brackets is also tight in $\mathbb{R}^\ell$. 
Therefore, we can apply\footnote{Along with standard arguments and the polarisation identity since the aforementioned theorem is stated in the real-valued case.} \citeauthor{memin2003stability} \cite[Theorem 11]{memin2003stability} for the sequence $([Y^k])_{k\in\overline{\mathbb{N}}}$.

%% file: Stability_BSDEJ_Induction_FirstStep.tex

Recall that in the proof of \cref{UniformPicardApproximation} we have set $\Upsilon^{k,(0)}=(0,0,\mathbf{0})$, for every $k\in\overline{\mathbb{N}}$. 
Now we provide some useful lemmata that we will then use for proving the first step of the induction in \cref{prop:first-step-Induction}.

\begin{lemma}\label{lem:Ck-PUT}
The sequences $(C^k)_{k\in\overline{\mathbb{N}}}$ and $( L^{k,(0)})_{k\in\overline{\mathbb{N}}}$ possess the \textup{P--UT} property.
\end{lemma}
\begin{proof}
First , recall the fact that $(C^k)_{k\in\overline{\mathbb{N}}}$ is a sequence of increasing processes.
Therefore, $\textup{Var}(C^k)=C^k$, for every $k\in\overline{\mathbb{N}}.$
Secondly, Condition \ref{Bintegrator}.\ref{BintegratorJ1P} implies that $(C^k)_{k\in\overline{\mathbb{N}}}$ is tight in $\mathbb D(\mathbb{R}_+;\mathbb{R})$, which in turn implies that $\textup{Var}(C^k)_t$ is tight in $\mathbb{R}$ for every $t\in\mathbb{R}_+.$
Now, we can conclude by \cite[Proposition VI.6.12]{jacod2003limit}. 

\vspace{0.5em}
For the sequence $(L^{k,(0)})_{k\in\overline{\mathbb{N}}}$ it is sufficient to prove that $\big(\textup{Var}\big[ L^{k,(0)}\big]_\infty\big)_{k\in\overline{\mathbb{N}}}$ is tight in $\mathbb{R}^{\ell}$; see \cite[Remark VI.6.6]{jacod2003limit}.
To this end, we will prove that the sequence is $\mathbb{L}^2$-bounded.
Indeed, by the following inequality $\big($here $f^{k,i}$ denotes the $i$-element of $f^{k}$, for every $k\in\overline{N}$, $i\in\{1,\dots,\ell\}\big)$
\begin{align*}
\textup{Var}[L^{k,(0),i}]_\infty \le \int_{(0,\infty)}  \bigl\vert f^{k,i}\bigl(s,Y_s^{k,(p)}, Z_s^{k,(p)}, U^{k,(p)}(s,\cdot)\bigr) \bigr\vert \d C_{s}^k,
\end{align*}
for every $i\in\{1,\dots,\ell\}$
and by Cauchy--Schwarz's inequality, applied as in \cite[Inequality (3.16)]{papapantoleon2016existence}, we derive
\begin{align*}
\sup_{k\in\overline{\mathbb{N}}}
\mathbb{E}\Big[\big\Vert\textup{Var}\big[ L^{k,(0)}\big]_\infty\big\Vert_2^2\Big]
	&\le 
	\frac{1}{\hat{\beta}}\sup_{k\in\overline{\mathbb{N}}}
	\mathbb{E} \big[ \Gamma^{k,(0)}\big]
	\overset{\text{\ref{BgeneratorUI}}}{<}\infty. \qedhere
\end{align*}
\end{proof}
In the remainder of this article,  we denote by $\mu_{C^{i}}$ the (random) measure on $(\mathbb{R}_+,\mathcal{B}(\mathbb{R}_+))$ associated to the increasing and c\`adl\`ag process $C^{i}$, for $i\in\overline{\mathbb{N}}$. 
\begin{lemma}\label{lem:Conv-Ck}
For any subsequence $(C^{k_l})_{l\in\mathbb{N}}$ there exists a further subsequence $(C^{k_{l_i}})_{i\in\mathbb{N}}$ such that 
\begin{equation*}
C^{k_{l_i}}\xrightarrow[i\to\infty]{\hspace{0.2em}\textup{J}_1(\mathbb{R})\hspace{0.2em}}C^\infty,\; \mathbb P\text{\rm --a.s.,} 
\text{ \rm as well as }
C^{k_{l_i}}_\infty\xrightarrow[i\to\infty]{\hspace{0.2em}| \cdot |\hspace{0.2em}} C^\infty_\infty, \; \mathbb P\text{\rm--a.s.}
\end{equation*}
Moreover, $\mu_{C^{k_{l_i}}}\xrightarrow[i\to\infty]{\hspace{0.2em}\textup{w}\hspace{0.2em}} \mu_{C^\infty},$ $\mathbb P${\rm--a.s.}, and $\sup_{i\in\mathbb{N}} \mu_{C^{k_{l_i}}}(\mathbb{R}_+)<\infty$, $\mathbb{P}${\rm--a.s.}
\end{lemma}
\begin{proof}
The first statement is direct by \citeauthor{dudley2002real} \cite[Theorem 9.2.1]{dudley2002real}.
Indeed, the fact that $(\mathbb D(\mathbb{R}_+;\mathbb{R}),\delta_{\textup{J}_1})$\footnote{The metric compatible with the $\text{J}_1$-convergence will be denoted by $\delta_{\text{J}_1}$ and the metric associated to a norm $\Vert\cdot\Vert$ will be denoted by $\delta_{\Vert\cdot\Vert}$.} and $(\mathbb{R},\delta_{|\cdot|})$ are both Polish spaces, together with Conditions \ref{BintegratorJ1P}, \ref{BintegratorInfinityProb} of \ref{Bintegrator} allow us to verify the statement.
Passing possibly to a further subsequence we can assume without loss of generality that both convergent sequences are indexed by $(k_{l_i})_{i\in\mathbb{N}}$.
The second statement is also true in view of Condition \ref{Bintegrator}.\ref{BintegratorPhiBounds} (the condition $\Phi^\infty=0$ implies that $C^\infty$ is continuous), \cref{JtoUConv} and \cref{EquivWeakConv}.
\end{proof}

\begin{lemma}\label{lem:Lk0k-UI}
The sequences $\big( \sup_{t\in\mathbb{R}_+}\bigl\Vert L^{k,(0)}_t\bigr\Vert_2^2\big)_{k\in\overline{\mathbb{N}}}$ and  
$\big(\Vert L_\infty^{k,(0)}\Vert_2^2\big)_{k\in\overline{\mathbb{N}}}$ are uniformly integrable.
\end{lemma}

\begin{proof}
Using Cauchy--Schwarz's inequality as in \cite[Inequality 3.16]{papapantoleon2016existence}, we can obtain for every $t\in[0,\infty)$ and $k\in\overline{\mathbb{N}}$
\begin{gather*}
	\bigl\Vert L^{k,(0)}_t\bigl\Vert_2^2 
		\le  \frac{1}{\hat\beta}\int_{(0,t]} \textup{e}^{\hat{\beta} A^k_s} \frac{\bigl\Vert f^k(s,0,0,\mathbf{0})\bigr\Vert_2^2}{(\alpha_s^{k})^2} \d C_s^k
		\le  \frac{1}{\hat\beta}\int_{(0,\infty)} \textup{e}^{\hat{\beta} A^k_s} \frac{\bigl\Vert f^k(s,0,0,\mathbf{0})\bigr\Vert_2^2}{(\alpha_s^{k})^2} \d C_s^k,
\shortintertext{as well as}
\big\Vert L_{\infty}^{k,(0)}\big\Vert_2^2
		\le  \frac{1}{\hat\beta}\int_{(0,\infty)} \textup{e}^{\hat{\beta} A^k_s} \frac{\bigl\Vert f^k(s,0,0,\mathbf{0})\bigr\Vert_2^2}{(\alpha_s^{k})^2} \d C_s^k.
\end{gather*}
In view of \ref{BgeneratorUI}, which states that the right-hand side is uniformly integrable, we obtain the required result by \citeauthor*{he1992semimartingale} \cite[Theorem 1.7.1]{he1992semimartingale}.
\end{proof}

\begin{proposition}\label{prop:first-step-Induction}
The first step of the induction is valid, that is
\begin{gather*}
\begin{multlined}[c][0.87\displaywidth]
\int_{(0,\cdot]} f^k\bigl(s,Y_s^{k,(0)}, Z_s^{k,(0)}, U^{k,(0)}(s,\cdot)\bigr)\d C_{s}^k 
\xrightarrow[k\to\infty]{\hspace{0.2em}\left(\textup{J}_1(\mathbb{R}^\ell),\mathbb L^2\right)}
\int_{(0,\cdot]} f^\infty\bigl(s,Y_s^{\infty,(0)}, Z_s^{\infty,(0)}, U^{\infty,(0)}(s,\cdot)\bigr)\d C_{s}^\infty.
\tag{LS$^{(0)}$}\label{LS0}
\end{multlined}
\shortintertext{}
\begin{multlined}[c][0.87\displaywidth]
\int_{(0,\infty)} f^k\bigl(s,Y_s^{k,(0)}, Z_s^{k,(0)}, U^{k,(0)}(s,\cdot)\bigr)\d C_{s}^k 
\xrightarrow[k\to\infty]{\hspace{0.2em}\mathbb L^2(\Omega,\mathcal{G},\mathbb{P};\mathbb{R}^\ell)\hspace{0.2em}}
\int_{(0,\infty)} f^\infty\bigl(s,Y_s^{\infty,(0)}, Z_s^{\infty,(0)}, U^{\infty,(0)}(s,\cdot)\bigr)\d C_{s}^\infty.
\tag{LS$^{(0)}_\infty$}\label{LS0inf}
\end{multlined}
\end{gather*}
\end{proposition}

\begin{proof}
Before we present the arguments, let us remind the reader that by definition $\bigl(Y^{k,(0)}, Z^{k,(0)}, U^{k,(0)}\bigr):=(0,0,\mathbf{0})$ for every $k\in\overline{\mathbb{N}}$.
We are going to apply Vitali's theorem, \emph{i.e.} we will prove initially Convergence \eqref{LS0} and \eqref{LS0inf} in probability and then that the sequences of the respective $\Vert \cdot \Vert_2^2$-norms are uniformly integrable.
The latter has been proved in \cref{lem:Lk0k-UI}.
The former will be proved by means of \citeauthor{dudley2002real} \cite[Theorem 9.2.1]{dudley2002real}.

\medskip
\noindent To this end, let us consider a subsequence $\bigl(L^{k_l,(0)}_{\infty}\bigr)_{l\in\mathbb{N}}$.
By \cref{lem:Conv-Ck} there exists a further subsequence $(k_{l_i})_{i\in\mathbb{N}}$ such that $\mu_{C^{k_{l_i}}}\xrightarrow[i\to\infty]{\hspace{0.2em}\textup{w}\hspace{0.2em}}\mu_{C^\infty}$, $\mathbb P$--almost surely.
Consequently, we have also that $\sup_{i\in\overline{\mathbb{N}}}\mu_{C^{k_{l_i}}}(\mathbb{R}_+)<\infty$, $\mathbb P$--almost surely.
In view of Condition \ref{Bgenerator}.\ref{Bgenerator:respects-J1-conv}, we can apply \Cref{prop:MOAppl} as well as \Cref{cor:MOAppl}.
Therefore, the subsequence $(L_\infty^{k_{l_i},(0)})_{i\in\mathbb{N}}$ converges $\mathbb P$--a.s. to $L_\infty^{\infty,(0)}$, and the subsequence 
$(L_\cdot^{k_{l_i},(0)})_{i\in\mathbb{N}}$ converges under the $\text{J}_1$-topology, $\mathbb P$--a.s., to $(L_\cdot^{\infty,(0)})_{i\in\mathbb{N}}$.
\end{proof}
\begin{corollary}\label{cor:Step-1}
The convergences 
\begin{gather*}
\big(Y^{k,(1)}, Z^{k,(1)}\cdot X^{k,\circ} +  U^{k,(1)}\star\widetilde{\mu}^{k,\natural}, N^{k,(1)}\big) 		
\xrightarrow{\hspace{0.2em}\left(\textup{J}_1(\mathbb R^{\ell \times 3}),\mathbb L^2\right)\hspace{0.2em}} 	
\big(Y^{\infty,(1)}, Z^{\infty,(1)}\cdot X^{\infty,\circ} + U^{\infty,(1)}\star \widetilde{\mu}^{\infty,\natural},0 \big),\\
\shortintertext{}
\begin{multlined}[c][\displaywidth]
\big([ Y^{k,(1)}],[ Z^{k,(1)}\cdot X^{k,\circ} + U^{k,(1)}\star\widetilde{\mu}^{k,\natural}], [N^{k,(1)}], [ Y^{k,(1)},X^{k,\circ}],[ Y^{k,(1)},X^{k,\natural}], [Y^{k,(1)},N^{k,(1)}] \big)\\		
\hspace{-0.5em} \xrightarrow{\hspace{0.1em}\left(\textup{J}_1(\mathbb R^{\ell \times (3\ell+m+n+\ell)} ),\mathbb L^1\right)\hspace{0.1em}}
\big([ Y^{\infty,(1)}], [ Z^{\infty,(1)}\cdot X^{\infty,\circ} +  U^{\infty,(1)}\star\widetilde{\mu}^{\infty,\natural}], 0, [ Y^{\infty,(1)},X^{\infty,\circ}],  [ Y^{\infty,(1)},X^{\infty,\natural}], 0\big),
\end{multlined}
\shortintertext{}
\begin{multlined}[c][\displaywidth]
\big(\langle Y^{k,(1)}\rangle,\langle Z^{k,(1)}\cdot X^{k,\circ} \rangle, \langle U^{k,(1)}\star\widetilde{\mu}^{k,\natural}\rangle, \langle N^{k,(1)}\rangle, \langle Y^{k,(1)},X^{k,\circ}\rangle,\langle Y^{k,(1)},X^{k,\natural}\rangle,\langle Y^{k,(1)},N^{k,(1)}\rangle \big)\\		
\hspace{-0.5em}\xrightarrow{\hspace{0.1em}\left(\textup{J}_1(\mathbb R^{\ell \times (4\ell+m+n+\ell)}),\mathbb L^1\right)\hspace{0.1em}}
\big(\langle Y^{\infty,(1)}\rangle, \langle Z^{\infty,(1)}\cdot X^{\infty,\circ} \rangle, \langle U^{\infty,(1)}\star\widetilde{\mu}^{\infty,\natural}\rangle, 0,\langle Y^{\infty,(1)},X^{\infty,\circ}\rangle,  \langle Y^{\infty,(1)},X^{\infty,\natural}\rangle, 0\big),
\end{multlined}
\end{gather*}
are valid, where $0$ denotes the zero process whose state space is a finite-dimensional Euclidean space.
\end{corollary}
\begin{proof}
Apply \cite[Corollary 3.10]{papapantoleon2019stability} for the sequence $(M^{k,(0)})_{k\in\overline{\mathbb{N}}}$, which allows us to conclude.
\end{proof}
%
%
%
%

\begin{remark}\label{YLimitInfty}
\begin{enumerate}[itemindent=0.cm,leftmargin=*]
\item Observe that $\displaystyle Y^{k,(1)}_{\infty}$ is well-defined $\mathbb P$--almost surely, for every $k\in\overline{\mathbb{N}}$.
Indeed, by \eqref{def:Mk}
\begin{gather*}
Y^{k,(1)}_t	= M^{k,(0)}_t - \int_{(0,t]} f^k\bigl(s,Y_s^{k,(0)}, Z_s^{k,(0)},U^{k,(0)}(s,\cdot)\bigr)\d C_{s}^k,\; \text{ for every }k\in\overline{\mathbb{N}}.
\end{gather*}
The limit as $t\longrightarrow\infty$ of the martingale part exists due to its square integrability and the limit of the Lebesgue--Stieltjes integral exists $\mathbb P$--almost surely$;$ see {\rm\cref{lem:Ck-PUT}}.

\medskip
\item Although we have accomplished the aim of this section, we will need to obtain an additional result.
Namely, we need to complement our induction hypothesis on the $p$-step with the assumption that the uniform integrability of the sequence $(\Gamma^{k,(p+1)})_{k\in\overline{\mathbb{N}}}$ $($recall {\rm Notation \ref{notation:LS-mpn}}$)$ is inherited by the uniform integrability of the sequences $(\Gamma^{k,(p)})_{k\in\overline{\mathbb{N}}}$ and of the standard data. 
The following lemma serves this aim.
\end{enumerate}
\end{remark}

\begin{lemma}\label{lem:UI-Gamma}
The sequence $(\Gamma^{k,(1)})_{k\in\overline{\mathbb{N}}}$ is uniformly integrable.
\end{lemma}
\begin{proof}
We will use the Lipschitz property of the generator $f^k$ as well as the fact that 
\begin{align}\label{prop:Lipschitz-bounds}
\frac{r^k}{(\alpha^k)^2} \le (\alpha^k)^2, \;
\frac{\theta^{k,\circ}}{(\alpha^k)^2} \le 1, \;
\text{and} \;
\frac{\theta^{k,\natural}}{(\alpha^k)^2} \le 1;
\end{align}
see the definition of $\alpha^k$ in \ref{Ffour}. Let $k\in\overline{\mathbb{N}}$. 
By definition of $\Gamma^{k,(1)}$, the Lipschitz property of the generator $f^k$ and the definition of $A^k$
\begin{align*}
\Gamma^{k,(1)}
&=
\int_{(0,\infty)}{}{} \frac{\bigl\Vert  f^k\big(s, Y^{k,(1)}_s, Z^{k,(1)}_s, U^{k,(1)}(s,\cdot)\big)\bigr\Vert_2^2}{(\alpha^k_s)^2}\d C^k_s\\
&\le
\int_{(0,\infty)}{}{} \Big[(\alpha^k_s)^2\bigl\Vert Y^{k,(1)}_s\bigr\Vert_2^2 + \bigl\Vert c^k_s Z^{k,(1)}_s\bigr\Vert_2^2 + \tnormbig{U^{k,(1)}_s(\cdot)}_{k,s}\Big]\d C^k_s
+ \int_{(0,\infty)}{}{} \frac{\bigl\Vert  f^k\big(s, 0, 0, \mathbf{0}\big)\bigr\Vert_2^2}{(\alpha^k_s)^2}\d C^k_s\\
&=
\int_{(0,\infty)}{}{} \bigl\Vert Y^{k,(1)}_s\bigr\Vert_2^2 \d A^k_s + 
\int_{(0,\infty)}{}{} \d  \textup{Tr}\Big[\big\langle Z^{k,(1)}\cdot X^{k,\circ} + U^{k,(1)}\star\widetilde{\mu}^{k,\natural}\big\rangle_s\Big]
+ \Gamma^{k,(0)}\\
&
\le
2\overline{A} \sup_{t\in\mathbb{R}_+}\bigl\Vert  \mathbb{E}[ \xi^k |\mathcal{G}^k_t]\bigr\Vert_2^2
+2\overline{A} \sup_{t\in\mathbb{R}_+}\bigg\Vert  \mathbb{E}\bigg[ \int_{(t,\infty)}f^k(u,0,0,\mathbf{0}) \d C^k_u\Big|\mathcal{G}^k_t\bigg] \bigg\Vert_2^2 
+ \textup{Tr}\Big[\langle M^{k,(0)}\rangle_\infty\Big]
+ \Gamma^{k,(0)}.
\numberthis{}\label{UI:Gammak1}
\end{align*}
In the first inequality, we used the Lipschitz property and the definition of $\alpha^k$ (see \ref{Ffour}) as well as Inequalities \eqref{prop:Lipschitz-bounds}.
In the second equality, we used the definition of $A^k,c^{k}$ and $\tnorm{\cdot}_{k,t}$; see also \cite[Identity (2.9)]{papapantoleon2016existence}.
In the second Inequality, we used Identity \eqref{PicardBSDEorthogonal}, \emph{i.e.}, 
\begin{align*}
Y^{k,(1)}_t = \mathbb{E}\bigg[ \xi^k + \int_{(t,\infty)} f^k(s, 0,0,\mathbf{0})\d C_s^k \bigg| \mathcal{G}_t^k\bigg].
\end{align*}
We only need to observe now that the summands on the right-hand side of \eqref{UI:Gammak1} form a uniformly integrable sequence, for $k\in\overline{\mathbb{N}}$, as sum of elements of uniformly integrable random variables; see \citeauthor{he1992semimartingale} \cite[Corollary 1.10]{he1992semimartingale}. 
Indeed,
\begin{itemize}
	\item for the sequence associated to the first summand  we have that $(\Vert \xi^k \Vert_2^2)_{k\in\overline{\mathbb{N}}}$ is uniformly integrable by Vitali's theorem.
	Then, we can conclude the uniform integrability of the required sequence by \cref{cor:DoobApplication};
	
	\smallskip
	\item for the sequence associated to the second summand, \emph{i.e.}, 
	\begin{align}\label{seq:norm-Lk-0}
	\bigg(\sup_{t\in\mathbb{R}_+}\bigg\Vert \mathbb{E}\bigg[  \int_{(t,\infty)} f^k(s, 0,0,\mathbf{0})\d C_s^k \bigg| \mathcal{G}_t^k\bigg]\bigg\Vert_2^2\bigg)_{k\in\overline{\mathbb{N}}},
	\end{align}
	we can prove by means of the conditional Jensen inequality and the conditional Cauchy--Schwarz inequality, that for all $k\in\overline{\mathbb{N}}$ holds
		\begin{align*}
	\bigg\Vert  \mathbb{E}\bigg[  \int_{(t,\infty)} f^k(s, 0,0,\mathbf{0})\d C_s^k \bigg| \mathcal{G}_t^k\bigg]\bigg\Vert_2^2
	\le \bigg\Vert  \mathbb{E}\bigg[  \int_{(0,\infty)} |f|^k(s, 0,0,\mathbf{0})\d C_s^k \Big| \mathcal{G}_t^k\bigg]\bigg\Vert_2^2
	\le  \mathbb{E}\big[ \Gamma^{k,(0)} \big| \mathcal{G}^k_t\big] . 
	\end{align*}
	Since $(\Gamma^{k,(0)})_{k\in\overline{\mathbb{N}}}$ is uniformly integrable (see Condition \ref{BgeneratorUI}) we can conclude the uniform integrability of \eqref{seq:norm-Lk-0} by \citeauthor{he1992semimartingale} \cite[Theorem 1.7, Theorem 1.8]{he1992semimartingale}.
	Then, we can conclude the uniform integrability of the required sequence by \cref{cor:DoobApplication};
	
	\smallskip
	\item for the sequence  associated to the third summand we use the fact (which is true in view of the validity of \eqref{LS0}) that
	\begin{align*}
		M^{k,(0)}_\infty \xrightarrow[k\to\infty]{\hspace{1em} \mathbb{L}^2(\Omega,\mathcal{G}, \mathbb{P})\hspace{1em}}  M^{\infty,(0)}_\infty,
	\end{align*}
	which implies\footnote{The uniform integrability of $\big(\text{Tr}\big[[M^{k,(0)}]_\infty\big]\big)_{k\in\overline{\mathbb{N}}}$ can also be deduced.} the uniform integrability of $\big(\text{Tr}\big[\langle M^{k,(0)}\rangle_\infty\big]\big)_{k\in\overline{\mathbb{N}}}$;
	
	\smallskip
	\item finally, the sequence associated to the last summand is uniformly integrable by \ref{BgeneratorUI}. \qedhere
\end{itemize}
\end{proof}

\begin{remark}\label{rem:ExtendInductionAssumption}
The validity of {\rm\cref{lem:UI-Gamma}} allows us indeed to complement our induction step with the statement 
\begin{equation*}
\text{\rm the sequence }	\big( \Gamma^{k,(p)})_{k\in\overline{\mathbb{N}}}\text{\rm is uniformly integrable.} \tag{$\textup{UI}^{(p)}$}
\end{equation*}
The reader may observe that the above property for $p=1$ transfers the uniform integrability to the sequence
\begin{align}\label{prop:UI-Lk1}
\bigg(\sup_{t\in\mathbb{R}_+}\bigl\Vert L_{t}^{k,(1)}\bigr\Vert_2^2\bigg)_{ k\in\overline{\mathbb{N}}}.
\end{align}
This is immediate by applying Cauchy--Schwarz's inequality to $\big(\bigl\Vert L_{t}^{k,(1)}\bigr\Vert_2^2\big)_{k\in\overline{\mathbb{N}}}$ $($as in {\rm\cite[Inequality (3.16)]{papapantoleon2016existence}}$)$ and observing that $(\Gamma^{k,(1)})_{k\in\overline{\mathbb{N}}}$ dominates the sequence $\big(\bigl\Vert L_{t}^{k,(1)}\bigr\Vert_2^2\big)_{k\in\overline{\mathbb{N}}}$, for every $t\in\mathbb{R}_+$.

\end{remark}

%% file: Stability_BSDEJ_Induction_pStep_Remarks.tex

In this sub-section, we assume that Convergences  
\begin{gather*}
L^{k,(p-1)}
\xrightarrow[k\to\infty]{\hspace{0.2em}\left(\textup{J}_1(\mathbb{R}^{\ell}),\mathbb L^2\right)\hspace{0.2em}}
L^{\infty,(p-1)},\tag{\text{LS$^{(p-1)}$}}\label{LSpminus1}
\shortintertext{}
L^{k,(p-1)}_\infty
\xrightarrow[k\to\infty]{\hspace{0.2em}\mathbb L^2(\mathcal G;\mathbb{R})\hspace{0.2em}}
L^{\infty,(p-1)}_\infty,
\tag{\text{LS$^{(p-1)}_\infty$}}\label{LSInfinitypminus1}
\end{gather*}
as well as the statement
\begin{align*}
\text{the sequence }\big(\Gamma^{k,(p-1)}\big)_{k\in\overline{\mathbb{N}}}\text{ is uniformly integrable,}
\tag{$\textup{UI}^{(p-1)}$}\label{Gamma-pminus1-UI}
\end{align*}
are true for some arbitrary but fixed $p\in\mathbb{N}$.
Then, we will prove that Convergences \eqref{LSp} and \eqref{LSpInfinity}, as well as the statement
\begin{equation*}\label{Gamma-p-UI}
\text{the sequence }    \big( \Gamma^{k,(p)})_{k\in\overline{\mathbb{N}}}\text{ is uniformly integrable} \tag{$\textup{UI}^{(p)}$}
\end{equation*}
are also true.

\medskip
Compared to the first step of the induction, the $p$-th step is more involved.
Let us thus briefly explain the approach we are going to follow in order to reduce the complexity.
In view of Vitali's theorem, it is sufficient to prove initially that Convergences \eqref{LSp} and \eqref{LSpInfinity} hold in probability and then we have to prove that the sequences are (sufficiently) uniformly integrable.
In order to obtain the aforementioned convergence in probability, we are going to use that 
$\bigl(\mathbb D(\mathbb{R}^{\ell}),\delta_{\textup{J}_1(\mathbb{R}^{\ell})}\bigr)$ and $\bigl(\mathbb{R}^{\ell},\Vert\cdot\Vert_2\bigr)$ are Polish spaces, as we did in the first step of the induction.
Therefore, in view of \citeauthor{dudley2002real} \cite[Theorem 9.2.1]{dudley2002real}, it is sufficient to prove that from every subsequence $\big(L^{k_l,(p)}\big)_{l\in\mathbb{N}}$, resp. $\big(L^{k_l,(p)}_\infty\big)_{l\in\mathbb{N}}$,
we can extract a further subsequence $\big(L^{k_{l_i},(p)}\big)_{i\in\mathbb{N}}$, resp. $\big(L^{k_{l_i},(p)}_\infty\big)_{i\in\mathbb{N}}$, such that
\begin{gather}
L^{k_{l_i},(p)}
\xrightarrow[i\to\infty]{\hspace{0.2em}\textup{J}_1(\mathbb{R}^{\ell})\hspace{0.2em}}
L^{\infty,(p)},\; \mathbb P\text{\rm--a.s.} ,
\label{LSp-Pas}
\shortintertext{resp.}
L^{k_{l_i},(p)}_\infty
\xrightarrow[i\to\infty]{\hspace{0.2em}\Vert\cdot\Vert_2\hspace{0.2em}}
L^{\infty,(p)}_\infty,\; \mathbb P\text{\rm--a.s.}
\label{LSpInfinity-Pas}
\end{gather}
Equivalently, for the given subsequence $(k_{l_i})_{i\in\mathbb{N}}$, there exists a set $\widetilde{\Omega}\in \mathcal{G}$ with $\mathbb{P}[\widetilde{\Omega}]=1$ such that
\begin{gather*}
\limsup_{i\to\infty} \Big\{\omega\in \widetilde{\Omega} : 
    \delta_{\textup{J}_1(\mathbb{R}^{\ell})}\big(L^{k_{l_i},(p)}(\omega), L^{\infty,(p)}(\omega)\big)>\varepsilon\Big\}=\emptyset,\; \forall \varepsilon >0,
\shortintertext{resp.}
\limsup_{i\to\infty} \Big\{\omega\in \widetilde{\Omega} : \Vert L^{k_{l_i},(p)}_\infty(\omega) - 
L^{\infty,(p)})_\infty(\omega)\Vert_2>\varepsilon\Big\}=\emptyset,\; \forall \varepsilon >0.
\end{gather*}

To this end, let us consider Convergence \eqref{LSpInfinity-Pas}, fix an $\varepsilon >0$ and assume that for $\mathbb P$--almost every $\omega \in \Omega$ there exist 
$\widetilde{Z}^{\varepsilon,(p)}(\omega,\cdot)\in D^{\circ,\ell\times m}$ and 
$\widetilde{U}^{\varepsilon,(p)}(\omega,\cdot) \in D^{\natural}$ such that 
\begin{equation}\label{Approx:Descr}
\bigg\Vert  \int_{(0,\infty)} f^\infty\bigl(s,Y^{\infty,(p)}_s,\widetilde{Z}^{\varepsilon,(p)}_s,\widetilde{U}^{\varepsilon,(p)}(s,\cdot)\bigr)\,\d C^\infty_s - L^{\infty,(p)}_\infty\bigg\Vert_2 <\frac{\varepsilon}{3},\; \mathbb P\text{\rm--a.s.} 
\end{equation}
Then, using the set inclusion, where $F^k$ are assumed to be $\mathbb{R}^{\ell}$-valued random functions for $k\in\{1,2,3\}$
\begin{equation*}
\bigg\{\omega\in \Omega: \bigg\Vert \sum_{k=1}^3F^k(\omega)\bigg\Vert_2 >\varepsilon\bigg\}
    \subseteq 
    \bigcup_{k=1}^3\bigg\{\omega\in \Omega:\bigl\Vert F^k(\omega)\bigr\Vert_2>\frac{\varepsilon}{3}\bigg\},
\end{equation*}
we can obtain Convergence \eqref{LSpInfinity-Pas} if we can find an $\widetilde{\Omega}\subset \Omega$ with $\mathbb{P}[\widetilde{\Omega}]=1$ such that\footnote{For notational convenience, we index the $k_{l_i}$-element in the next expression simply by $i$.}
\begin{gather}
\limsup_{i\to\infty}\biggl\{  \omega\in\widetilde{\Omega} :
\bigg\Vert  
L^{i,(p)}_\infty (\omega)- 
\bigg(\int_{(0,\infty)} f^{i}\bigl(s,Y_s^{i,(p)},\widetilde{Z}^{\varepsilon,(p)}_s,\widetilde{U}^{\varepsilon,(p)}(s,\cdot)\bigr)\,\d C^{i}_s\bigg)(\omega)\bigg\Vert_2 >\frac{\varepsilon}{3} \biggr\}=\emptyset, \label{limsup1}\\
\shortintertext{and}
\begin{multlined}[0.9\textwidth]
\limsup_{i\to\infty}\biggl\{  \omega\in\widetilde{\Omega} :  
\bigg\Vert  
\bigg(\int_{(0,\infty)} 
f^{i}\bigl(s,Y_s^{i,(p)},\widetilde{Z}^{\varepsilon,(p)}_s,\widetilde{U}^{\varepsilon,(p)}(s,\cdot)\bigr)\,\d C^{i}_s\bigg)(\omega)\\
-
\bigg(\int_{(0,\infty)} f^\infty\bigl(s,Y^{\infty,(p)}_s,\widetilde{Z}^{\varepsilon,(p)}_s,\widetilde{U}^{\varepsilon,(p)}(s,\cdot)\bigr)\,\d C^\infty_s\bigg)(\omega)\bigg\Vert_2 >\frac{\varepsilon}{3}\biggl\}=\emptyset.
\end{multlined}\label{limsup2}
\end{gather}
Hence, we prove \eqref{LSpInfinity-Pas} if \eqref{Approx:Descr}, \eqref{limsup1} and \eqref{limsup2} are true.
An analogous decomposition can be done for \eqref{LSp-Pas}, where the distance is measured by the $\delta_{\textup{J}_1(\mathbb{R}^{\ell})}$-metric.
Returning now to the uniform integrability that the sequences should satisfy, we will need to prove that the family $\big(\sup_{t\in\mathbb{R}_+}\Vert L^{k,(p)}_t\Vert_2^2\big)_{k\in\overline{\mathbb{N}}}$ is uniformly integrable, 
which is a sufficient condition for concluding both the Convergence \eqref{LSp} and \eqref{LSpInfinity}.


\begin{iassumption}
From now on we fix an arbitrary subsequence $\bigl(L^{k_l,(p)}\bigr)_{l\in\mathbb{N}}$, resp. $\bigl(L^{k_l,(p)}_\infty\bigr)_{l\in\mathbb{N}}$.
\end{iassumption}


\vspace{0.5em}
Let us conclude the description of our strategy by collecting all the information we have available for the next subsections.
We will state them as a remark so that they are easily referred to. 
Moreover, for notational convenience, we can assume that the sequence for which the forthcoming convergences are obtained $\mathbb P$--almost surely is indexed by $(k_{l_i})_{i\in\mathbb{N}}$.
This can be done without loss of generality, since we can pass to a further subsequence finitely many times. 
\begin{remark}\label{rem:Conv-Subseq-klm}
\begin{enumerate}[itemindent=0.cm,leftmargin=*]
    \item\label{rem:Conv-Subseq-klm-1} By {\rm\cref{lem:Conv-Ck}}, there exist a $\delta_{\textup{J}_1(\mathbb{R}^{\ell})}$-convergent $(\mathbb{P}\text{\rm--a.s.})$ subsequence $\bigl(C^{k_{l_i}}\bigr)_{i\in\overline{\mathbb{N}}}$ as well as a $\delta_{\vert \cdot \vert}$-convergent $(\mathbb{P}\text{\rm--a.s.})$ subsequence $\bigl(C^{k_{l_i}}_\infty\bigr)_{i\in\overline{\mathbb{N}}}$.
\medskip
    \item\label{rem:Conv-Subseq-klm-2}  By {\rm Condition \ref{Bgenerator}.\ref{Bgenerator:respects-J1-conv}}, we can assume in particular that 
\begin{equation*}
  f^{k_{l_i}}\bigl(\cdot,0,0,\mathbf 0\bigr)\xrightarrow[i\to\infty]{\hspace{0.2em}\textup{J}_1(\mathbb{R}^{\ell})\hspace{0.2em}}  f^{\infty}\bigl(\cdot,0,0,\mathbf 0\bigr),\; \mathbb P\text{\rm--a.s.}
\end{equation*}
    \item\label{rem:Conv-Subseq-klm-3}
The convergences in \eqref{LSpminus1} and \eqref{LSInfinitypminus1}, which are assumed true $($this is the induction assumption$)$, allow us to obtain that 
\begin{equation*}
L^{{k_{l_i},(p-1)}}
\xrightarrow[i\to\infty]{\hspace{0.2em}\textup{J}_1(\mathbb{R}^{\ell})\hspace{0.2em}}
L^{\infty,(p-1)}, \;
\text{\rm as well as}\;
L^{k_{l_i},(p-1)}_\infty
\xrightarrow[i\to\infty]{\hspace{0.2em}\Vert\cdot\Vert_2\hspace{0.2em}}
L^{\infty,(p-1)}_\infty, \; \mathbb P\text{\rm--a.s.}
\end{equation*} 

\item\label{rem:Conv-Subseq-klm-4} 
In view of the discussion made in the outline of the proof of {\rm\cref{BSDERobMainTheorem}}, see in particular on {\rm Page \pageref{conv:Picardfixedp-partial}}, 
the validity of the convergences in \eqref{LSpminus1} and \eqref{LSInfinitypminus1} allows us to apply {\rm\cite[Corollary 3.10]{papapantoleon2019stability}} for the martingale sequence $(M^{k,(p-1)})_{k\in\overline{\mathbb{N}}}.$
More precisely, the convergences
\begin{gather*}
\big(Y^{k,(p)}, Z^{k,(p)}\cdot X^{k,\circ} +  U^{k,(p)}\star\widetilde{\mu}^{k,\natural}, N^{k,(p)}\big)        
\xrightarrow{\hspace{0.2em}\left(\textup{J}_1(\mathbb R^{\ell \times 3}),\mathbb L^2\right)\hspace{0.2em}}  
\big(Y^{\infty,(p)}, Z^{\infty,(p)}\cdot X^{\infty,\circ} + U^{\infty,(p)}\star \widetilde{\mu}^{\infty,\natural},0 \big),
\\
\shortintertext{}
\begin{multlined}[c][\displaywidth]
\big([ Y^{k,(p)}],[ Z^{k,(p)}\cdot X^{k,\circ} + U^{k,(p)}\star\widetilde{\mu}^{k,\natural}], [N^{k,(p)}], [ Y^{k,(p)},X^{k,\circ}],[ Y^{k,(p)},X^{k,\natural}], [Y^{k,(p)},N^{k,(p)}] \big)\\       
\hspace{-0.3em}\xrightarrow{\hspace{0.1em}\left(\textup{J}_1(\mathbb R^{\ell \times (3\ell+m+n+\ell)}),\mathbb L^1\right)\hspace{0.1em}}
\big([ Y^{\infty,(p)}], [ Z^{\infty,(p)}\cdot X^{\infty,\circ} + U^{\infty,(p)}\star\widetilde{\mu}^{\infty,\natural}], 0, [ Y^{\infty,(p)},X^{\infty,\circ}],  [ Y^{\infty,(p)},X^{\infty,\natural}], 0\big),
\end{multlined}
\shortintertext{}
\begin{multlined}[c][\displaywidth]
\big(\langle Y^{k,(p)}\rangle,\langle Z^{k,(p)}\cdot X^{k,\circ} \rangle, \langle U^{k,(p)}\star\widetilde{\mu}^{k,\natural}\rangle, \langle N^{k,(p)}\rangle, \langle Y^{k,(p)},X^{k,\circ}\rangle,\langle Y^{k,(p)},X^{k,\natural}\rangle,\langle Y^{k,(p)},N^{k,(p)}\rangle \big)\\      
\hspace{-0.3em}\xrightarrow{\hspace{0.1em}\left(\textup{J}_1(\mathbb R^{\ell \times (4\ell+m+n+\ell)}),\mathbb L^1\right)\hspace{0.1em}}
\big(\langle Y^{\infty,(p)}\rangle, \langle Z^{\infty,(p)}\cdot X^{\infty,\circ} \rangle, \langle U^{\infty,(p)}\star\widetilde{\mu}^{\infty,\natural}\rangle, 0,\langle Y^{\infty,(p)},X^{\infty,\circ}\rangle,  \langle Y^{\infty,(p)},X^{\infty,\natural}\rangle, 0\big),
\end{multlined}
\end{gather*}
are valid.
Here $0$ denotes the zero process whose state space is a finite-dimensional Euclidean space.
For later reference, we state only the results we are going to make use of
\begin{gather*}
\begin{multlined}[c][0.85\displaywidth]
\big(\langle Z^{k,(p)}\cdot X^{k,\circ}\rangle, \langle U^{k,(p)}\star \widetilde\mu^{k,\natural}\rangle, \langle Z^{k,(p)}\cdot X^{k,\circ}, X^{k,\circ}\rangle, \langle U^{k,(p)}\star \widetilde\mu^{k,\natural}, X^{k,\natural}\rangle\big) 
\xrightarrow[k\to\infty]{\hspace{0.2em}\left(\textup{J}_1(\mathbb{R}^{\ell\times (2\ell+m+n)}),\mathbb L^1\right)\hspace{0.2em}}\\
\big(\langle Z^{\infty,(p)}\cdot X^{\infty,\circ}\rangle, \langle U^{\infty,(p)}\star \widetilde\mu^{\infty,\natural}\rangle, \langle Z^{\infty,(p)}\cdot X^{\infty,\circ}, X^{\infty,\circ}\rangle, \langle U^{\infty,(p)}\star \widetilde\mu^{\infty,\natural}, X^{\infty,\natural}\rangle\big).
\numberthis\label{Conv:Angle-Brack-p-step}
\end{multlined}
\end{gather*}
Additionally, we can apply {\rm\cref{cor:DoobApplication}} for the sequence $\big(\big\Vert M^{k,(p-1)}_\infty\big\Vert_2^2\big)_{k\in\overline{\mathbb{N}}}$ in order to obtain that 
\begin{equation}\label{UI:supM-pminus1}
    \bigg(\sup_{t\in\mathbb{R}_+} \big\Vert M^{k,(p-1)}_t\big\Vert_2^2\bigg)_{k\in\overline{\mathbb{N}}} \text{\rm is uniformly integrable.}
\end{equation}
\item\label{UI:Delta_p}
The sequence $(\Delta^{k,(p)})_{k\in\overline{\mathbb{N}}}$ is uniformly integrable, since it is strongly majored\footnote{See \cite[Definition VI.3.34]{jacod2003limit}.} by $(\textup{Tr}[\langle M^{k,(p-1)}\rangle_\infty])_{k\in\overline{\mathbb{N}}}$ which is uniformly integrable. 
The uniform integrability of the latter sequence is derived by the uniform integrability of the sequence $\big(\Vert M^{k,(p-1)}_\infty \Vert_2^{2}\big)_{k\in\overline{\mathbb{N}}}$. 
\smallskip
\item\label{rem:Conv-Subseq-klm-7}
In view of {\rm Conditions \ref{MFilqlc}--\ref{Mfinalrv}}, we can apply {\rm\citeauthor{memin2003stability} \cite[Corollary 12]{memin2003stability}}.
Therefore, we may assume that 
\begin{align}\label{conv:Memin-Cor}
\big(\overline{X}^{i}, [\overline{X}^i],\langle\overline{X}^i\rangle\big)
\xrightarrow[i\to\infty]{\hspace{0.2em}\textup{J}_1(\mathbb{R}^{(m+n)\times\{1+2(m+n)\}})\hspace{0.2em}}
\big(\overline{X}^{\infty}, [\overline{X}^\infty],\langle\overline{X}^\infty\rangle\big),\; \mathbb{P}\text{\rm--a.s.}
\end{align}
In particular $X^{i,\natural}
\xrightarrow[i\to\infty]{\hspace{0.2em}\textup{J}_1(\mathbb{R}^{n})\hspace{0.2em}}
X^{\infty,\natural},\; \mathbb{P}\text{\rm--a.s.},$
which, in conjunction with {\rm\citeauthor{jacod2003limit} \cite[Corollary VI.2.8]{jacod2003limit}}, further implies 
\begin{align*}
g * \mu^{i,\natural}
\xrightarrow[i\to\infty]{\hspace{0.2em}\textup{J}_1(\mathbb{R}^{\ell})\hspace{0.2em}}
g * \mu^{\infty,\natural},\; \mathbb{P}\text{\rm--a.s.,}\;
\text{\rm for every } g\in D^{\natural}.
\end{align*}
On the other hand, for every $i\in\overline{\mathbb{N}}$ holds
\begin{align*}
g\star \widetilde{\mu}^{i,\natural} = g* \mu^{i,\natural} - g*\nu^{i,\natural},\; \mathbb{P}\text{\rm--a.s.},
\text{\rm for every } g\in D^{\natural}.
\end{align*}
Since $D^{\natural}$ is countable, we can assume that this is true for every $\omega \in \overline{\Omega}$, for some $\overline{\Omega}\in\mathcal{G}$ with $\mathbb{P}[\overline{\Omega}]=1$. Finally, by the last convergence and {\rm \citeauthor{memin2003stability} \cite[Theorem 11, Corollary 12]{memin2003stability}}, we have that for every $g\in D^{\natural}$
\begin{align}\label{conv:g-compensator}
\big(g \star \widetilde{\mu}^{i,\natural},g * \nu^{i,\natural}, \langle g \star \widetilde{\mu}^{i,\natural} \rangle\big)
\xrightarrow[i\to\infty]{\hspace{0.2em}\textup{J}_1(\mathbb{R}^{\ell \times (2+\ell) })\hspace{0.2em}}
\big(g \star \widetilde{\mu}^{\infty,\natural},g * \nu^{\infty,\natural},\langle g \star \widetilde{\mu}^{\infty,\natural} \rangle\big),\; \mathbb{P}\text{\rm--a.s.}
\end{align}

\medskip
\item\label{rem:Conv-Subseq-klm-5} In view of \ref{rem:Conv-Subseq-klm-4} and \ref{rem:Conv-Subseq-klm-7} we can have 
    that
\begin{gather}
\begin{multlined}[c][0.85\displaywidth]
\big(M^{i,(p-1)},\langle Z^{i,(p)}\cdot X^{i,\circ}\rangle, \langle U^{i,(p)}\star \widetilde{\mu}^{i}\rangle\big) 
\xrightarrow[i\to\infty]{\hspace{0.2em}\textup{J}_1(\mathbb{R}^{\ell\times(1+2\ell )})\hspace{0.2em}}\\
\big(M^{\infty,(p-1)},\langle Z^{\infty,(p)}\cdot X^{\infty,\circ}\rangle, \langle U^{\infty,(p)}\star \widetilde\mu^{\infty,\natural}\rangle\big),\; \mathbb{P}\text{\rm--a.s.},
\label{Conv:Angle-Brack-kp-Pas}
\end{multlined}
\shortintertext{}
\begin{multlined}[c][0.85\displaywidth]
(\langle Z^{k_{l_i},(p)}\cdot X^{k_{l_i},\circ}, X^{k_{l_i},\circ}\rangle, \langle U^{k_{l_i},(p)}\star \widetilde\mu^{k_{l_i},\natural}, g\star \widetilde\mu^{k_{l_i},\natural}\rangle) 
\xrightarrow[i\to\infty]{\hspace{0.2em}\textup{J}_1(\mathbb{R}^{\ell\times (m+n)})\hspace{0.2em}}\\
(\langle Z^{\infty,(p)}\cdot X^{\infty,\circ}, X^{\infty,\circ}\rangle, 
\langle U^{\infty,(p)}\star \widetilde\mu^{\infty,\natural}, g\star \widetilde\mu^{\infty,\natural}\rangle),\; \mathbb{P}\text{\rm--a.s.}
\label{Conv:Angle-Brack-kp-Pas-2}
\end{multlined}
\end{gather}
\item\label{rem:Conv-Subseq-klm-4-a}
Recall \eqref{def:Mk}, \emph{i.e.} $Y^{k,(p)}=M^{k,(p-1)} - L^{k,(p-1)}$ for every $k\in\overline{\mathbb{N}}.$
We claim that  
\begin{equation}\label{UI:supYp}
    \bigg(\sup_{t\in\mathbb{R}_+}\bigl\Vert Y^{k,(p)}_t\bigr\Vert_2^2\bigg)_{k\in\overline{\mathbb{N}}}\text{ is uniformly integrable.}
\end{equation} 
This can be concluded as follows.
We have 
\begin{equation*}
    \sup_{t\in\mathbb{R}_+}\bigl\Vert Y^{k,(p)}_t\bigr\Vert_2^2\le 2 \sup_{t\in\mathbb{R}_+}\bigl\Vert M^{k,(p-1)}_t\bigr\Vert_2^2 + 2 \sup_{t\in\mathbb{R}_+}\bigl\Vert L^{k,(p-1)}_t\bigr\Vert_2^2, \text{ for every } k\in\overline{\mathbb{N}}.
\end{equation*} 
For the sequence $\big(\sup_{t\in\mathbb{R}_+}\big\Vert L^{k,(p-1)}_t\big\Vert_2^2\big)_{k\in\overline{\mathbb{N}}}$, we can conclude its uniform integrability by arguing analogously to {\rm\Cref{prop:UI-Lk1}}, since $(\Gamma^{k,(p-1)})_{k\in\overline{\mathbb{N}}}$ has been  assumed uniformly integrable; see \eqref{Gamma-pminus1-UI}.
Now, we can conclude our initial statement by \eqref{UI:supM-pminus1}. 
\medskip
\item\label{rem:Conv-Subseq-klm-6}
In view of {\rm Condition \ref{Bintegrator}.\ref{BintegratorPhiBounds}}, which implies that $L^{\infty,(p-1)}$ is a continuous process, parts \ref{rem:Conv-Subseq-klm-3} and \ref{rem:Conv-Subseq-klm-5} above and {\rm\citeauthor{jacod2003limit} \cite[Proposition VI.1.123]{jacod2003limit} }
we can conclude that 
\begin{equation}\label{conv:Y-i-p}
    Y^{k_{l_i},(p)}
\xrightarrow[i\to\infty]{\hspace{0.2em}\textup{J}_1(\mathbb{R}^{\ell})\hspace{0.2em}}
    Y^{\infty,(p)},\; \mathbb P\text{\rm--a.s.}
\end{equation}

\item The sets $D^\circ$ and $D^\natural$ are countable. 
Therefore, in view of {\rm Condition} \ref{Bgenerator}.\ref{Bgenerator:respects-J1-conv} we can assume that 
\begin{align*}
    \Big(f^{k_{l_i}}\big(t,Y^{k_{l_i},(p)}_t, \widetilde{Z}^{\varepsilon,(p)}_t,\widetilde{U}^{\varepsilon,(p)}(t,\cdot)\big) \Big)_{t\in\mathbb{R}_+}
\xrightarrow[i\to\infty]{\hspace{0.2em}\textup{J}_1(\mathbb{R}^{\ell})\hspace{0.2em}}
    \Big(f^{\infty}\big(t,Y^{\infty,(p)}_t, \widetilde{Z}^{\varepsilon,(p)}_t,\widetilde{U}^{\varepsilon,(p)}(t,\cdot)\big) \Big)_{t\in\mathbb{R}_+},\;
\mathbb{P}\text{\rm--a.s.}
\end{align*}
\item\label{rem:Conv-Subseq-klm-stopping-time}
In view of \ref{BFiniteStoppingTime}, we may assume that $\tau^{k_{l_i}}\longrightarrow \tau^{\infty}$, $\mathbb{P}\text{\rm--a.s.}$
The finiteness of $\tau^{\infty}$ implies that $\tau^{k_{l_i}}$ is {finally} finite $\mathbb{P}\text{\rm--a.s.}$
\end{enumerate}
\end{remark}

\begin{remark}\label{rem:Omega-sub}
The purpose of the above remarks was not only to collect all the available information,
but also to provide us with a set $\Omega_\textup{sub}\in \mathcal{G}$ with $\mathbb P[\Omega_\textup{sub}]=1$ such that the above properties hold for every $\omega \in \Omega_\textup{sub}$.
In the next subsubsections, whenever we say that a property holds $\mathbb P$--almost surely for a subsequence indexed by $(k_{l_i})_{i\in\overline{\mathbb{N}}}$, the reader should understand that it holds for every $\omega\in\Omega_\textup{sub}.$
Moreover, we will index the elements of the subsequence $(k_{l_i})_{i\in\overline{\mathbb{N}}}$ simply by $i\in\overline{\mathbb{N}}$; when we use the index $k\in\overline{\mathbb{N}}$ we will refer to the initial sequence.
\end{remark}

%% file: Stability_BSDEJ_Induction_pStep_Convergence.tex

\subsubsection{The claim \texorpdfstring{\eqref{Approx:Descr}}{(3.21)} is true}

\begin{proposition}[Lusin approximation]\label{Lusin_approximation}
For every $\varepsilon>0$, there exist $\widetilde{Z}^{\varepsilon,(p)}:\Omega_\textup{sub}\times\mathbb{R}_+ \longrightarrow\mathbb{R}^{\ell\times m} $ and $\widetilde{U}^{\varepsilon,(p)}:\Omega_\textup{sub}\times\mathbb{R}_+\times\mathbb{R}^{n}\longrightarrow\mathbb{R}^\ell$ $($which are defined $\omega$-by-$\omega)$ with the following properties
\begin{enumerate}
\item $\widetilde{Z}^{\varepsilon,(p)}(\omega,\cdot)\in D^{\circ,\ell\times m};$
\item $\widetilde{U}^{\varepsilon,(p)}(\omega,\cdot) \in D^{\natural};$  
\item the pair $\bigl(\widetilde{Z}^{\varepsilon,(p)},\widetilde{U}^{\varepsilon,(p)}\bigr)$\footnote{We will omit the $\omega$ in order to simplify notation. Moreover, we will denote $Z^{\varepsilon,(p)}(\omega,s)$ by $Z^{\varepsilon,(p)}_s$.} satisfies
\begin{gather*}
\int_{(0,\infty)} 
  \textup{Tr}\big[ (Z^{\infty,(p)}_s - \widetilde{Z}^{\varepsilon,(p)}_s)\frac{\d \langle X^{\infty,\circ}\rangle_s}{\d C^\infty_s}(Z^{\infty,(p)}_s - \widetilde{Z}^{\varepsilon,(p)}_s)^\top \big]
  \d C^\infty_s
  <\frac{\varepsilon^2}{18}   \frac{\hat{\beta}}{\textup{e}^{\hat{\beta}\overline{A}}},\; \forall \omega \in \Omega_{\textup{sub}},
\shortintertext{}
  \int_{(0,\infty)} 
  \tnorm{U^{\infty,(p)}(s,\cdot) - \widetilde{U}^{\varepsilon,(p)}(s,\cdot)}_{\infty,s}^2 \d C^{\infty}_s
  <\frac{\varepsilon^2}{18}  \frac{\hat{\beta}}{\textup{e}^{\hat{\beta}\overline{A}}},\; \forall \omega \in \Omega_{\textup{sub}}.
\end{gather*}
\end{enumerate}
In particular, it holds
\begin{align}\label{ineq:LusinApproximation}
  \bigg\Vert  \int_{(0,\infty)} f^\infty\bigl(s,Y_s^{\infty,(p)},\widetilde{Z}^{\varepsilon,(p)}_s,\widetilde{U}^{\varepsilon,(p)}(s,\cdot)\bigr)
  \d C^\infty_s 
   - L^{\infty,(p)}_{\infty}
  \bigg\Vert_2 <\frac{\varepsilon}{3},\; \forall \omega \in \Omega_{\textup{sub}}.
\end{align}
\end{proposition}

\begin{proof}
Let us fix an $\varepsilon >0$. 
For the first inequality we apply \cref{lem:Construction_DCirc} and for the second we apply \cref{lem:Construction_DNatural};
recall that  
\begin{gather*}
\int_{(0,\infty)} 
  \textup{Tr}\bigg[ (Z^{\infty,(p)}_s - \widetilde{Z}^{\varepsilon,(p)}_s)\frac{\d \langle X^{\infty,\circ}\rangle_s}{\d C^\infty_s}(Z^{\infty,(p)}_s - \widetilde{Z}^{\varepsilon,(p)}_s)^\top \bigg]
  \d C^\infty_s
  =\int_{(0,\infty)} 
\big\Vert   c^\infty_s(\widetilde{Z}^{\varepsilon,(p)}_s - Z_s^{\infty,(p)})\big\Vert^2_2 
\d C^\infty_s,
\shortintertext{and}
  \int_{(0,\infty)} \Vert U^{\infty,(p)}(s,x) - \widetilde{U}^{\varepsilon,(p)}(s,x)\Vert_2^2 \nu^{\infty,\natural}(\d s,\d x)
  =  
  \int_{(0,\infty)} \tnorm{U^{\infty,(p)}(s,\cdot) - \widetilde{U}^{\varepsilon,(p)}(s,\cdot)}_{\infty,s}^2 \d C^{\infty}_s.
\end{gather*}
Finally, we can prove the validity of \eqref{ineq:LusinApproximation} in view of the following set inclusions\footnote{We regard the following sets as subsets of $\Omega_{\text{sub}}$.}
\begin{align*}
&\biggl\{
  \bigg\Vert  \int_{(0,\infty)} f^\infty\bigl(s,Y_s^{\infty,(p)},\widetilde{Z}^{\varepsilon,(p)}_s,\widetilde{U}^{\varepsilon,(p)}(s,\cdot)\bigr)
  \d C^\infty_s - 
  L^{\infty,(p)}_\infty\bigg\Vert_2 \ge\frac{\varepsilon}{3} \biggr\}\\
&=\biggl\{
  \bigg\Vert  \int_{(0,\infty)} f^\infty\bigl(s,Y_s^{\infty,(p)},\widetilde{Z}^{\varepsilon,(p)}_s,\widetilde{U}^{\varepsilon,(p)}(s,\cdot)\bigr)
  \d C^\infty_s
  -L^{\infty,(p)}_\infty\bigg\Vert^2_2 \ge \frac{\varepsilon^2}{9} \biggr\}\\
&\subseteq
\biggl\{\frac{1}{\hat{\beta}}  \int_{(0,\infty)} \textup{e}^{\hat{\beta} A^\infty_s}
\Bigl[ \big\Vert   c^\infty_s(\widetilde{Z}^{\varepsilon,(p)}_s - Z_s^{\infty,(p)})\big\Vert^2_2 
+\tnormbig{\widetilde{U}^{\varepsilon,(p)}(s,\cdot) - U^{\infty,(p)}(s,\cdot)}_{\infty,s}^2\Bigr]
\d C^\infty_s \ge\frac{\varepsilon^2}{9} \biggr\}\footnotemark\\
&
\subseteq
\biggl\{\frac{\textup{e}^{\hat{\beta} \overline{A}}}{\hat{\beta}}  
\int_{(0,\infty)} 
\big\Vert   c^\infty_s(\widetilde{Z}^{\varepsilon,(p)}_s - Z_s^{\infty,(p)})\big\Vert^2_2 
\d C^\infty_s \ge\frac{\varepsilon^2}{18} \biggr\} \\
&\qquad \qquad \qquad \qquad \qquad 
\bigcup \biggl\{\frac{\textup{e}^{\hat{\beta} \overline{A}}}{\hat{\beta}}  
\int_{(0,\infty)} 
\tnormbig{\widetilde{U}^{\varepsilon,(p)}(s,\cdot) - U^{\infty,(p)}(s,\cdot)}_{\infty,s}^2
\d C^\infty_s \ge\frac{\varepsilon^2}{18} \biggr\}
\subseteq\emptyset. \qedhere
\end{align*}
\footnotetext{We have used the fact that $\theta^{\infty,c},\theta^{\infty,d}\le \alpha^\infty$ (see \ref{Ffour}) as well as \cite[Inequality 3.16]{papapantoleon2016existence}.}
\end{proof}

\subsubsection{The claim \texorpdfstring{\eqref{limsup1}}{(3.22)} is true}\label{subsubsec:limsup1}
\begin{lemma}\label{lem:limsup1}
It holds
\begin{equation*}
  \limsup_{i\to\infty} \bigg\{\omega\in\widetilde{\Omega}:\bigg\|
L^{i,(p)}_\infty(\omega) - 
\bigg(\int_{(0,\infty)} 
     f^{i}\bigl(s,Y_s^{{i},(p)},\widetilde{Z}^{\varepsilon,(p)}_s,\widetilde{U}^{\varepsilon,(p)}(s,\cdot)\bigr)\d C^{i}_s\bigg)(\omega)\bigg\|_2
     \ge \frac{\varepsilon}{3}\bigg\}=\emptyset.
 \end{equation*}
\end{lemma}

\begin{proof}
Using the Cauchy--Schwarz's inequality as in \cite[Inequality 3.16]{papapantoleon2016existence}, Condition \ref{BIConditionAk}, the Lipschitz property of the generator and Inequalities \eqref{prop:Lipschitz-bounds}, we have 
\begin{align*}
&\biggl\|
  L^{i,(p)}_\infty - 
  \int_{(0,\infty)} 
  f^{i}\bigl(s,Y_s^{{i},(p)},\widetilde{Z}^{\varepsilon,(p)}_s,\widetilde{U}^{\varepsilon,(p)}(s,\cdot)\bigr)\d C^{i}_s
\bigg\|_2^2\\
&\leq
\frac{\textup{e}^{\hat{\beta}\overline{A}}}{\hat{\beta}}
\int_{(0,\infty)} 
  \big\Vert (Z^{i,(p)}_s - \widetilde{Z}^{\varepsilon,(p)}_s) c^i_s\big\Vert_2^2 
  + \tnorm{U^{i,(p)}_s(\cdot)- \widetilde{U}^{\varepsilon,(p)}(s,\cdot)}^2_{i,s} 
\d C^{i}_s\\
& =
\frac{\textup{e}^{\hat{\beta}\overline{A}}}{\hat{\beta}}
\int_{(0,\infty)} 
  \textup{Tr}\big[ (Z^{i,(p)}_s - \widetilde{Z}^{\varepsilon,(p)}_s) \frac{\d  \langle X^{i,\circ}\rangle_s}{\d  C^i_s}(Z^{i,(p)}_s - \widetilde{Z}^{\varepsilon,(p)}_s)^\top\big] \d C^i_s \\
&  \quad +
\frac{\textup{e}^{\hat{\beta}\overline{A}}}{\hat{\beta}}
\int_{(0,\infty)} 
 \tnorm{U^{i,(p)}_s(\cdot)- \widetilde{U}^{\varepsilon,(p)}(s,\cdot)}^2_{i,s}\d C^{i}_s.
\numberthis \label{eq:last-convergence}
\end{align*}

In order to argue about the validity of the aforementioned convergence, we initially return to Convergence \eqref{conv:Memin-Cor}, on which we will apply the usual measure-theoretic arguments in order to derive
for every $Z\in D^{\circ,\ell\times m}$
\begin{align}
\label{conv:ZRNZ-compact-support}
  \int_{(0,\infty)} Z_s\frac{\d \langle X^{i,\circ}\rangle_s}{\d C^i_s} Z_s \d C^i_s
  \xrightarrow[i\to\infty]{\hspace{1em}\vert \cdot\vert\hspace{1em}}
  \int_{(0,\infty)} Z_s\frac{\d \langle X^{\infty,\circ}\rangle_s}{\d C^\infty_s} Z_s \d C^\infty_s,
  \; \mathbb{P}\text{\rm--a.s.}
\end{align}
Let us fix $(p,q) \in \{1,\dots, m\}^2$.
The first step towards proving our aim is to assume an open interval $I\subset \mathbb{R}_+$.
Then, from Convergence \eqref{conv:Memin-Cor} we have for the $pq$-elements of $(\langle X^{i,\circ} \rangle)_{i\in\overline{\mathbb{N}}}$
\begin{align*}
  \int_{(0,\infty)} \mathds{1}_I(s)\d  \langle X^{i,\circ,p},X^{i,\circ,q}\rangle_s
  \xrightarrow[i\to\infty]{\hspace{1em}\vert \cdot\vert\hspace{1em}}
  \int_{(0,\infty)} \mathds{1}_I(s)\d  \langle X^{i,\circ,p},X^{i,\circ,q}\rangle_s,\; \mathbb{P}\text{\rm--a.s.}
\end{align*}
The reader may recall that $\langle X^{\infty,\circ} \rangle$ is a continuous process, so no doubts are raised on points of discontinuity because they do not exist.
The next step is to assume a step function, \emph{i.e.}, a finite family of real numbers $(\alpha^{k})_{k\in\{1,\dots,n\}}$ and a finite family of disjoint open intervals $(I^{k})_{k\in\{1,\dots,n\}}$, for which 
\begin{align*}
  \int_{(0,\infty)} \sum_{k=1}^n\alpha^{k}\mathds{1}_I^{k}(s)\d  \langle X^{i,\circ,p},X^{i,\circ,q}\rangle_s
  \xrightarrow[i\to\infty]{\hspace{1em}\vert \cdot\vert\hspace{1em}}
  \int_{(0,\infty)} \sum_{k=1}^n\alpha^{k}\mathds{1}_I^{k}(s)\d  \langle X^{i,\circ,p},X^{i,\circ,q}\rangle_s,\; \mathbb{P}\text{\rm--a.s.},
\end{align*}
is immediate from the previous convergence.
The third step is to assume a function $Z$ with compact support which has discontinuities of the first kind, \emph{i.e.}, the left- and right-limits exist.
Since every such function is the uniform limit of step-functions, we can conclude from the second step the validity of
\begin{align*}
  \int_{(0,\infty)} Z(s)\d  \langle X^{i,\circ,p},X^{i,\circ,q}\rangle_s
  \xrightarrow[i\to\infty]{\hspace{1em}\vert \cdot\vert\hspace{1em}}
  \int_{(0,\infty)} Z(s)\d  \langle X^{i,\circ,p},X^{i,\circ,q}\rangle_s,\; \mathbb{P}\text{\rm--a.s.}
\end{align*}
Finally, we recall that $D^{\circ,\ell\times m}$ is countable and consists of continuous functions with compact support.
Hence, from the last step we can immediately derive (by summing suitably the respective elements) Convergence \eqref{conv:ZRNZ-compact-support}. Next, from Convergence \eqref{conv:g-compensator} for $U\in D^{\natural}$
\begin{gather*}
\begin{multlined}[c][0.85\displaywidth]
  \big[U^\top U\big]* \nu^{i,\natural}_\infty
    - \sum_{s\ge 0} \int_{\mathbb{R}^n} \big(U(s,x)\big)^\top \nu^{i,\natural}(\{s\}\times\d x) 
                    \int_{\mathbb{R}^n} U(s,x) \nu^{i,\natural}(\{s\}\times\d x) 
  \xrightarrow[i\to\infty]{\hspace{1em}\vert \cdot\vert\hspace{1em}}
  [U^\top U] *\nu^{\infty,\natural}_\infty,\; \mathbb{P}\text{\rm--a.s.}
\end{multlined}
\end{gather*}
Returning to  \eqref{Conv:Angle-Brack-kp-Pas-2}, and using again the usual measure-theoretic arguments, we have in particular for $Z\in D^{\circ,\ell\times  m}$\footnote{We use the fact that $Z$ and $U$ have compact supports and the pointwise convergence on the continuity points of the limit (see \cite[Proposition VI.2.1]{jacod2003limit}, in order to attain the value of the integral over $\mathbb{R}_+$ on a compact subinterval.}
\begin{gather*}
\int_{(0,\infty)} Z^{i,(p)}_s\frac{\d \langle X^{i,\circ}\rangle_s}{\d C^i_s} Z_s \d C^i_s
\xrightarrow[i\to\infty]{\hspace{1em}\vert \cdot\vert\hspace{1em}}
\int_{(0,\infty)} Z^{\infty,(p)}_s\frac{\d \langle X^{\infty,\circ}\rangle_s}{\d C^\infty_s} Z_s \d C^\infty_s,\\
\shortintertext{}
\begin{multlined}[c][0.85\displaywidth]
\big[\big(U^{i,(p)}\big)^\top U\big] *\nu^{i,\natural}_\infty
  - \sum_{s\ge 0}
\int_{\mathbb{R}^n} \big(U^{i,(p)}(s,x)\big)^\top \nu^{i,\natural}(\{s\}\times\d x) 
\int_{\mathbb{R}^n} U(x) \nu^{i,\natural}(\{s\}\times\d x) \\
\xrightarrow[i\to\infty]{\hspace{1em}\vert \cdot\vert\hspace{1em}}
\big[\big(U^{\infty,(p)}\big)^\top U\big]* \nu^{\infty,\natural}_\infty.
\end{multlined}
\end{gather*}
We return to \eqref{Conv:Angle-Brack-kp-Pas} and at this point we use Condition \ref{BFiniteStoppingTime}, which implies\footnote{The reader may recall \cref{rem:Conv-Subseq-klm}.\ref{rem:Conv-Subseq-klm-stopping-time}.}
\begin{gather*}
\int_{(0,\infty)} Z^{i,(p)}_s\frac{\d \langle X^{i,\circ}\rangle_s}{\d C^i_s} (Z^{i,(p)}_s)^\top \d C^i_s
\xrightarrow[i\to\infty]{\hspace{1em}\vert \cdot\vert\hspace{1em}}
\int_{(0,\infty)} Z^{\infty,(p)}_s\frac{\d \langle X^{\infty,\circ}\rangle_s}{\d C^\infty_s} (Z^{\infty,(p)}_s)^{\top} \d C^\infty_s,\\
\shortintertext{}
\begin{multlined}[c][0.85\displaywidth]
\langle U^{i,(p)}\star \widetilde{\mu}^{{i,\infty} }\rangle_\infty
\xrightarrow[i\to\infty]{\hspace{1em}\vert \cdot\vert\hspace{1em}}
\langle U^{\infty,(p)}\star \widetilde{\mu}^{\infty,\infty}\rangle_\infty.
\end{multlined}
\end{gather*}
Overall, from all the above convergence and Inequality \eqref{eq:last-convergence} we derive the desired statement.
\end{proof}
\subsubsection{The claim \texorpdfstring{\eqref{limsup2}}{(3.23)} is true}
We will apply \cref{prop:MOAppl}.
In order to proceed further, we will provide the following  preparatory results.
Recall that in \cref{rem:Omega-sub} we have fixed an $\Omega_\textup{sub}$ with $\mathbb P[\Omega_\textup{sub}]=1.$ 
%
%
%
%
%
For the remainder of the section we will fix an arbitrarily small $\varepsilon>0$; we will assume that $\varepsilon<1$.
Moreover, $\widetilde{Z}^{\varepsilon,(p)}$, $\widetilde{U}^{\varepsilon,(p)}$ will be assumed to satisfy the properties described in \cref{Lusin_approximation} and will be hereinafter fixed.
%
%
%
%
%
%
%
\begin{lemma}\label{lem:uniform_bound_Yip}
For every $\omega\in\Omega_{\text{\emph{sub}}}$, we have $\sup_{i\in\overline{\mathbb{N}}}\Vert Y^{i,(p)}(\omega)\Vert_\infty <\infty$.
\end{lemma}

\begin{proof}
We will apply \cref{UniformlyBounded}.
To this end, we will ensure that the three required conditions are satisfied. 
The first, resp. second, condition is indeed true in view of \ref{rem:Conv-Subseq-klm-4-a}, resp. \ref{rem:Conv-Subseq-klm-6}, of \cref{rem:Conv-Subseq-klm}.
Regarding the validity of the last condition, we can argue as follows
\begin{align*}
  \mathbb E\biggl[\limsup_{(i,t)\to(\infty,\infty)} \big\|Y^{i,(p)}_t\big\|_2 \biggr] 
    &=  \mathbb E\biggl[\lim_{j\to\infty} \big\|Y^{i_j,(p)}_{t_j}\big\|_2 \biggr]\footnotemark 
    \le \mathbb E\biggl[\lim_{j\to\infty} \sup_{t\in\mathbb R_+}\big\|Y^{i_j,(p)}_{t}\big\|_2 \biggr] \\
    &\le \liminf_{j\to\infty} \mathbb E\bigg[ \sup_{t\in\mathbb R_+}\big\|Y^{i_j,(p)}_{t}\big\|_2 \bigg]
    \le \sup_{i\in\overline{\mathbb N}} \mathbb E\bigg[ \sup_{t\in\mathbb R_+}\big\|Y^{i,(p)}_{t}\big\|_2 \bigg]<\infty,
\end{align*}%
\footnotetext{The equality holds for some subsequence $\bigl((i_j,t_j)\bigr)_{j\in\mathbb N}$, which depends on the path, such that $(i_j,t_j)\to(\infty,\infty)$.}%
where in the last step we used \cref{rem:Conv-Subseq-klm}.\ref{rem:Conv-Subseq-klm-4-a} again.  
Therefore, $\mathbb P\big[\big\{ \limsup_{(i,t)\to(\infty,\infty)} \big\|Y^{i,(p)}_t\big\|_2  = \infty\bigr\}\bigr]=0$, which allows us to conclude that 
\begin{equation*}
  \limsup_{(i,t)\to(\infty,\infty)} \|Y^{i,(p)}(\omega)\|_\infty<\infty, \text{ for every }\omega\in\Omega_{\text{sub}}. \qedhere
\end{equation*}
\end{proof}

\begin{proposition}\label{prop:p-Step-WeakApprox}
For every $\omega\in\Omega_{\text{\emph{sub}}}$  holds
\begin{gather*}
\int_{(0,\cdot]} \hspace{-0.2em}f^i\bigl(s,Y_s^{i,(p)},\widetilde{Z}^{\varepsilon,(p)}_s,\widetilde{U}^{\varepsilon,(p)}(s,\cdot)\bigr)\d C^i_s
\xrightarrow[i\to\infty]{\hspace{0.2em}\textup{J}_1\hspace{0.2em}}
\int_{(0,\cdot]}\hspace{-0.2em} f^\infty\bigl(s,Y^{\infty,(p)}_s,\widetilde{Z}^{\varepsilon,(p)}_s,\widetilde{U}^{\varepsilon,(p)}(s,\cdot)\bigr)\d C^\infty_s ,
\shortintertext{}
\int_{(0,\infty)} \hspace{-0.2em}f^i\bigl(s,Y_s^{i,(p)},\widetilde{Z}^{\varepsilon,(p)}_s,\widetilde{U}^{\varepsilon,(p)}(s,\cdot)\bigr)\d C^i_s \xrightarrow[i\to\infty]{\hspace{0.2em}\Vert \cdot \Vert_2\hspace{0.2em}}
\int_{(0,\infty)}\hspace{-0.2em} f^\infty\bigl(s,Y^{\infty,(p)}_s,\widetilde{Z}^{\varepsilon,(p)}_s,\widetilde{U}^{\varepsilon,(p)}(s,\cdot)\bigr)\d C^\infty_s.
\end{gather*}
\end{proposition}

\begin{proof}
We will apply \cref{cor:MOAppl} for the former convergence and \cref{prop:MOAppl} for the latter.
For both of these, the required assumptions are satisfied
\begin{itemize}
  \item by \cref{rem:Conv-Subseq-klm} we have the weak convergence for the measures associated to the sequence $(C^i)_{i\in\overline{\mathbb{N}}}$;
  \item by \ref{Bgenerator} we have that the integrands converge in the Skorokhod topology;
  \item \cref{lem:uniform_bound_Yip} in conjunction with \ref{Bgenerator}.\ref{Bgenerator:respects-J1-conv} imply that the integrands are uniformly bounded.\qedhere
\end{itemize}
\end{proof}
\subsubsection{The claim \texorpdfstring{\eqref{LSp-Pas}}{(3.19)} is true}
Since the limit process is continuous, the metric $\delta_{\textup{J}_1}$ is identical to the metric induced by the locally uniform convergence; see \citeauthor{jacod2003limit} \cite[Proposition VI.1.17]{jacod2003limit}. 
Hence, we need only to prove that for every $N\in\mathbb{N}$
\begin{align*}
\sup_{t\in [0,N]} \big\Vert L^{i,(p)}_t - L^{\infty,(p)}_t\big\Vert_2 \xrightarrow[i\to\infty]{}0.
\end{align*}
On the other hand, a closer examination of the properties we have presented in the previous sub-subsections leads us to the required conclusion. 
Indeed, starting from the choice of $\widetilde{Z}^{\varepsilon,(p)}$ and $\widetilde{U}^{\varepsilon,(p)}$ in \cref{Lusin_approximation} 
we have for every $N\in \mathbb N$
\begin{gather*}
\begin{multlined}[\textwidth]
\bigg\Vert  \int_{(0,N]} f^\infty\bigl(s,Y^{\infty,(p)}_s,\widetilde{Z}^{\varepsilon,(p)}_s,\widetilde{U}^{\varepsilon,(p)}(s,\cdot)\bigr)\d C^\infty_s - L^{\infty,(p)}_N\bigg\Vert_2^{2}\\
\leq
\frac{\textup{e}^{\hat{\beta} \overline{A}}}{\hat{\beta}}  
\int_{(0,\infty)} 
\big\Vert   c^\infty_s(\widetilde{Z}^{\varepsilon,(p)}_s - Z_s^{\infty,(p)})\big\Vert^2_2 
+\tnormbig{\widetilde{U}^{\varepsilon,(p)}(s,\cdot) - U^{\infty,(p)}(s,\cdot)}_{\infty,s}^2
\d C^\infty_s 
 <\frac{\varepsilon}{3},\; \mathbb P\text{\rm--a.s.}
\end{multlined}
\end{gather*}
The analogous to Convergence \eqref{limsup1}, resp. \eqref{limsup2}, can be proven in view of the comments in \cref{rem:Conv-Subseq-klm} and using the usual measure theoretic arguments, resp. has been proved in \cref{prop:p-Step-WeakApprox}.

%% file: Stability_BSDEJ_Induction_pStep_UI.tex

\subsubsection{Uniform integrability} 

In this sub-subsection, we will prove that statement \eqref{Gamma-p-UI} is true,
\emph{i.e.} $\big( \Gamma^{k,(p)})_{k\in\overline{\mathbb{N}}}$ is uniformly integrable.
This property will simultaneously ensure that the sequence $\big( \sup_{t\in\mathbb R_+}\big|L^{k,(p)}_t\big|^2\big)_{k\in\overline{\mathbb{N}}}$ is uniformly integrable and the validity of the induction step.
Recall that we need the uniform integrability of  $\big( \sup_{t\in\mathbb R_+}\big|L^{k,(p)}_t\big|^2\big)_{k\in\overline{\mathbb{N}}}$ in order to apply Vitali's theorem and prove that Convergences \eqref{LSp} and \eqref{LSpInfinity} hold in $\mathbb L^2$-mean and not only in probability, which we proved in the previous sub-subsections.

\begin{lemma}
The sequence $\big( \Gamma^{k,(p)})_{k\in\overline{\mathbb{N}}}$ is uniformly integrable.
\end{lemma}
\begin{proof}
We essentially follow the same steps as in \cref{lem:UI-Gamma} in order to derive for every $k\in\mathbb{N}$
\begin{align*}
\Gamma^{k,(p)}
	& \le
2\overline{A} \sup_{t\in\mathbb{R}_+}\bigl\Vert  \mathbb{E}[ \xi^k |\mathcal{G}^k_t]\bigr\Vert_2^2
+2\overline{A} \sup_{t\in\mathbb{R}_+}\Bigl\Vert  \mathbb{E}\Big[
L^{k,(p-1)}_\infty-L^{k,(p-1)}_t
 \Big|\mathcal{G}^k_t\Big] \Bigr\Vert_2^2
+ \textup{Tr}\big[\langle M^{k,(p-1)}\rangle_\infty\big]
+ \Gamma^{k,(0)}.
\end{align*}
We only need to observe now that the summands of the right hand side form a uniformly integrable sequence, for $k\in\overline{\mathbb{N}}$, as sum of elements of uniformly integrable random variables; see \citeauthor{he1992semimartingale} \cite[Corollary 1.10]{he1992semimartingale}. 
Indeed
\begin{itemize}[itemindent=0.cm,leftmargin=*]
	\item for the sequence associated to the first summand, we have that $(\Vert \xi^k \Vert_2^2)_{k\in\overline{\mathbb{N}}}$ is uniformly integrable by Vitali's theorem.
	Then, we can conclude the uniform integrability of the required sequence by \cref{cor:DoobApplication};
	
	\smallskip
	\item for the sequence associated to the second summand, we can prove by means of the conditional Jensen Inequality and the conditional Cauchy--Schwarz's inequality that for all $k\in\overline{\mathbb{N}}$ holds
	\begin{align*}
	\big\Vert  \mathbb{E}\big[  L^{k,(p-1)}_\infty - L_t^{k,(p-1)} \big| \mathcal{G}_t^k\big]\big\Vert_2^2
	\le  2\mathbb{E}\big[ \Gamma^{k,(p-1)} \big| \mathcal{G}^k_t\big] . 
	\end{align*}
	Since $(\Gamma^{k,(p-1)})_{k\in\overline{\mathbb{N}}}$ is uniformly integrable (by the induction assumption) we can conclude the uniform integrability of  $\big(  \big\Vert  \mathbb{E}\big[  L^{k,(p-1)}_\infty - L_t^{k,(p-1)} \big| \mathcal{G}_t^k\big]\big\Vert_2^2\big)_{k\in\overline{\mathbb{N}}}$ by \citeauthor{he1992semimartingale} \cite[Theorem 1.7]{he1992semimartingale} and \cref{cor:DoobApplication};
	
\smallskip
	\item for the sequence  associated to the third summand we use the fact that (which is true in view of the validity of \eqref{LSpminus1})
	\begin{align*}
		M^{k,(p-1)}_\infty \xrightarrow[k\to\infty]{\hspace{1em} \mathbb{L}^2(\Omega,\mathcal{G}, \mathbb{P})\hspace{1em}}  M^{\infty,(p-1)}_\infty,
	\end{align*}
	which implies the uniform integrability of 
	$\big(\text{Tr}\big[\langle M^{k,(p-1)}\rangle_\infty\big]\big)_{k\in\overline{\mathbb{N}}}$%
	\footnote{The uniform integrability of $\big(\text{Tr}\big[[M^{k,(0)}]_\infty\big]\big)_{k\in\overline{\mathbb{N}}}$ can also be deduced.}; recall also \cref{rem:Conv-Subseq-klm}.\ref{UI:Delta_p};
	
	\smallskip
	\item finally, the sequence associated to the last summand is uniformly integrable by \ref{BgeneratorUI}.\qedhere
\end{itemize}

\end{proof}
We can now deduce the following.
\begin{lemma}
The sequence $\Big(\sup_{t\in\mathbb R_+}\big\|L^{k,(p)}_t\big\|_2^2\Big)_{k\in\overline{\mathbb{N}}}$ is uniformly integrable.
\end{lemma}
\begin{proof}
Apply Cauchy--Schwarz's inequality as in \cite[Inequality 3.16]{papapantoleon2016existence} and use the previous lemma.
\end{proof}

%% file: Stability_BSDEJ_Comments_on_Framework.tex

In this subsection, we would like to discuss the nature of the conditions we have imposed in order to set the framework for \cref{BSDERobMainTheorem}.

\medskip
Let us start with Conditions \ref{MFilqlc}--\ref{Mfinalrv}, which are required for \citeauthor{papapantoleon2019stability} \cite[Corollary 3.10]{papapantoleon2019stability}. 
We have repeatedly stated (and it should be clear from the outline of the proof presented in Subsection \ref{subsec:BodyProof}) that the stability of martingale representations plays a crucial role here.  
However, one may wonder about the necessity of the convergence imposed in \ref{MSIprime}
\begin{equation*}
\widebar X^k_\infty
\xrightarrow[k\to\infty]{\hspace{0.2em}\mathbb L^2(\mathcal G;\mathbb R^{m+n})\hspace{0.2em}}
\widebar X^\infty_\infty,
\end{equation*}
which is stronger compared to \cite[Condition {(M2)}]{papapantoleon2019stability}; the latter, assuming that $m=n$, reads as
\begin{equation*}
X^{k,\circ}_\infty + X^{k,\natural}_\infty
\xrightarrow[k\to\infty]{\hspace{0.2em}\mathbb L^2(\mathcal G;\mathbb R^n)\hspace{0.2em}}
X^{\infty,c}_\infty + X^{\infty,d}_\infty.
\end{equation*}
The answer has two parts. 
Initially, we would like to allow the driving martingales $X^{k,\circ}$ and $X^{k,\natural}$ to have different dimensions, for $k\in\overline{\mathbb{N}}$. 
Afterwards, in sub-subsection \ref{subsubsec:limsup1}, we need to use the convergence
\begin{equation*}
 \bigl(\langle X^{i,\circ}\rangle, \langle X^{i,\natural} \rangle \bigr)
 	\xrightarrow[i\to\infty]{\hspace{0.2em}\textup{J}_1(\mathbb R^{m\times m}\times\mathbb R^{n\times n})\hspace{0.2em}}
 \bigl(\langle X^{\infty,c}\rangle, \langle X^{\infty,d} \rangle \bigr),\; \mathbb P\text{\rm--a.s.},
\end{equation*} 
as well as the respective convergences from \cref{rem:Conv-Subseq-klm}.\ref{rem:Conv-Subseq-klm-4}.
These convergences are not guaranteed by \cite[Theorem 3.3]{papapantoleon2019stability}.
Nevertheless, we still have the flexibility to approximate the continuous martingale part of $\oX^\infty$\footnote{Recall that $X^{\infty,\circ}\in\mathcal{H}^{2,c}(\mathbb{G}^{\infty},\tau^{\infty};\mathbb{R}^{m})$.} with purely discontinuous martingales, a property which is essential when the discussion comes to numerical schemes.

\medskip
Regarding Condition \ref{BIConditionAk}, we have already provided some comments in \cref{rem:BSDE-solution}. 
Let us  provide a more technical remark at this point. 
The aforementioned condition comes essentially into play for proving the uniform integrability of the sequence $(\Gamma^{k,(p)})_{k\in\overline{\mathbb{N}}}$, for every $p\in\mathbb{N}$. 
There, we need in particular to prove the uniform integrability of the sequence $\big(\int_{(0,\infty)} \Vert Y^{k,(p)}_s \Vert_2^{2} \textup{d}A^{k}_s\big)_{k\in\overline{\mathbb{N}}}$. 
From \cref{cor:DoobApplication} we can prove the uniform integrability of $(\sup_{t\in\mathbb{R}_+} \Vert Y^{k,(p)}_t \Vert_2^{2})_{k\in\overline{\mathbb{N}}}$, which is equivalent to the desired property, once Condition \ref{BIConditionAk} is enforced. 
In a more general situation, it is not clear if, or how, the desired property can be obtained.

\medskip
Let us proceed now to Condition \ref{Bgenerator}.
The first part of this condition provides a regularity property for the paths of the generators under relatively weak assumptions and is in accordance with \citeauthor*{briand2002robustness} \cite[Condition (H3) (ii)]{briand2002robustness}, which considers the case where $X^{\infty,\circ}$ is the Brownian motion and $X^{\infty,\natural}=0$. 
One could have assumed a weaker regularity property but, as a trade-off, the convergence in the second part should be strengthened, \emph{e.g.}, to pathwise convergence under the supremum norm. 
For the second part of the condition, similar conditions are considered when we need to guarantee the convergence of compositions under the $\textup{J}_1$-topology, \emph{e.g.} see  \citeauthor*{kurtz1991weak} \cite[Lemma 2.1]{kurtz1991weak}.
In the case where $(f^{k})_{k\in\overline{\mathbb{N}}}$ are globally Lipschitz with the same Lipschitz constant, one need not be as abstract as in \cite{kurtz1991weak} since one can exploit the Lipschitz property of the generators.
For example, in \citeauthor*{briand2002robustness} \cite[Proposition 11]{briand2002robustness} under the aforementioned uniform equi-continuity assumption on the generators, the uniform (in time) convergence \emph{for every fixed point $(y,z)$}, in our notation,
(see \cite[Convergence (7)]{briand2002robustness}) is equivalent to the uniform (in time) convergence \emph{on compacts}.
Hence, Condition \ref{Bgenerator}.\ref{Bgenerator:respects-J1-conv} is a relaxation of \cite[Convergence (7)]{briand2002robustness}. 
%
%
%
%
%
%
%
%
%
%
%

\medskip
Let us further proceed to the discussion about Condition \ref{Bintegrator}, which we will analyse part by part.
\begin{enumerate}[label=\textup{{(\roman*)}}, itemindent=0.cm, leftmargin=*]
	\item 
The pair $(\oX^k,C^k)$ satisfies \cite[Assumption 2.10]{papapantoleon2016existence}, for each $k\in\overline{\mathbb{N}}$.
The existence of the integrator $C^k$ is ensured by arguments analogous to \citeauthor{jacod2003limit} \cite[Proposition II.2.9]{jacod2003limit}.
In other words, we have set \[C^k:=h^{k}\big(\langle X^{k,\circ}\rangle,\textup{Id}*\nu^{(X^\natural,\mathbb G^k)}\big),\] for some function $h^{k}$, for every $k\in\overline{\mathbb{N}}$.
On the other hand, in view of Conditions \ref{MFilqlc}--\ref{Mfinalrv}, we know by \cite[Corollary 3.10]{papapantoleon2019stability} that 
\begin{align*}
	\langle \oX^{k}\rangle \xrightarrow[k\to\infty]{(\textup{J}_1(\mathbb{R}^{(m+n)\times(m+n)}),\mathbb{L}^{1})}\langle \oX^{\infty} \rangle.
\end{align*}
Hence, we can choose the sequence $(h^k)_{k\in\mathbb N}$, which may allow dependence on time, as soon as it respects the $\textup{J}_1$-topology. 
In particular, we may choose $h^{k}=h$ for every $k\in\overline{\mathbb{N}}$, where $h:(\mathbb{D}(\mathbb{R}^{}),\delta_{\textup{J}_1})\longrightarrow(\mathbb{R}_+,|\cdot|)$ is a continuous function.
\medskip
\item This is almost analogous to the previous point. 
Since $\langle \oX^{k}\rangle_\infty\in\mathbb{L}^{1}(\mathcal{G})$, for every $k\in\overline{\mathbb{N}}$, the finiteness of $C^{k}_\infty$, for every $k\in\overline{\mathbb{N}}$, is immediate by the comments above.
However, there is a subtle point in the convergence part that we will take care of with the help of Condition \ref{BFiniteStoppingTime}.
In general, we \emph{cannot} prove that 
$\langle \oX^{k} \rangle_\infty \xrightarrow{\hspace{0.2em}\mathbb L^1(\mathcal{G})\hspace{0.2em}} \langle \oX^{\mage \infty } \rangle_\infty$, although we \emph{can} prove that  
$[\oX^{k} ]_\infty \xrightarrow{\hspace{0.2em}\mathbb L^1(\mathcal{G})\hspace{0.2em}} [\oX^{\mage \infty}]_\infty$; see \citeauthor{memin2003stability} \cite[Proof of Corollary 12]{memin2003stability} for the one dimensional case.
However, the process $\langle \oX^{\infty} \rangle$ is continuous.
Hence, for a subsequence $(k_l)_{l\in\mathbb{N}}$ such that 
\begin{align*}
	\langle \oX^{k_l}\rangle \xrightarrow[l\to\infty]{\textup{J}_1(\mathbb{R}^{(m+n)\times(m+n)})}\langle \oX^{\infty} \rangle, \; \mathbb{P}\text{\rm--a.s.},
\end{align*}
it holds for every $t\in \mathbb{R}_+$  
\begin{align*}
	\langle \oX^{k_l}\rangle_t \xrightarrow[l\to\infty]{\Vert\cdot\Vert_2}\langle \oX^{\infty} \rangle_t,\; \mathbb{P}\text{\rm--a.s.}
\end{align*}
Hence, if Condition \ref{BFiniteStoppingTime} is enforced, we can derive for the diagonal elements of the predictable quadratic covariation processes (which are increasing) that 
\begin{align}\label{conv:AngleBracketInfty}
	\langle \oX^{k_{l_i},jj}\rangle_{\infty}=\langle \oX^{k_{l_i},jj}\rangle_{\tau^k} \xrightarrow[i\to\infty]{|\cdot|}\langle \oX^{\infty,jj} \rangle_{\tau^{k}}=\langle \oX^{\infty,jj} \rangle_\infty,\; \mathbb{P}\text{\rm--a.s.}
\end{align}
By the polarisation identity, we can prove the convergence for all the off-diagonal elements as well.
Hence, the assumption $C^{k}_\infty \xrightarrow[k\to\infty]{\mathbb{P}} C^{\infty}_\infty$ is reasonable.
\item In view of \ref{BIConditionAk} and the previous comments, this assumption is completely natural.
\end{enumerate}\smallskip

Finally, let us comment on Condition \ref{BFiniteStoppingTime}, which was used only for the proof\footnote{Actually, we have not used Condition \ref{BFiniteStoppingTime} when an alternative valid argument could have been used, \emph{e.g.}, the use of the compact supports of integrands $Z$ and $U$.} of \cref{lem:limsup1} and for arguing about the reasonableness of Condition \ref{Bintegrator}.{(ii)}.  
The finiteness of $\tau^{\infty}$ seems unavoidable given our need for Convergence \eqref{conv:AngleBracketInfty}. 
Regarding the convergence of the stopping times, on the one hand it is trivially satisfied when the terminal times are deterministic and not necessarily identical.
This case corresponds to the majority of the articles in the literature. 
On the other hand, when the terminal times are random, the required convergence can be proved for specific debut times of a convergent sequence of processes; see the proof of \citeauthor{jacod2003limit} \cite[Theorem IX.1.17]{jacod2003limit}.
However, there is no general way (to the best of our knowledge) for constructing convergent sequences of stopping times when the associated filtrations change. 
Answers to this problem in special cases may be offered by \citeauthor{coquet2007convergence} \cite[Proposition 20]{coquet2007convergence} or \citeauthor{kchia2011semimartingales} \cite[Lemma 20]{kchia2011semimartingales}.
As a last comment, the convergence we assumed does not require any integrability condition on the stopping times; as a comparison in \citeauthor{toldo2006stability} \cite{toldo2006stability} a uniform integrability condition is required.

%% file: Examples.tex

In this section, we present a few examples and applications that demonstrate the power of our framework.
We assume that the sequence of processes $(\overline{X}^k)$ for $k\in\mathbb{N}$ are indexed either in a discrete- or in a continuous-time set; in case the processes are indexed in a discrete-time set, they can be embedded in our framework in the obvious way.
In other words, the processes $(\overline{X}^k)_{k\in\mathbb{N}}$ could be, for example, random walks or discretisations or even other continuous-time approximations of the driving martingales.
In case of random walks, which is the natural choice when considering numerical schemes for BSDEs, each martingale will be defined with respect to its own filtration, and the convergence of these filtrations is a natural requirement.
Otherwise, our framework allows great flexibility in choosing a convergent sequence of suitable standard data. \medskip

Before we proceed, let us remind the reader that the Skorokhod space endowed with the $\textrm{J}_1$-topology is Polish.
Consequently, we can use Skorokhod's representation theorem in order to assume, without loss of generality, that the processes converge in probability, instead of in law. 
Of course, there is a crucial point which needs some care: the weak convergence of the associated filtrations; the reader should recall that this notion requires the filtrations to be defined on the same probability space, which is the case throughout this work.

\subsection{Continuous martingales}

As a first example, we consider the case where the limit-BSDE is driven by a non-trivial, continuous, square-integrable martingale with independent increments, \emph{i.e.}, $X^{\infty,\natural}=0$.
In this case, we may assume that, for every $k\in \overline{\mathbb{N}}$, $X^{k,\natural}=0$ and, consequently, the domain of the generator $f^{k}$ is 
$\Omega \times \mathbb{R}_+\times \mathbb{R}^{\ell} \times \mathbb{R}^{\ell \times m}$, \emph{i.e.}, there is no dependence on elements of $\mathfrak{H}^{k}$.
Consider a sequence $(X^{k,\circ})_{k\in\overline{\mathbb{N}}}$ of square-integrable martingales with independent increments such that 
\begin{align*}
X^{k,\circ} \xrightarrow[k\to \infty]{\hspace{0.2em} (\textrm{J}_1,\mathbb{P})\hspace{0.2em}} X^{\infty,\circ}
\,\, \text{ and } \,\,
X^{k,\circ}_{\infty} \xrightarrow[k\to \infty]{\hspace{0.2em} \mathbb{L}^2(\Omega,\mathcal{G},\mathbb{P}; \mathbb{R}^{\ell})\hspace{0.2em}} X^{\infty,\circ}_{\infty}.
\end{align*}
Then, we obtain the weak convergence of the associated natural filtrations
\begin{align*}
\mathbb{F}^{X^{k,\circ}}\xrightarrow[k\to \infty]{\textrm{w}} \mathbb{F}^{X^{\infty,\circ}},
\end{align*}
by \citeauthor*{coquet2001weak} \cite[Proposition 2]{coquet2001weak}.
Once Conditions \ref{BgeneratorUI}--\ref{BFiniteStoppingTime} are satisfied, this yields a slight, but strict, generalisation of the stability results in \citeauthor*{briand2002robustness} \cite{briand2002robustness}, coming from the fact that we allow for unbounded Lipschitz constants and integrators of Lebesgue--Stieltjes integrals, which suitably combined form a bounded family of processes $(A^{k})_{k\in\overline{\mathbb{N}}}$, and the fact that we allow for unbounded terminal times.
The natural candidate for the It\^o integrator is Brownian motion, however more general \textit{continuous}, square-integrable martingales can also be considered.
This was already observed in \citeauthor*{briand2002robustness} \cite[Section 6]{briand2002robustness}.

\begin{remark}
Let us point out that the generality of our framework allows to conclude the convergence for It\^o integrals defined with respect to general $($\textit{i.e.} non-continuous$)$ martingales, once we assume that $X^{k,\circ}=0$, that $X^{k,\natural}$ is a square-integrable martingale $($with independent increments$)$, and we restrict the class of integrators to $U^k\equiv Z^k\cdot x$, for every $k\in\overline{\mathbb{N}}$.
In this case, {\rm Condition \ref{MXWPRP}} should be replaced by the \emph{strong} predictable representation property; see {\rm\citeauthor*{he1992semimartingale} \cite[Chapter XIII]{he1992semimartingale}}.
\end{remark}

\subsection{Probabilistic numerical schemes for deterministic equations}

As a second example, we are interested in approximating the solution of a (deterministic) integral equation with terminal condition via a sequence of solutions of BSDEs. 
To this end, for simplicity, we will modify the framework of the previous example by considering the stochastically trivial case for the limit-BSDE, \emph{i.e.},  $X^{\infty,\circ}=X^{\infty,\natural}=0$.
The filtration that we associate to the limit-BSDE is the constantly trivial filtration, \emph{i.e.}, $\mathcal{G}^{\infty}_t:=\{\emptyset,\Omega\}$ for every $t\in\mathbb{R}_+$. 
As a consequence, in the limit-BSDE all the random elements are trivialised.
The reader can verify that, for any sequence $(X^{k,\circ})_{k\in\mathbb{N}}$ such that 
$X^{k,\circ}_\infty \xrightarrow[k\to \infty]{\hspace{0.2em} \mathbb{L}^2(\Omega,\mathcal{G},\mathbb{P}; \mathbb{R}^{\ell}) \hspace{0.2em}} 0$ and for any sequence of filtrations associated to $(X^{k,\circ})_{k\in\mathbb{N}}$, {\mage Conditions \ref{MFilqlc}--\ref{Mfinalrv} are trivially satisfied.}
In particular, Condition \ref{MFilweak} is trivially satisfied because the only `random variables' one can choose from $\mathcal{G}^{\infty}_{\infty}$ are constant numbers.   
Let us point out, that we can freely choose the integrator $C^{\infty}$ of the Lebesgue--Stieltjes integral to be any c\`adl\`ag, increasing function. 
Almost all the notions of convergence are clear, \emph{e.g.}, the sequence of processes converges to a constant function and the sequence of random variables converges to a number.
However, the convergence of the generators $f^k$ to the generator $f^{\infty}$ of the (deterministic) integral limit-equation needs to be modified. 
If we want to have a unique solution for the limit-equation, we should choose the domain of $f^{\infty}$ to be $\Omega \times \mathbb{R}_+\times \mathbb{R}^{\ell}$, where of course there is no (true) randomness, \textit{i.e.}, no dependence on $\omega\in \Omega$.
Now, one should modify Condition \ref{Bgenerator}.\ref{Bgenerator:respects-J1-conv} as follows:
if $(W^{k})_{k\in\overline{\mathbb{N}}}$ is a sequence of c\`adl\`ag maps such that $W^k\xrightarrow{\hspace{0.5em}\textup{J}_1\hspace{0.5em}}W^\infty$, then for every $Z\in D^{\circ,\ell\times m}$
\begin{align*}
 \big(f^k(t, W^k_t,Z_t\big)_{t\in\mathbb{R}_+} 
 \xrightarrow[k\to\infty]{\hspace{1em}\textup{J}_1(\mathbb{R}^\ell)\hspace{1em}}
 \big(f^\infty(t, W^\infty_t)\big)_{t\in\mathbb{R}_+},\; \mathbb{P}\text{\rm--a.s.}
\end{align*}
In other words, one would expect the $Z$-variable to vanish in the limit.
Now, \Cref{BSDERobMainTheorem} yields that the sequence of solutions $(Y^k, Z^k, N^k)$ converges to $(Y,0,0)$, where $Y$ is the solution of the deterministic integral equation.
This allows us to design probabilistic schemes for the numerical solution of deterministic equations, which could be interesting in high-dimensional situations.
This example may also be compared to the results obtained in \citeauthor*{backhoff2020nonexponential} \cite{backhoff2020nonexponential}, where the authors consider the case where the $Z$-variable does not vanish and, thus, uncover interesting phenomena under suitable scaling.

\subsection{PII martingales}

Let us turn our interest to more general cases and consider a setting which is not covered, to the best of our knowledge, in the existing literature. 
More precisely, we assume additionally to \ref{MFilqlc}--\ref{BFiniteStoppingTime} that the sequence of integrators consists of martingales with independent increments and the associated filtrations are their natural ones.
Since we do not necessarily require the increments to be stationary, we consider a class which is broader than the class of L\'evy martingales; these are called PII martingales in \citeauthor{jacod2003limit} \cite{jacod2003limit}.
We underline that, in view of the independence of the increments, Condition \ref{MFilweak} is satisfied due to \citeauthor*{coquet2001weak} \cite[Proposition 2]{coquet2001weak}. 
Thus, if the sequence $(\overline{X}^{k})_{k\in\mathbb{N}}$ consists of processes indexed on a continuous time-set, then Theorem \ref{BSDERobMainTheorem} allows to conclude the convergence of perturbations of the limit-BSDE when we `wiggle' all the elements of the standard data.
On the other hand, if the sequence $(\overline{X}^{k})_{k\in\mathbb{N}}$ consists of processes indexed on a discrete time-set, then Theorem \ref{BSDERobMainTheorem} can be applied to derive the convergence of numerical schemes for the limit-BSDE.
The reader may immediately verify that the work of \citeauthor*{madan2015convergence} \cite{madan2015convergence} is a special case of this example.    

{ 
\subsection{BSDEs as dual problems}
   
As a final example, let us recall that in many applications BSDEs can provide an alternative characterisation for the solution of the problem at hand; a prominent example is stochastic optimal control problems. 
The current work covers Lipschitz BSDEs, so one can immediately derive the stability of, say, stochastic optimal control problems via \Cref{BSDERobMainTheorem}, as long as their dual problem corresponds to a BSDE with a Lipschitz generator. 
Such an approach has been used for the control of population dynamics by \citeauthor*{jusselin2019scaling} \cite{jusselin2019scaling}.
In \cite[Theorem 2]{jusselin2019scaling} the authors use stronger conditions, compared to \Cref{BSDERobMainTheorem}, with the only exception of \ref{Bgenerator}.
More precisely, the notion of convergence for the generators in \cite[Theorem 2]{jusselin2019scaling} requires suitable scaling of the $Z$-variable in order to obtain the desired limit.
Thus, \cite[Theorem 2]{jusselin2019scaling} cannot be recovered by \Cref{BSDERobMainTheorem}. 
The reader may recall from a previous example that the approach of using a suitable scaling has also been used by \citeauthor*{backhoff2020nonexponential} \cite{backhoff2020nonexponential}.
}   

%% file: Appendix_Notation.tex

\subsection{Definitions and complete notation}
\label{appendix-A1}

\begin{assumption}{2.10 of \cite{papapantoleon2016existence}}\label{assumptionC}
Let $\oX:=(X^\circ,X^\natural)\in\mathcal H^2(\mathbb R^m)\times\mathcal H^2(\mathbb R^n)$ and $C$ be a predictable, c\`adl\`ag and increasing process. 
The pair $(\oX,C)$ satisfies {\rm \cite[Assumption 2.10]{papapantoleon2016existence}} if 
each component of $\langle X^\circ \rangle$ is absolutely continuous with respect to $C$ and if the disintegration property given $C$ holds for the compensator $\nu^{X^\natural}$, \emph{i.e.},
there exists a transition kernel $K^{\oX}:(\Omega\times \mathbb R_+,\mathcal P)\longrightarrow \mathcal{R}(\mathbb R^n,\mathcal B(\mathbb R^n))$, where $\mathcal{R}\big(\mathbb R^n,\mathcal B(\mathbb R^n)\big)$ is the space of Radon measures on $\big(\mathbb R^n,\mathcal B(\mathbb R^n)\big)$, such that
\begin{align*}
  \nu^{X^{\natural}}(\omega;\d  t,\d  x)  =   K_t^{\oX}(\omega;\d  x)\, \d  C^{\oX}_t.
\end{align*}
\end{assumption}%
Following the notation of \cite{papapantoleon2016existence}, we define
\begin{gather*}
c^k := \left( \frac{\d  \langle X^{k,\circ}\rangle}{\d  C^k}\right)^{\frac{1}{2}},
\end{gather*}
in conjunction with the notation introduced in Section \ref{sec:Framework}, for $t\in\mathbb{R}_+$, $V:\Omega\times\mathbb{R}_+\times\mathbb{R}^n\rightarrow \mathbb{R}^p$ be such that $\mathbb{R}^n \ni x \mapsto V(\omega,t,x)\in\mathbb{R}^p$ is Borel-measurable
\begin{align}
\widehat K_{t}^k(V_s)(\omega)&:= \int_{\mathbb{R}^n}V(\omega,s,x) K_t^k(\omega;\d  x),
\label{def:hatK}
\end{align}
if $\displaystyle \int_{\mathbb{R}^n} (\Vert V(\omega,s,x)\Vert_2^2 \wedge 1) K_t^k(\omega;\d x)<\infty$, otherwise $\widehat K_{t}^k(V_s)(\omega):=\infty$;
\begin{gather*}
\zeta_t^{k,\natural}
  := \int_{\mathbb{R}^n}\nu^{k,\natural}(\{t\}\times\d x)
  =\widehat{\nu}_t^{k,\natural}(\mathds{1}_{\mathbb{R}^n})\le 1
\end{gather*}
and 
\begin{gather*}
\widehat{\nu}_t^{k,\natural} (V_s)(\omega)
  := \int_{\mathbb{R}^n} V(\omega,s,x) \nu^{k,\natural}(\omega;\{t\}\times\d x) ,
\end{gather*}
if $\displaystyle \int_{\mathbb{R}^n} \Vert V(\omega,s,x)\Vert_1 \nu^{k,\natural}(\omega;\{t\}\times\d x)<\infty,$
otherwise $\widehat{\nu}_t^{k,\natural} (V_s)(\omega):=\infty$;
and
\begin{align*}
\left(\tnorm{V(t,\cdot)}_{k,t}(\omega)\right)^2
  := \widehat K_{t}^k(\Vert V_s - \widehat{\nu}_t^{k,\natural}(V_t)\Vert_2^2)(\omega) 
  +  (1-\zeta_t^{k,\natural}(\omega))\Delta C^k_t(\omega) \Vert \widehat K_{t}^k(V_t)(\omega)\Vert_2^2.
\end{align*}
Finally, we define the following spaces for $\d  \mathbb P\otimes \d  C^{k}$--a.e. 
$(\omega,t)\in \Omega \times \mathbb{R}_+$ 
\[
\mathfrak{H}^{k}_{\omega,t}
  :=\overline{\left\{\mathcal U:(\mathbb R^n,\mathcal B(\mathbb 
    R^n))\longrightarrow (\mathbb{R}^{\ell},\mathcal B(\mathbb{R}^{\ell})),\
    \tnorm{\ \mathcal U(\cdot)}_{k,t}(\omega)<\infty\right\}}.\]
as well as
\[
\mathfrak{H}^{k}:=\big\{U:\Omega\times\mathbb{R}_+\times\mathbb 
  R^n\longrightarrow\mathbb{R}^{\ell},\ U(\omega,t,\cdot)\in \mathfrak H_{\omega,t}^{k}, \text{ for $\d  \mathbb P\otimes \d  C^{k}{\text{\rm--a.e.}}$ 
  $(\omega,t)\in \Omega \times \mathbb{R}_+$}\big\}.
  \]

For every $k\in\overline{\mathbb{N}}$ the following are true:
 the generator $(f^k)_{k\in\overline{\mathbb{N}}}:\Omega\times\mathbb{R}_{+}\times\mathbb{R}^{\ell}\times\mathbb{R}^{\ell\times m}\times\mathfrak{H}^k\longrightarrow\mathbb{R}^{\ell}$ satisfies Condition \ref{Fthree}, \emph{i.e.}, it is such that for any 
$(y,z,u)\in\mathbb{R}^{\ell}\times \mathbb{R}^{\ell \times m}\times \mathfrak{H}^k$, 
the map 
\[
(t,\omega)\longmapsto f^{k}(t,\omega,y,z,u_t(\omega;\cdot))\ \text{is $\mathcal F_t^{k}\otimes\mathcal B([0,t])$-measurable.}
\]
Moreover, $f$ satisfies a stochastic Lipschitz condition, that is to say there exist 
\[
  r^{k}:\big(\Omega\times\mathbb{R}_+,\mathcal{P}^{k}\big)\longrightarrow(\mathbb{R}_+,\mathcal{B}(\mathbb{R}_+))\ \textrm{ and }\
  \vartheta^{k} =(\theta^{k,\circ},\theta^{k,\natural}):\big(\Omega\times\mathbb{R}_+,\mathcal{P}^{k}\big)\longrightarrow(\mathbb{R}_+^2,\mathcal{B}(\mathbb{R}_+^{2})),
\] such that, for $\d  C \otimes \d  \mathbb{P}$--a.e. $(t,\omega)\in\mathbb R_+\times\Omega$
\begin{align*}
\begin{split}
  & \quad \big|f^{k}(t,\omega,y,z,u_t(\omega;\cdot))-f^{k}(t,\omega,y',z',u_t'(\omega;\cdot))\big|^2\\
  &\hspace{2em}\leq   r_t^{k}(\omega) |y-y'|^2  + \theta^{k,\circ}_t(\omega)\Vert c_t^{k}(\omega) (z-z')\Vert^2 + \theta^{k,\natural}_t(\omega) \left(\tnorm{u_t(\omega;\cdot)-u'_t(\omega;\cdot)}_{k,t}(\omega)\right)^2.    
\end{split}
\end{align*}
Moreover, to every generator $f^k$ we associate the predictable processes $\alpha^k_\cdot:=\max\{\sqrt{r^{k}_\cdot},\theta^{k,\circ},\theta^{k,\natural}_\cdot\}$, $A^k_\cdot := \big((\alpha^{k})^{2}\cdot C^{k}\big)_\cdot$ and the bound $\Phi^k$ such that $\Delta A^{k} \le \Phi^{k}$. 
All of them are described in Condition \ref{Ffour}.

\medskip
The spaces that are necessary for the existence and uniqueness result are introduced in \citep[Section 2.3]{papapantoleon2016existence}. 
For $\beta\geq0$, $p\in\mathbb{N}$ and for each $k\in\overline{\mathbb{N}}$ we use the following spaces
\label{def-spaces}
\begin{align*}
\mathbb{L}^{2,k}_{\beta}& : =   \Big\{\xi: \text{$\mathbb{R}^\ell$-valued, $\mathcal{G}^k_{\tau^k}$-measurable},
\; \Vert\xi\Vert^2_{\mathbb{L}^{2,k}_{\beta}}  := \mathbb{E}\Big[ \textrm{e}^{\beta A^k_{\tau^k}} \Vert \xi\Vert_2^2\Big] <\infty\Big\}, \\
\mathcal{H}^{2,k}_{\beta}(\mathbb{R}^{p})
  &:=\bigg\{ M\in\mathcal{H}^{2}(\mathbb{G}^k,\tau^k;\mathbb{R}^p):
  \Vert M\Vert_{\mathcal{H}^{2,k}_{\beta}(\mathbb{R}^{{}^p})}^2 :=\mathbb{E} \bigg[\int_{(0,\tau^k]} \textrm{e}^{\beta A^k_t} \d  \textrm{Tr} \bigg[\langle M\rangle_t^{\mathbb{G}^k}\bigg]\bigg]<\infty \bigg\}, \\
\mathbb H^{2,k}_{\beta}(\mathbb{R}^{p}) &:= 
\Bigg\{\phi: \text{$\mathbb{R}^p$-valued 
  $\mathbb{G}^k$-optional semimartingale with c\`adl\`ag paths and}\\
  &\qquad\quad \Vert \phi\Vert_{\mathbb H^{2,k}_{\beta}(\mathbb{R}^{{}^p})}^2:= 
  \E\bigg[\int_{(0,\tau^k]}\textrm{e}^{\beta A^k_t} \Vert\phi_t\Vert_2^2 \d  C^k_{t}\bigg]< \infty\bigg\},\\
\mathcal S^{2,k}_{\beta}(\mathbb{R}^{p})&:= 
  \bigg\{\phi: \text{$\mathbb{R}^\ell$-valued 
  $\mathbb{G}^k$-optional semimartingale with c\`adl\`ag paths and}\\
  &\qquad\quad \Vert \phi\Vert_{\mathcal S^{2,k}_{\beta}(\mathbb{R}^{{}^p})}^2:= 
    \E\bigg[ \sup_{t\in\llbracket0,\tau^k\rrbracket}\textrm{e}^{\beta A^k_t} \Vert \phi_t\Vert_2^2\bigg]<\infty\bigg\},\\
\mathbb H^{2,k,\circ} &:=
\bigg\{ Z:(\Omega\times\mathbb{R}_+,\mathcal{P}^{k}) 
\longrightarrow (\mathbb{R}^{\ell\times m},\mathcal B(\mathbb{R}^{\ell\times m})): 
\E\bigg[\int_{(0,\tau^{k}]} \textrm{Tr}\bigg[Z_t \frac{\d  \langle X\rangle_s}{\d  C_s}Z_t^\top\bigg]\d  C_t\bigg]<\infty \bigg\},\\
\mathbb H^{2,k,\circ}_{\beta} 
  &:= \bigg\{Z\in \mathbb H^{2,k,\circ}:
    \Vert Z\Vert_{\mathbb H^{2,k,\circ}_{\beta}}^2 :=\E\bigg[
    \int_{(0,\tau^k]} \textrm{e}^{\beta A^k_t} \d  \textrm{Tr}\big[\langle Z\cdot X^{k,\circ} \rangle^{\mathbb{G}^k}_t\big] \bigg]<\infty\bigg\}, \\
\mathbb H^{2,k,\natural} &:=\bigg\{ U:\big(\widetilde{\Omega},\widetilde{\mathcal{P}}^{k}\big) 
  \longrightarrow \big(\mathbb{R}^{\ell},\mathcal{B}(\mathbb{R}^{\ell})\big): 
  \E\bigg[\int_{(0,\tau^{k}]} \d  {\rm Tr} 
  \Big[\langle U\star\widetilde{\mu}^{k,\natural}\rangle^{\mathbb{G}^{k}}_t\Big]\bigg]<\infty  \bigg\},\\
\mathbb H^{2,k,\natural}_{\beta} 
  &:= \bigg\{U\in \mathbb H^{2,k,\natural}:  
  \norm{U}_{\mathbb H^{2,k,\natural}_{\beta}}:=\mathbb{E} \bigg[\int_{(0,\tau^k]}  \textrm{e}^{\beta A^k_t}\d  \textrm{Tr}\big[\langle U\star\widetilde{\mu}^{k,\natural}\rangle^{\mathbb{G}^k}_t\big]\bigg]<\infty\bigg\},\\
\mathcal{H}^{2,k,\perp}_{\beta} 
  &:=\Big\{M\in\mathcal{H}^{2,k}:\langle X^{k,\circ},\; M\rangle=0,\; M_{\mu^{k,\natural}}[\Delta M| \widetilde{\mathcal{P}}^{k}]=0, \text{ and }\Vert M\Vert^2_{\mathcal{H}^{2,k}_{\beta}}<\infty\Big\}.
\end{align*}

\vspace{0.5em}

Finally, for $(Y,Z,U,N)\in  \mathbb H^{2,k}_{\beta}\times
        \mathbb H^{2,k,\circ}_{\beta}\times
        \mathbb H^{2,k,\natural}_{\beta}\times
        \mathcal H^{2,k,\perp}_{\beta}$
we define
\begin{align*}
\norm{(Y,Z,U,N)}_{k,\beta}^2 &:= 
          \norm {\alpha Y}_{\mathbb H^{2,k}_{\beta}}^2  
          + \norm{Z}_{\mathbb H^{2,k,\circ}_{\beta}}^2 
          + \norm{U}_{\mathbb H^{2,k,\natural}_{\beta}}^2 
          + \norm N _{\mathcal H^{2,k,\perp}_{\beta}}^2,
\end{align*}
      and for
        $(Y,Z,U,N)\in   \mathcal S^{2,k}_{\beta}\times
        \mathbb H^{2,k,\circ}_{\beta}\times
        \mathbb H^{2,k,\natural}_{\beta}\times
        \mathcal H^{2,k,\perp}_{\beta}$, 
we define
\begin{align*}
\norm{(Y,Z,U,N)}_{\star,k,\beta}^2 &:=
          \norm {Y}_{\mathcal S^{2,k}_{\beta}}^2 
          + \norm{Z}_{\mathbb H^{2,k,\circ}_{\beta}}^2 
          + \norm{U}_{\mathbb H^{2,k,\natural}_{\beta}}^2 
          + \norm N _{\mathcal H^{2,k,\perp}_{\beta}}^2.
\end{align*}%
%

%% file: Appendix_SBSDE_2.tex

\subsection{Moore--Osgood's theorem}\label{App:MooreOsgood}  

Moore--Osgood's theorem, whose well-known form is \cref{MooreOsgoodA}, provides sufficient conditions for the existence of the limit of a doubly-indexed sequence and can be seen as a special case of \citeauthor*{rudin1976principles} \cite[Theorem 7.11]{rudin1976principles}.
Here we provide a second form, since the second time\footnote{The first one is in the proof of \cref{BSDERobMainTheorem}, while the second is in the proof of \cref{prop:MOAppl}} we need to apply the aforementioned theorem we need to relax the existence of the pointwise limits; compare the second conditions of \cref{MooreOsgoodA} and \cref{MooreOsgoodB}.  
The interested reader should consult \citeauthor*{hobson1907theory} \cite[Chapter VI, Sections 336--338]{hobson1907theory} for more details on the existence of the iterated limits and of the joint limit of a doubly-indexed sequence. 
Specifically for the validity of \cref{MooreOsgoodB} see \cite[Chapter VI, Section 337, p.~466]{hobson1907theory} and the reference therein.
\begin{theorem}\label{MooreOsgoodA}
Let $(\Gamma,d_\Gamma)$ be a metric space and $(\gamma_{k,p})_{k,p\in\mathbb N}$ be a sequence such that $\gamma_{\infty,p}:=\lim_{k\to\infty}\gamma_{k,p}$ exists for every $p\in\mathbb{N}$ and $\gamma_{k,\infty}:=\lim_{p\to\infty}\gamma_{k,p}$ exists for every $k\in\mathbb{N}$. 
If
\begin{enumerate}
\item\label{MOAi} $\lim\limits_{p\to\infty}\sup\limits_{k\in\mathbb{N}}d_\Gamma(\gamma_{k,p},\gamma_{k,\infty})=0$,
\item\label{MOAii} $\lim\limits_{k\to\infty}d_\Gamma(\gamma_{k,p},\gamma_{\infty,p})=0,$ for all $p\in {\mathbb N}$,
\end{enumerate}
then the joint limit $\lim\limits_{k,p\to\infty}\gamma_{n,k}$  exists. 
In particular holds $\lim\limits_{k,p\to\infty}\gamma_{k,p}=\lim\limits_{p\to\infty}\gamma_{\infty,p}=\lim\limits_{k\to\infty}\gamma_{k,\infty}.$
\end{theorem}

\begin{theorem}\label{MooreOsgoodB}
Let $(\mathbb{R},|\cdot|)$ and $(\gamma_{k,p})_{k,p\in\mathbb N}$ be a sequence such that $\gamma_{k,\infty}:=\lim_{p\to\infty}\gamma_{k,p}$ exists for every $k\in\mathbb{N}$.
If
\begin{enumerate}
\item\label{MOBi} $\lim\limits_{p\to\infty}\sup\limits_{k\in\mathbb{N}} |\gamma_{k,p}-\gamma_{k,\infty}|=0,$
\item\label{MOBii} $\lim\limits_{p\to\infty}\big\{\limsup\limits_{k\to\infty}\gamma_{k,p} - \liminf\limits_{k\to\infty}\gamma_{k,p}\big\}=0$,
\end{enumerate}
then the joint limit $\lim\limits_{k,p\to\infty}\gamma_{k,p}$  exists. 
In particular, it holds $\lim\limits_{k,p\to\infty}\gamma_{k,p}=\lim\limits_{p\to\infty}\lim\limits_{k\to\infty}\gamma_{k,p}=\lim\limits_{k\to\infty}\limsup\limits_{p\to\infty}\gamma_{k,p}$.
\end{theorem}
\subsection{Weak convergence of measures on the positive real line}

The aim of the current appendix is to provide a characterization of weak convergence of finite measures\footnote{We only consider positive measures, \emph{i.e.}, not signed ones.} defined on the positive real line in case the limit measure is atomless.
This characterization uses relatively compact sets of the Skorokhod space instead of relatively compact sets of the space of continuous
functions defined on $\mathbb{R}_+$ endowed with the supremum norm $\Vert \cdot\Vert_{\infty}$.
This result is of independent interest.

\begin{definition}\label{def:weak-Conv-finite-measure}
Let $(\mu^{k})_{k\in\overline{\mathbb N}}$ be a countable family of measures on $\bigl(\mathbb{R}_+, \mathcal B(\mathbb{R}_+)\bigr)$.
We will say that the sequence $(\mu^{k})_{k\in\mathbb N}$ converges weakly to the measure $\mu^{\infty}$ if for every $f:\bigl(\mathbb{R}_+,\|\cdot\|_\infty\bigr)\rightarrow \bigl(\mathbb{R}, |\cdot|\bigr)$ continuous and bounded holds
\begin{equation*}
\lim_{k\to\infty}\biggl|\int_{[0,\infty)}f(x)\,\mu^{k}(\d x) - \int_{[0,\infty)}f(x)\,\mu^{\infty}(\d x)\biggr|=0.
\end{equation*}
We denote the weak convergence of $(\mu^{k})_{k\in\mathbb N}$ to $\mu^{\infty}$ by $\mu^{k}\xrightarrow{\hspace{0.2em}\textup{w}\hspace{0.2em}} \mu^{\infty}.$
\end{definition}
In this section we will use the set
\begin{align*}
\mathcal W_{0,\infty}^+:=\Big\{F\in \mathbb{D}(\mathbb{R}_+): F \text{ is increasing with } F_0=0 \text{ and } \lim_{t\to\infty}F_t\in\mathbb{R}_+\Big\}.
\end{align*}

\begin{remark}\label{rem:W0infty-extension}
We can extend every element of $\mathcal W_{0,\infty}^+$, name $F$ an arbitrary element, on $[0,\infty]$ such that it is left-continuous at the symbol $\infty$ by defining $F_\infty:=\lim_{t\to\infty}F_t.$
\end{remark}

We provide in the following proposition some convenient equivalence results for the weak convergence of finite measures.
The statement is tailor-made to our needs, but the interested reader may consult \citeauthor*{bogachev2007measure} \cite[Sections 8.1--8.3]{bogachev2007measure}.
Then we provide in \cref{CharactWeakConv} a new, to the best of our knowledge, characterisation of weak convergence of finite measures to an atomless measure defined on the positive real line.

\begin{proposition}\label{EquivWeakConv}
Let $(\mu^{k})_{k\in\overline{\mathbb N}}$ be sequence of finite measures on $\bigl(\mathbb{R}_+, \mathcal B(\mathbb{R}_+)\bigr)$ with associated distribution functions $(F^{k})_{k\in\overline{\mathbb{N}}}$, where $F^{k}\in\mathcal{W}_{0,\infty}^{+}$, for $k\in\overline{\mathbb N}$.\footnote{For simplicity we assume that $\mu^{k}(\{0\})=0$.}
We assume the following
\begin{enumerate}
	\item the sequence $(\mu^{k})_{k\in\mathbb N}$ is bounded, \emph{i.e.} $\sup_{k\in\mathbb N}\mu^{k}(\mathbb{R}_+)=\sup_{k\in\mathbb N}F^{k}_\infty<\infty;$ 
	
	\smallskip
	\item the sequence $(\mu^{k})_{k\in\mathbb N}$ is tight. 
		In other words, for every $\varepsilon>0$ there exists $N_{\varepsilon}\in\mathbb{N}$ such that 
		\begin{align*}
		\sup_{k\in\mathbb N}\mu^{k}(\mathbb{R}_+\setminus [0,N_{\varepsilon}])=\sup_{k\in\mathbb N}|F^{k}_\infty - F^{k}_{N_{\varepsilon}}|<\varepsilon;
		\end{align*}
	\item the limit-measure $\mu^{\infty}$ is \emph{atomless}, \emph{i.e.} $\mu^{\infty}(\{t\})=0$ for every $t\in\mathbb{R}_+$.  
		Equivalently, $F^{\infty}$ is continuous.
\end{enumerate}
Then, the following are equivalent
\begin{enumerate}[label={\rm(\alph*)}]
\item\label{W1} $\displaystyle \mu^{k}\xrightarrow{\hspace{0.2em}\textup{w}\hspace{0.2em}} \mu;$

\smallskip
\item\label{W2} $\displaystyle F^{k}_t\xrightarrow[k\to\infty]{\hspace{0.2cm}\hspace{0.2cm}} F^{\infty}_t$, for every $t\in\mathbb{R}_+;$

\smallskip
\item\label{W3} $\displaystyle F^{k}\xrightarrow{\hspace{0.2em}\textup{J}_1(\mathbb{R})\hspace{0.2em}} F^{\infty};$

\smallskip
\item\label{W4} $\displaystyle F^{k}\xrightarrow{\hspace{0.2em}\textup{lu}\hspace{0.2em}} F^{\infty}$, where \textup{lu} stands for the locally uniform convergence$;$

\smallskip
\item\label{W5} $\displaystyle \underset{I\in\mathcal{I}_N}{\sup}|\mu^{k}(I)- \mu^{\infty}(I)|\xrightarrow[k\to\infty]{\hspace{0.2cm}\hspace{0.2cm}} 0$, for every $N\in\mathbb{N}$, where $\mathcal{I}_N:=\bigl\{I\subset[0,N]: I \text{ is interval}\bigr\}.$
\end{enumerate}
\end{proposition} 

\begin{proof}
The equivalence between \ref{W1} and \ref{W2} is a classical result, \emph{e.g.} see \cite[Proposition 8.1.8]{bogachev2007measure}.
The equivalences between \ref{W2}, \ref{W3} and \ref{W4} are provided by \citeauthor{jacod2003limit} \cite[Theorem VI.2.15.c.(i)]{jacod2003limit}. 
We obtain the equivalence between \ref{W4} and \ref{W5} in view of the validity of the following inequalities for every $N\in\mathbb{N}$:
\begin{align*}
&\sup_{t\in[0,N]}|F^{k}_t - F^{\infty}_t|
	\leq \underset{I\in\mathcal{I}_N}{\sup}|\mu^{k}(I) - \mu^{\infty}(I)|\\
&\hspace{1em}\leq
		\sup_{s,t\in[0,N]}|F^{k}_{t-} - F^{k}_s - F^{\infty}_{t-} + F^{\infty}_s| 
	+ 	\sup_{s,t\in[0,N]}|F^{k}_t - F^{k}_s - F^{\infty}_t + F^{\infty}_s|\\
&\hspace{2em}
	+	\sup_{s,t\in[0,N]}|F^{k}_{t-} - F^{k}_{s-} - F^{\infty}_{t-} + F^{\infty}_{s-}| \\
&\hspace{1em}\leq 
		3\sup_{t\in[0,N]}|F^{k}_t - F^{\infty}_t| 
	+	3\sup_{t\in[0,N]}|F^{k}_{t-} - F^{\infty}_{t-}| \\
&\hspace{1em}\leq 3\sup_{t\in[0,N]}|F^{k}_t - F^{\infty}_t|
				+ 3\sup_{t\in[0,N]}|F^{k}_{t-} - F^{k}_{t}| 
				+ 3\sup_{t\in[0,N]}|F^{k}_{t} - F^{\infty}_t| \\
&\hspace{1em}\leq 6\sup_{t\in[0,N]}|F^{k}_t - F^{\infty}_t| 
				+ 3\sup_{t\in[0,N]}|F^{k}_{t-} - F^{k}_{t}|.
\end{align*}
Then 
\begin{align*}
\hspace{-1em}\limsup_{k\to\infty}\underset{I\in\mathcal{I}_N}{\sup}|\mu^{k}(I) - \mu^{\infty}(I)|
&\leq 6	\lim_{k\to\infty}\sup_{t\in[0,N]}|F^{k}_t - F^{\infty}_t| 
	+ 3 \limsup_{k\to\infty}\sup_{t\in[0,N]}|F^{k}_{t-} - F^{k}_{t}|\\
&\leq 3 \sup_{t\in[0,N]}|F^{\infty}_{t-}-F^{\infty}_t| = 0,
\end{align*}
where we used \cite[Lemma VI.2.5]{jacod2003limit} in the last inequality. 
Hence, for every $N\in\mathbb{N}$,
\begin{equation*}
\lim_{k\to\infty} \sup_{t\in[0,N]}|F^{k}_t - F^{\infty}_t| = 0, \text{ if and only if } 
\lim_{k\to\infty}\underset{I\in\mathcal{I}_N}{\sup}|\mu^{k}(I) - \mu^{\infty}(I)|=0. \qedhere
\end{equation*}
\end{proof}
%
	Assume the framework of the previous proposition and that its Part \ref{W1} is true. Then, for $\mathds{1}_{\mathbb{R}_+}$ we get 
	\[\mu^{k}(\mathbb{R}_+)\xrightarrow[k\to\infty]{\hspace{0.2cm}|\cdot|\hspace{0.2cm}}\mu^{\infty}(\mathbb{R}_+),\]
	or, alternatively written, $F^{k}_\infty\xrightarrow[k\to\infty]{\hspace{0.2cm}|\cdot|\hspace{0.2cm}} F^{\infty}_{\infty}$. 
	The following lemma allows us to use a stronger convergence for the sequence of distribution functions in the special case we consider, namely the convergence holds under the supremum norm. 
	In other words, we have the convergence of the distribution functions under the Kolmogorov--Smirnov distance.
	For the sake of completeness we present the details.
%
%
\begin{lemma}\label{JtoUConv}
Let $(F^{k})_{k\in\overline{\mathbb N}} \subset  \mathcal W_{0,\infty}^+$ which satisfies the following properties
\begin{enumerate}
\item \label{LemUCJ1} $F^{k}{\xrightarrow[k\to\infty]{\hspace{0.2em}\textup{J}_1(\mathbb{R})\hspace{0.2em}}} F^{\infty};$
\item \label{LemUCinf} $F^{k}_\infty\xrightarrow[k\to\infty]{\hspace{0.2cm}|\cdot|\hspace{0.2cm}} F^{\infty}_\infty;$
\item \label{LemUCCont} the limit function $F^{\infty}$ is continuous.
\end{enumerate} 
Then it holds $F^{k}\xrightarrow[k\to\infty]{\hspace{0.2cm}\|\cdot\|_\infty\hspace{0.1cm}} F^{\infty}$.
\end{lemma}
\begin{proof}
Let us fix an arbitrary $\varepsilon>0.$
We start by exploiting the fact that $F^{\infty}\in\mathcal W_{0,\infty}^+$, \emph{i.e.}, there exists $K>0$ such that 
\begin{equation}\label{JtoUconv-1}
\sup_{t\geq K}\big\{F^{\infty}_\infty-F^{\infty}_t\big\} < \frac{\varepsilon}{8}.
\end{equation}
By \ref{LemUCinf}, there exists $k_1=k_1(\varepsilon)\in \mathbb{N}$ such that 
\begin{equation}\label{JtoUconv-2}
\sup_{k\ge k_1}|F^{k}_\infty-F^{\infty}_\infty| < \frac{\varepsilon}{8}.
\end{equation}
By \ref{LemUCCont} there exists $\bar{s}>K$ such that $F^{\infty}_{\bar{s}}=\frac{1}{2}(F^{\infty}_K+F^{\infty}_\infty)$, \emph{i.e.}, $|F^{\infty}_{\infty} - F^{\infty}_{\bar{s}}|<\frac{\varepsilon}{16}$. 
By \ref{LemUCJ1}, \ref{LemUCCont} and \cref{EquivWeakConv}.\ref{W4}, there exists $k_2=k_2(\varepsilon,\bar{s})$ such that 
\begin{equation}\label{JtoUconv-3}
\sup_{k\ge k_2} \sup_{0\leq t\leq\bar{s}}|F^{k}_t-F^{\infty}_t|
<\frac{1}{2}(F^{\infty}_\infty-F^{\infty}_K)
<\frac{\varepsilon}{16},
\end{equation}
where the last inequality is valid in view of \eqref{JtoUconv-1}, since $\bar s>K.$
Using the fact that the functions are increasing we derive 
\begin{align*}
&\sup_{k\ge k_1\vee k_2}\bigg\{\sup_{\bar{s}\leq t<\infty} F^{k}_\infty-F^{k}_t\bigg\}
= \sup_{k\ge k_1\vee k_2}\big\{F^{k}_\infty-F^{k}_{\bar{s}} \big\}\\
&\hspace{1em}\overset{\phantom{\eqref{JtoUconv-1}}}{\le}
\sup_{k\ge k_1\vee k_2}\bigl|F^{k}_\infty-F^{\infty}_\infty\bigr| 
+\sup_{k\ge k_1\vee k_2}\bigl|F^{\infty}_\infty-F^{\infty}_{\bar{s}} \bigr| 
+\sup_{k\ge k_1\vee k_2}\bigl|F^{\infty}_{\bar{s}}-F^{k}_{\bar{s}} \bigr|
\overset{\eqref{JtoUconv-2}}{\underset{\eqref{JtoUconv-3}}{<}}\frac{\varepsilon}{8}+\frac{\varepsilon}{16}+\frac{\varepsilon}{16}
=\frac{\varepsilon}{4}.
\numberthis\label{JtoUconv-4}
\end{align*} 
Therefore, by combining the above, we obtain for every $k\ge k_1\vee k_2$
\begin{equation*}
\sup_{t\in\mathbb{R}_+}|F^{k}_t-F^{\infty}_t| 
\leq \sup_{0\leq t\leq\bar{s}}|F^{k}_t-F^{\infty}_t| + \sup_{\bar{s}\leq t \leq\infty}|F^{k}_t-F^{\infty}_t|
\overset{\eqref{JtoUconv-3},\eqref{JtoUconv-4}}{\underset{\eqref{JtoUconv-2}, \eqref{JtoUconv-1}}{<}} 
\frac{\varepsilon}{16} +\frac{\varepsilon}{4} + \frac{\varepsilon}{4} +\frac{\varepsilon}{8}<\varepsilon. \qedhere
\end{equation*}
\end{proof}

The next step is to provide a characterisation of the aforementioned weak convergence of measures using relatively compact sets of the Skorokhod space $\mathbb{D}(\mathbb{R}_+;\mathbb{R}^{p})$.
Two remarks are in order. 
The first remark is that a relatively compact set of the Skorokhod space is, in general, only locally uniformly bounded, see \citeauthor{jacod2003limit} \cite[Theorem VI.1.14.b]{jacod2003limit}.
Therefore, we will restrict ourselves to those relatively compact subsets which are uniformly bounded.
The second remark is that proving \cref{CharactWeakConv} for integrands in $\mathbb{D}(\mathbb{R}_+;\mathbb{R})$, allows us to generalise the statement for integrands in $\mathbb{D}(\mathbb{R}_+;\mathbb{R}^{p})$, where the norm in convergence \eqref{UnifIntConv} should be substituted by $\Vert \cdot \Vert_1$.
Indeed, one can verify that the modulus $w'$%
\footnote{For $\alpha\in\mathbb{D}(\mathbb{R}_+;\mathbb{R}^p)$, $\theta>0$, $N\in\mathbb{N}$ we use the modulus $w'_N(\alpha,\theta)$ as defined in \cite[Section VI.1]{jacod2003limit}.} 
 takes care of the `uniform local behaviour' of elements of a relatively compact set, independently of the dimension of the state space.

\medskip    

The following theorem can be regarded as an extension of \citeauthor*{parthasarathy1972probability} \cite[Theorem II.6.8]{parthasarathy1972probability}, in the special case that the limiting measure is an atomless measure on $\mathbb{R}_+$.
The proof of \cref{CharactWeakConv} is inspired by the proof of the aforementioned theorem.
However, we need to underline that in \cite[Theorem II.6.8]{parthasarathy1972probability} the equi-continuity of the integrands is used in order to determine for each arbitrary, but fixed, measure $\mu$ a countable cover of $\mathbb{R}_+$ consisting of $\mu$-continuity sets on which the integrands are uniformly arbitrarily small.
This is a property which is not possible to be obtained for relatively compact subsets of $\mathbb{D}(\mathbb{R}_+;\mathbb{R})$ and an arbitrary measure $\mu$ on $(\mathbb{R}_+,\mathcal{B}(\mathbb{R}_+))$.
Hence, we will need to properly use the `equi--right-continuity' of the elements of a relatively compact subset of $\mathbb{D}(\mathbb{R}_+;\mathbb{R})$.
\begin{theorem}\label{CharactWeakConv}
Consider the measurable space $(\mathbb{R}_+, \mathcal B(\mathbb{R}_+))$, let $(\mu^{k})_{k\in\overline{\mathbb N}}$ be a sequence of finite measures such that $\mu^{\infty}$ is an atomless finite measure and $\mathcal{A}\subset\mathbb{D}(\mathbb{R}_+;\mathbb{R})$.
Assume that the following conditions are true
\begin{enumerate}
\item\label{WeakConvergence}	$\mu^{k} \xrightarrow{\hspace{0.2em}\textup{w}\hspace{0.2em}} \mu^{\infty};$
\item\label{UnifBound} 			$\mathcal A$ is uniformly bounded, \emph{i.e.}, $\sup_{\alpha\in\mathcal{A}} \|\alpha\|_{\infty}<\infty;$  
\item\label{equircontin} 		$\mathcal A$ is uniformly `equi-right-continuous', \emph{i.e.}, $\lim\limits_{\zeta\downarrow 0}\sup\limits_{\alpha\in\mathcal{A}}  w'_N(\alpha,\zeta)=0$ for every $N\in\mathbb{N}.$\footnote{The combination of parts \ref{UnifBound} and \ref{equircontin} guarantee that $\mathcal{A}$ is relatively compact.}
\end{enumerate}
Then
\begin{align}\label{UnifIntConv}
\lim_{k\to\infty}\sup_{\alpha\in\mathcal{A}}\;  \biggl|\int_{[0,\infty)} \alpha(x)\,\mu^{k}(\d x) - \int_{[0,\infty)} \alpha(x)\,\mu^{\infty}(\d x)\biggr|=0.
\end{align}
\end{theorem}
\begin{proof}
We will decompose initially the quantity whose convergence we intent to prove as follows
\begin{align}\label{CharactWeakConv-1}
\notag&\sup_{\alpha\in\mathcal{A}}  \biggl|\int_{[0,\infty)} \alpha(x)\,\mu^{k}(\d x) - \int_{[0,\infty)} \alpha(x)\,\mu^{\infty}(\d x)\biggr|\\
&\le \sup_{\alpha\in\mathcal{A}}  \biggl|\int_{[0,N]} \alpha(x)\,\mu^{k}(\d x) - \int_{[0,N]} \alpha(x)\,\mu^{\infty}(\d x)\biggr|
+\sup_{\alpha\in\mathcal{A}}  \biggl|\int_{[N,\infty)} \alpha(x)\,\mu^{k}(\d x)\biggr| 
+\sup_{\alpha\in\mathcal{A}}   \biggl|\int_{[N,\infty)} \alpha(x)\,\mu^{\infty}(\d x)\biggr|,
\end{align}
for some $N>0$ to be determined.
Then, the first summand of the right-hand side will be decomposed as follows
\begin{flalign*}
&\sup_{\alpha\in\mathcal{A}}  \biggl|\int_{[0,N]} \alpha(x)\,\mu^{k}(\d x) - \int_{[0,N]} \alpha(x)\,\mu^{\infty}(\d x)\biggr|&\\
&\begin{multlined}[c][0.9\displaywidth]\numberthis\label{CharactWeakConv-2}
\hspace{1em}\le \sup_{\alpha\in\mathcal{A}}		\biggl|		\int_{[0,N]}\alpha(x)\, \mu^{k}(\d x) 	-	\int_{[0,N]}\alpha(x)\, \widetilde{\mu}^{k,\alpha}(\d x)\biggr|	\\
+\sup_{\alpha\in\mathcal{A}}		\biggl|		\int_{[0,N]}\alpha(x)\,  \widetilde{\mu}^{k,\alpha}(\d x) 	-	\int_{[0,N]}\alpha(x)\, \widetilde{\mu}^{\infty,\alpha}(\d x)\biggr|\\
+\sup_{\alpha\in\mathcal{A}}		\biggl|		\int_{[0,N]}\alpha(x)\,  \widetilde{\mu}^{\infty,\alpha}(\d x) 	-	\int_{[0,N]}\alpha(x)\, \mu^{\infty}(\d x)\biggr|,
\end{multlined}&
\end{flalign*}
where the measures $\widetilde{\mu}^{k,\alpha}$, for $k\in\overline{\mathbb N}$ and $\alpha\in\mathcal A$, will be constructed given the measure $\mu^{k}$ and the function $\alpha$.
For the second and third summand of \eqref{CharactWeakConv-1} we will prove that they become arbitrarily small for large $N>0$.
Then, we will conclude once we obtain the convergence to $0$ (as $k\to\infty$) of the summands of \eqref{CharactWeakConv-2}.

\medskip
To this end, let us fix an $\varepsilon>0$.

\medskip
\raisebox{0.5pt}{\tiny$\blacksquare$}
Condition \ref{WeakConvergence} implies that $\sup_{k\in\overline{\mathbb N}}\mu^{k}(\mathbb{R}_+)<\infty$, 
as well as the tightness of the sequence since $\mu^{k}(\mathbb{R}_+)\xrightarrow[k\to\infty]{}\mu^{\infty}(\mathbb{R}_+)$.
Hence, for $M:=\sup_{\alpha\in\mathcal{A}}\|\alpha\|_\infty<\infty$, which is Condition \ref{UnifBound}, there exists $N=N(\varepsilon,M)>0$ such that 
\begin{align*}
\sup_{k\in\overline{\mathbb N}} \mu^{k}\bigl([N,\infty)\bigr)<\frac{\varepsilon}{2M}, \text{ which implies }
\sup_{k\in\overline{\mathbb N},\alpha\in\mathcal A}\bigg|\int_{[N,\infty)}\alpha(x)\,\mu^{k}(\d x)\bigg|<\frac{\varepsilon}{2}.
\end{align*}
In other words, we have proven that the second and third summand of the right-hand side of \eqref{CharactWeakConv-1} become arbitrarily small for large $N.$
In the following $N$ is assumed fixed, but large enough so that the above hold.

\medskip
\raisebox{0.5pt}{\tiny$\blacksquare$}
We proceed now to construct, for every $k\in\overline{\mathbb N}$ and $\alpha\in\mathcal A$, a suitable measure $\widetilde{\mu}^{k,\alpha}$ on $\bigl(\mathbb{R}_+,\mathcal B(\mathbb{R}_+)\bigr)$ such that we can obtain the convergence of the summands in \eqref{CharactWeakConv-2}.
Define $K_N:={\sup_{k\in\overline{\mathbb N}}}\ \mu^{k}\left([0,N]\right)<\infty$.
By \ref{equircontin} there exists $\zeta=\zeta(\varepsilon, K_N)>0$ such that for every $\zeta'<\zeta$ it holds ${\sup_{\alpha\in\mathcal{A}}} w'_N(\alpha,\zeta')<\frac{\varepsilon}{4K_N}$.
By the definition of $w'$ (recall we use  \cite[Section VI.1]{jacod2003limit}) for every $\alpha\in\mathcal{A}$ and every $\zeta'<\zeta$, there exists a $\zeta'$-sparse set $P^{N,\zeta'}_\alpha$, \emph{i.e.} 
\begin{multline*}
P^{N,\zeta'}_\alpha:= 
	\Bigl\{0=t(0;P^{N,\zeta'}_\alpha)<t(1;P^{N,\zeta'}_\alpha)<\ldots<t(\kappa^\alpha;P^{N,\zeta'}_\alpha)=N:\\
		\min_{i\in\{1,\ldots,\kappa^\alpha-1\}} \bigl(t(i;P^{N,\zeta'}_\alpha)-t(i-1;P^{N,\zeta'}_\alpha)\bigr)>\zeta' \Bigr\},
\end{multline*}
such that
\begin{equation}\label{uniformrightcontinuity}
\begin{multlined}[c][0.85\displaywidth]
\max_{i\in\{1,\ldots,\kappa^\alpha\}} \sup\Bigl\{ |\alpha(s)-\alpha(u)|:s,u\in\bigl[t(i-1;P^{N,\zeta'}_\alpha),t(i;P^{N,\zeta'}_\alpha)\bigr)\Bigr\}
<w'_N(\alpha,\zeta')+\frac{\varepsilon}{4K_N}\\
<{\sup_{\alpha\in\mathcal{A}}}w'_N(\alpha,\zeta')+\frac{\varepsilon}{4K_N}
<\frac{\varepsilon}{2K_N}.
\end{multlined}
\end{equation}
To sum up, for every $(\alpha,\zeta')\in\mathcal{A}\times (0,\zeta)$ we can find a partition of $[0,N]$,
$P^{N,\zeta'}_\alpha$, satisfying \eqref{uniformrightcontinuity}. 
For the following, given a finite measure on $\bigl(\mathbb{R}_+, \mathcal B(\mathbb{R}_+)\bigr)$, 
called $\lambda$,  and a sparse set $P^{N,\zeta'}_\alpha$ we define
\begin{enumerate}[label=\text{\tiny$\blacktriangleright$}]
\item 	the interval $A(i;{P^{N,\zeta'}_\alpha}):=[t(i-1;P^{N,\zeta'}_\alpha),t(i;P^{N,\zeta'}_\alpha))$ for $i\in\{1,\ldots, \kappa^\alpha\}$ and
\item 	the operator $\lambda \overset{\widetilde{\phantom{\lambda}}^{\alpha}}{\longmapsto}\widetilde{\lambda}^\alpha$  where the finite measure $\widetilde{\lambda}^\alpha:\bigl(\mathbb{R}_+, \mathcal B(\mathbb{R}_+)\bigr)\rightarrow[0,\infty)$ is defined by 
	\begin{align*}
	\widetilde{\lambda}^\alpha(\cdot):=\sum_{i=1}^{\kappa^\alpha}\lambda \bigl(A(i;{P^{N,\zeta'}_\alpha})\bigr)\delta_{t(i-1;P^{N,\zeta'}_\alpha)} (\cdot).
	\end{align*}
	
\end{enumerate}
Here, $\delta_{x}$ is the Dirac measure sitting at the point $x$. 
Using the notation introduced above, we have
\begin{align*}
&\Biggl|\int_{[0,N]}\alpha(x)\,\lambda(\d x) - \int_{[0,N]}\alpha(x)\,\widetilde{\lambda}^\alpha(\d x)\Biggr|
= \Biggl|\sum_{i=1}^{\kappa^\alpha}\int_{A(i;{P^{N,\zeta'}_\alpha})}\alpha(x)\,\lambda(\d x) - \sum_{i=1}^{\kappa^\alpha}\int_{A(i;{P^{N,\zeta'}_\alpha})}\alpha(x)\,\widetilde{\lambda}^\alpha(\d x)\Biggr|\\
&\hspace{1em}\overset{\phantom{\eqref{uniformrightcontinuity}}}{\le}  
\sum_{i=1}^{\kappa^\alpha}\Biggl|\int_{A(i;{P^{N,\zeta'}_\alpha})}\alpha(x)\,\lambda(\d x) - \int_{A(i;{P^{N,\zeta'}_\alpha})}\alpha(x)\,\widetilde{\lambda}^\alpha(\d x)\Biggr|\\
&\hspace{1em}\overset{\phantom{\eqref{uniformrightcontinuity}}}{=} 
\sum_{i=1}^{\kappa^\alpha} \Biggl|\int_{A(i;{P^{N,\zeta'}_\alpha})}\alpha(x)\,\lambda(\d x) -	\alpha\bigl(t(i-1;P_\alpha^{N,\zeta'}\bigr) \lambda\left({[t(i-1;P^{N,\zeta'}_\alpha),t(i;P^{N,\zeta'}_\alpha))}\right)\Biggr|\\
&\hspace{1em}\overset{\phantom{\eqref{uniformrightcontinuity}}}{=} 
\sum_{i=1}^{\kappa^\alpha}\biggl|\int_{A(i;{P^{N,\zeta'}_\alpha})}\alpha(x) -	\alpha\bigl(t(i-1;P_\alpha^{N,\zeta'}\bigr) \lambda(\d x) \biggr|
\le \sum_{i=1}^{\kappa^\alpha}	\int_{A(i;{P^{N,\zeta'}_\alpha})}\Bigl|\alpha(x) -	\alpha\bigl(t(i-1;P_\alpha^{N,\zeta'}\bigr) \Bigr|\lambda(\d x)\\
&\hspace{1em}\overset{\eqref{uniformrightcontinuity}}{\leq}
\frac{\varepsilon}{2K_N} \sum_{i=1}^{\kappa^\alpha}	\int_{A(i;{P^{N,\zeta'}_\alpha})}\lambda(\d x) =
\frac{\varepsilon}{2K_N} \lambda\left([0,N]\right).
\numberthis\label{CharactWeakConv-3}
\end{align*}
For every $k\in\overline{\mathbb N}$ and $\alpha\in\mathcal A$ we define $\widetilde{\mu}^{k,\alpha} :=\widetilde{(\mu^{k})}^\alpha$.
Now, using the approximation \eqref{CharactWeakConv-3} for $\mu^{k}$, for $k\in\overline{\mathbb N}$, we obtain
\begin{align}
\sup_{k\in\overline{\mathbb N},\alpha\in\mathcal{A}}
\biggl|		\int_{[0,N]}\alpha(x)\, \mu^{k}(\d x) 	-	\int_{[0,N]}\alpha(x)\, \widetilde{\mu}^{k,\alpha}(\d x)\biggr|<\frac{\varepsilon}{2}.
\label{Cnintappr}
\end{align}
To sum up, we have constructed the measures $\widetilde{\mu}^{k,\alpha}$ such that the first and third summand of \eqref{CharactWeakConv-2} become arbitrarily small uniformly on $k\in\overline{\mathbb N}$ and on $\alpha\in\mathcal A.$

\vspace{0.5em}
\raisebox{0.5pt}{\tiny$\blacksquare$} We can conclude \eqref{UnifIntConv} once we obtain the validity of 
\begin{align*}
\lim_{k\to\infty}\sup_{\alpha\in\mathcal{A}} \biggl|	\int_{[0,N]}\alpha(x)\, \widetilde{\mu}^{k,\alpha}(\d x)  	-	\int_{[0,N]}\alpha(x)\, \widetilde{\mu}^{\infty,\alpha}(\d x)\biggr|=0.
\end{align*}
Indeed, for $\zeta'\in(0,\zeta)$,
\begin{align*}
&\biggl|	\int_{[0,N]}\alpha(x)\,\widetilde{\mu}^{k,\alpha}(\d x)  	-	\int_{[0,N]}\alpha(x)\,\widetilde{\mu}^{\infty,\alpha}(\d x)\biggr|\\
&\hspace{1em}=\bigg|\sum_{i=1}^{\kappa^\alpha}\int_{A(i;{P^{N,\zeta'}_\alpha})}\alpha(x)\,\widetilde{\mu}^{k,\alpha}(\d x) - \sum_{i=1}^{\kappa^\alpha}\int_{A(i;{P^{N,\zeta'}_\alpha})}\alpha(x)\,\widetilde{\mu}^{\infty,\alpha}(\d x)\bigg|\\
&\hspace{1em}= \biggl|\sum_{i=1}^{\kappa^\alpha}\Bigl( \int_{A(i;{P^{N,\zeta'}_\alpha})}\alpha(x)\,\widetilde{\mu}^{k,\alpha}(\d x) 	-	\int_{A(i;{P^{N,\zeta'}_\alpha})}\alpha(x)\,\widetilde{\mu}^{\infty,\alpha}(\d x)\Bigr)\bigg|\\
&\hspace{1em}= \bigg|\sum_{i=1}^{\kappa^\alpha}\alpha\bigl(t(i;P^{N,\zeta'}_\alpha)\bigr)	 \Bigl[\mu^{k} (A(i;{P^{N,\zeta'}_\alpha})) 		- \mu   \bigl(A(i;{P^{N,\zeta'}_\alpha})\bigr)\Bigr]\bigg|\\
&\hspace{1em}\overset{\ref{UnifBound}}{\leq}
 \|\alpha\|_{\infty}\sum_{i=1}^{\kappa^\alpha}
 \Bigl|	\mu^{k} \bigl(A(i;{P^{N,\zeta'}_\alpha})\bigr) -	\mu^{\infty}  \bigl(A(i;{P^{N,\zeta'}_\alpha})\bigr)\Bigr|\\
&\hspace{1em}\le
\sup_{\alpha\in\mathcal{A}} \|\alpha\|_\infty\cdot\kappa^\alpha\cdot 
 \sup_{i\in\{1,\ldots,\kappa^\alpha\}}\Bigl|	\mu^{k} \bigl(A(i;{P^{N,\zeta'}_\alpha})\bigr) - 	\mu^{\infty}   \bigl(A(i;{P^{N,\zeta'}_\alpha})\bigr)\Bigr|\\
&\hspace{1em}\le
\sup_{\alpha\in\mathcal{A}} \|\alpha\|_\infty\cdot\kappa^\alpha\cdot  \sup_{I\in\mathcal{I}_N}\bigl|\mu^{k} (I)- 	\mu^{\infty}  (I)\bigr|.
\numberthis\label{CharactWeakConv-4}
\end{align*}
Let us now define $m(P^{N,\zeta'}_\alpha):=\min\bigl\{t(i;{P^{N,\zeta'}_\alpha}) - t(i-1;{P^{N,\zeta'}_\alpha}),\;i\in\{1,\ldots,\kappa^\alpha\}\bigr\}>\zeta'$ by the definition of $P^{N,\zeta'}_\alpha$.
Therefore,
\begin{equation*}
\sup_{\alpha\in\mathcal{A}}\ \kappa^\alpha
	\leq \sup_{\alpha\in\mathcal{A}}\bigg\lceil \frac{N}{m(P^{N,\zeta'}_\alpha)} \bigg\rceil\footnotemark
		\leq \bigg\lceil \frac{N}{\zeta'} \bigg\rceil<\infty,
\end{equation*}
where $\lceil x\rceil$ denotes the least integer greater than or equal to $x.$
Hence, for every fixed sparse set $P^{N,\zeta'}_\alpha$ which satisfies \eqref{uniformrightcontinuity} we obtain
\footnotetext{We denote by $\lceil x\rceil$ the least integer greater than or equal to $x$.}
\begin{equation*}
\sup_{\alpha\in\mathcal{A}} 
\biggl|	\int_{[0,N]}\alpha(x)\,\widetilde{\mu}^{k,\alpha}(\d x) -	\int_{[0,N]}\alpha(x)\,\widetilde{\mu}^{\infty,\alpha}(\d x)\biggr|
\le
\sup_{\alpha\in\mathcal{A}} \|\alpha\|_\infty\times 
\sup_{\alpha\in\mathcal{A}} \kappa^\alpha\times 
\sup_{I\in\mathcal{I}_N}
\bigl|	\mu^{k} \left(I\right) - 	\mu   \left(I\right)\big|\xrightarrow[k\to\infty]{\ref{WeakConvergence}}0,
\end{equation*}
where we have used \cref{EquivWeakConv}.\ref{W5} which is equivalent to \ref{WeakConvergence}.
\end{proof}
%
%
%
%
We have presented the previous theorem for $\mathcal A\subset \mathbb{D}(\mathbb{R}_+;\mathbb{R})$.
However, as we mentioned before, the result can be readily adapted for $\mathcal A\subset \mathbb{D}(\mathbb{R}_+;\mathbb{R}^p)$.
This will be presented in the following proposition.
%
%
%
%
%
%
\begin{proposition}\label{prop:MOAppl}
Consider the measurable space $\bigl(\mathbb{R}_+, \mathcal B(\mathbb{R}_+)\bigr)$ and let $(\mu^{k})_{k\in\overline{\mathbb N}}$ be a sequence of finite measures such that $\mu^{k}(\{0\})=0$ for every $k\in\overline{\mathbb N}$ and $\mu^{\infty}$ is atomless.
Additionally, let  $\mathcal A:=(\alpha^k)_{k\in\overline{\mathbb N}}\subset\mathbb{D}(\mathbb{R}_+;\mathbb{R}^p)$.
Assume the following to be true
\begin{enumerate}
\item\label{WeakConvergence-prop}	$\mu^{k} \xrightarrow{\hspace{0.2em}\textup{w}\hspace{0.2em}} \mu^{\infty};$
\item\label{UnifBound-prop} 		$\mathcal A$ is uniformly bounded, \emph{i.e.} $\sup_{\alpha\in\mathcal{A}} \|\alpha\|_{\infty}<\infty;$
\item\label{equircontin-prop} 		$\mathcal A$ is $\delta_{\textup{J}_1(\mathbb{R}^p)}$-convergent.
									In other words, $\alpha^k\xrightarrow[k\to\infty]{\hspace{0.2em}\textup{J}_1(\mathbb{R}^p)\hspace{0.2em}}\alpha^\infty.$
\end{enumerate}
Then
\begin{equation*}
\lim_{k,m\to\infty} \biggl\|\int_{[0,\infty)} \alpha^k(x)\,\mu^{m}(\d x) - \int_{[0,\infty)} \alpha^\infty(x)\,\mu^{\infty}(\d x)\biggr\|_1=0.
\end{equation*}
\end{proposition}
\begin{proof}
Let us define, for $k\in\overline{\mathbb N}$  and $m\in\mathbb N$,
\begin{equation*}
	\gamma^{k,m}:=
	\biggl\|\int_{[0,\infty)}\alpha^k(x)\, \mu^{m}(\d x) - \int_{[0,\infty)}\alpha^k(x)\, \mu^{\infty}(\d x)\biggr\|_1 .
\end{equation*}
We are going to apply \cref{MooreOsgoodB} in order to obtain the required result.
By \cref{CharactWeakConv} we have that 
\begin{equation*}
	\lim_{m\to\infty}\sup_{k\in\overline{\mathbb N}}\big|\gamma^{k,m}\big|=0.
\end{equation*}
In other words, Condition \ref{MOBi} of \cref{MooreOsgoodB} is satisfied.
Let us prove, now, that Condition \ref{MOBii} of the aforementioned theorem is also satisfied, \emph{i.e.} we need to prove that 
\begin{equation*}
	\lim_{m\to\infty}\bigg(\limsup_{k\to\infty}\gamma^{k,m} - \liminf_{k\to\infty}\gamma^{k,m} \bigg)=0.
\end{equation*}
However,  it is sufficient to prove that $\lim_{m\to\infty}\limsup_{k\to\infty}\gamma^{k,m}=0,$ since the elements of the doubly-indexed sequence are positive.
To this end, we have
\begin{flalign*}
&\hspace{1em}\limsup_{k\to\infty}\gamma^{k,m} &\\
&\hspace{2em}\begin{multlined}[c][0.9\displaywidth]
\le\limsup_{k\to\infty}\biggl\|\int_{[0,\infty)}\alpha^k(x) \,\mu^{m}(\d x) - \int_{[0,\infty)}\alpha^\infty(x) \,\mu^{m}(\d x)   \biggr\|_1\\
+\limsup_{k\to\infty}\biggl\|\int_{[0,\infty)}\alpha^\infty(x) \,\mu^{m}(\d x) - \int_{[0,\infty)}\alpha^\infty(x) \,\mu^{\infty}(\d x)   \biggr\|_1\\
+\limsup_{k\to\infty}\biggl\|\int_{[0,\infty)}\alpha^\infty(x) \,\mu^{\infty}(\d x) - \int_{[0,\infty)}\alpha^k(x) \,\mu^{\infty}(\d x)   \biggr\|_1.
\end{multlined}&
\end{flalign*}
Initially, we are going to prove now that the two last summands of the right-hand side of the above inequality are equal to zero.
We start with the second summand and we realise that we need only to use that $\mu^{k} \xrightarrow{\hspace{0.2em}\textup{w}\hspace{0.2em}} \mu^{\infty}$ and \cref{CharactWeakConv} for the special case of a singleton.
The third summand is also equal to zero as an outcome of the bounded convergence theorem.
Indeed, we have that the sequence $\mathcal A$ is uniformly bounded and by \cite[Proposition VI.2.1]{jacod2003limit} we have the pointwise convergence $\alpha^k(x)\xrightarrow{\hspace{1em}} \alpha^\infty(x)$ for $k\longrightarrow\infty,$ at every point $x$ which is a point of continuity of $\alpha^\infty$.
Since the set $\{x\in\mathbb{R}_+: \Delta \alpha^\infty(x)\ne 0\}$ is at most countable, we can conclude it is an $\mu^{\infty}$--null set.
Therefore, since every element of $\mathcal A$ is Borel-measurable, we can conclude.

\medskip
Finally, we deal with the first summand.
Let us provide initially some auxiliary results.
Using again \cite[Proposition VI.2.1]{jacod2003limit} we have that 
\begin{equation}\label{prop:MOAppl-1}
	\limsup_{k\to\infty} \big\|\alpha^k(x) - \alpha^\infty(x)\big\|_1
	\le \big\|\Delta \alpha^\infty(x)\big\|_1 \text{ for every }x\in\mathbb{R}_+,
\end{equation}
because the only two possible accumulation points of the sequence $\big(\alpha^{k,i}(x)\big)_{k\in\mathbb N}$ are $\alpha^{\infty,i}(x)$ and $\alpha^{\infty,i}(x-)$, for every $i\in\{1,\dots,p\}$ and every $x\in\mathbb{R}_+$.
Now observe that the function $\mathbb{R}_+\ni x\longmapsto \big\|\Delta \alpha^\infty(x)\big\|_1 \in\mathbb{R}_+$ is Borel measurable and $\mu^{\infty}$--almost everywhere equal to the zero function.
This observation allows us to apply \citeauthor{mazzone1995characterization} \cite[Theorem 1]{mazzone1995characterization} in order to obtain the convergence 
\begin{equation}\label{prop:MOAppl-2}
	\int_{[0,\infty)}\big\|\Delta \alpha^\infty(x)\big\|_1 \,\mu^{m}(\d x)
	\xrightarrow[m\to\infty]{|\cdot|}
	\int_{[0,\infty)}\big\|\Delta \alpha^\infty(x)\big\|_1 \,\mu_{\infty}(\d x)=0.
\end{equation}
In view of the above and the boundedness of the sequence $\mathcal A$, we can apply (the reverse) Fatou's lemma and we obtain for every $m\in\mathbb N$
\begin{multline*}
\limsup_{k\to\infty}
\biggl\|\int_{[0,\infty)}\alpha^k(x) \,\mu^{m}(\d x) - \int_{[0,\infty)}\alpha^\infty(x) \,\mu^{m}(\d x)   \biggr\|_1\\
\le\hspace{1em} 
\int_{[0,\infty)}
\limsup_{k\to\infty}
\Big\|\alpha^k(x) - \alpha^\infty(x)\Big\|_1 \,\mu^{m}(\d x)
	\overset{\eqref{prop:MOAppl-1}}{\le} \int_{[0,\infty)}\big\|\Delta \alpha^\infty(x)\big\|_1 \,\mu^{m}(\d x).
\end{multline*}
Now, we can conclude that the Condition \ref{MOBii} of \cref{MooreOsgoodB} is indeed satisfied by combining the above bound with Convergence \eqref{prop:MOAppl-2}.
\end{proof}
\begin{remark}
Let us adopt the notation of {\rm\cref{CharactWeakConv}} and {\rm\cref{prop:MOAppl}} and assume that the measurable space is $\bigl([0,T], \mathcal B([0,T])\bigr)$, for some $T\in\mathbb{R}_+$.
We claim that we can adapt the aforementioned results without loss of generality, which can be justified as follows.
The limit measure $\mu^{\infty}$ is atomless, so the generality is not harmed if in {\rm\cref{CharactWeakConv}} the integrands $\alpha^k,$ for some $k\in\overline{\mathbb N}$, have a jump at point $T$.
Indeed, by the weak convergence $\mu^{k}\xrightarrow{\textup{w}}\mu^{\infty}$ and {\rm\cref{EquivWeakConv}} we have that the sequence of distribution functions $(F^{k})_{k\in\overline{\mathbb N}}$ which is associated to the sequence $(\mu^k)_{k\in\overline{\mathbb N}}$ converges uniformly.
Therefore, we can easily conclude that 
\begin{equation*}
 	\limsup_{k\to\infty}\big\|\Delta \alpha^k(T)\big\|_2\mu^{k}(\{T\}) 
 	= \limsup_{k\to\infty}\big\|\Delta \alpha^k(T)\big\|_2\Delta F^{k}_T 
 	\le \big\|\Delta \alpha^\infty(T)\big\|_2\Delta F^{\infty}_T =0. 
\end{equation*} 
In other words, we can simply assume that the distribution functions are constant after time $T$ in order to reduce the general case to the compact-interval case.
\end{remark}
The following corollary is almost evident due to the fact that the limit measure is atomless.
However we state it in the form we will need it in \cref{prop:p-Step-WeakApprox} and we provide its complete (rather trivial) proof.
\begin{corollary}\label{cor:MOAppl}
Consider the measurable space $\bigl(\mathbb{R}_+, \mathcal B(\mathbb{R}_+)\bigr)$ and let $(\mu^{k})_{k\in\overline{\mathbb N}}$ be a sequence of finite measures, where $\mu^k(\{0\})=0$ for every $k\in\overline{\mathbb N}$ and $\mu^{\infty}$ is atomless.
Let, moreover, $\mathcal A:=(\alpha^k)_{k\in\overline{\mathbb N}}\subset\mathbb{D}(\mathbb{R}_+;\mathbb{R}^p)$ and the following to be true
\begin{enumerate}
\item\label{WeakConvergence-cor}	$\mu^{k} \xrightarrow{\hspace{0.2em}\textup{w}\hspace{0.2em}} \mu^{\infty};$
\item\label{UnifBound-cor} 			$\mathcal A$ is uniformly bounded, \emph{i.e.} $\sup_{\alpha\in\mathcal{A}} \|\alpha\|_{\infty}<\infty;$ 
\item\label{equircontin-cor} 		$\mathcal A$ is $\delta_{\textup{J}_1(\mathbb{R}^p)}$-convergent. 
\end{enumerate}
Then 
\begin{equation}
\int_{[0,\cdot]} \alpha^k(x)\,\mu^{k}(\d x)
\xrightarrow[k\to\infty]{\hspace{0.2em}\textup{J}_1(\mathbb{R}^p)\hspace{0.2em}}
 \int_{[0,\cdot]} \alpha^\infty(x)\,\mu^{\infty}(\d x).
\label{UnifIntConv-cor}
\end{equation}
\end{corollary}
\begin{proof}
Let $F^{k}$ denote the distribution function associated to the measure $\mu^{k}$ for every $k\in\overline{\mathbb N}$ and  
\begin{equation*}
	\gamma^k(t):=\int_{[0,t]}\alpha^k(x)\,\mu^{k}(\d x)=\int_{[0,t]}\alpha^k(s)\,\d F^{k}_s\in\mathbb{D}(\mathbb{R}_+;\mathbb{R}^p) \text{ for }k\in\overline{\mathbb N}.
\end{equation*}
We are going to apply \citeauthor*{ethier2005markov} \cite[Proposition 3.6.5]{ethier2005markov} in order to prove the required convergence.
To this end, let us fix a $t^\infty\in\mathbb{R}_+$ and a sequence $(t^k)_{k\in\mathbb N}$ such that $t^k\xrightarrow{\hspace{1em}}t^\infty$ as $k\to\infty.$
We need to prove that the requirements \textup{(a)}--\textup{(c)} of the aforementioned proposition are satisfied.
However,  due to the continuity of the limit, these three requirements reduce to a single one: for any $t^{\infty}\in\mathbb{R}_+$ and for any sequence $(t^{k})_{k\in\mathbb{N}}$ such that $t^{k}\xrightarrow[k\to\infty]{|\cdot|}t^{\infty}$ holds
\begin{align}\label{conv:appCorollary}
	\lim_{k\to\infty}\big\|\gamma^k(t^k) - \gamma^\infty(t^\infty)\big\|_1=0.
\end{align}
Let us fix a $t^{\infty}\in\mathbb{R}_+$ and consider a sequence $(t^{k})_{k\in\mathbb{N}}$ such that $t^{k}\xrightarrow[k\to\infty]{|\cdot|}t^{\infty}$. 
Then, by \cite[Theorem VI.2.15.c)]{jacod2003limit}, we can easily conclude that $F^{k}_{t^{k}\wedge \cdot}\xrightarrow[k\to\infty]{\textup{J}_1(\mathbb{R}_+;\mathbb{R})}F^{k}_{t^{k}\wedge \cdot}$.
We can now immediately verify by means of \cref{prop:MOAppl} that \eqref{conv:appCorollary} is true. 
\end{proof}

\begin{remark}\label{rem:MOApplication}
Let us borrow the notation of {\rm\cref{prop:MOAppl}}, with $F^{k}$ the distribution function of the measure $\mu^{k}$, for $k\in\overline{\mathbb{N}}$, where $F^{\infty}$ is not necessarily continuous.
It is well-known in the literature of the Skorokhod $\textup{J}_1$-convergence that a sufficient condition for the convergence of the integral functions is the joint convergence of the integrands and the integrators, \emph{i.e.}, if 
$(\alpha^{k},F^{k})\xrightarrow[k\to\infty]{\textup{J}_1}(\alpha^{\infty},F^{\infty})$, then 
$\alpha^{k}\cdot F^{k} \xrightarrow[k\to\infty]{\textup{J}_1} \alpha^{\infty}\cdot F^{\infty}$, where $\cdot$ denotes integration here.
When $F^{\infty}$ is continuous, the joint convergence is trivially obtained because of the continuity of the limit integrator $F^{\infty}$.
In other words, we could have proved {\rm\cref{cor:MOAppl}} without using {\rm\cref{prop:MOAppl}}.
However, the purpose of presenting this approach is twofold: firstly to present the characterisation of weak convergence of measures as stated in {\rm\cref{CharactWeakConv}} and secondly to prove that under the assumption $(\alpha^{k},F^{k})\xrightarrow[k\to\infty]{\textup{J}_1}(\alpha^{\infty},F^{\infty})$, it is true that
$\alpha^{k}\cdot F^{l} \xrightarrow[(k,l)\to(\infty,\infty)]{\textup{J}_1} \alpha^{\infty}\cdot F^{\infty}$, which is stronger than 
$\alpha^{k}\cdot F^{k} \xrightarrow[k\to\infty]{\textup{J}_1} \alpha^{\infty}\cdot F^{\infty}$.
In particular
\begin{align*}
\lim_{k\to\infty}\lim_{l\to\infty} \alpha^{k}\cdot F^{l} 
=\lim_{l\to\infty}\bigg\{\limsup_{k\to\infty} \alpha^{k}\cdot F^{l} - \liminf_{k\to\infty} \alpha^{k}\cdot F^{l}\bigg\}
=\lim_{k\to\infty}\alpha^{k}\cdot F^{k} 
= \alpha^{\infty}\cdot F^{\infty}.
\end{align*}
The next natural question is whether a characterisation similar to {\rm\cref{CharactWeakConv}} holds if we consider a weakly-convergent sequence of measures to a limit-measure with atoms.
In view of {\rm\citeauthor{mazzone1995characterization} \cite[Theorem 1]{mazzone1995characterization}}, one should not expect a characterisation based on arbitrary relatively compact sets of the Skorokhod space. 

\end{remark}
%
%
%
The following lemma will be helpful in proving that a relatively compact set of $\mathbb{D}(\mathbb{R}_+;\mathbb{R}^p)$ is uniformly bounded.
Observe that every relatively compact subset of $\mathbb{D}(\mathbb{R}_+;\mathbb{R}^p)$  satisfies \ref{LocBound} of the following lemma; see \cite[Theorem VI.1.14.b).(i)]{jacod2003limit}.
\begin{lemma}\label{UniformlyBounded}
Let $(\alpha^k)_{k\in\mathbb{N}}\subset\mathbb{D}\bigl(\mathbb{R}^p\bigr)$ such that
\begin{enumerate}
\item\label{LimitAtInfty} 	$\alpha^k(\infty):=\lim_{t\to\infty}\alpha^k(t)$ exists and $\big\|\alpha^k(\infty)\|_2<\infty$ for every $k\in\mathbb{N};$ 
\item\label{LocBound} 		$\sup_{k\in\mathbb{N}}\sup_{t\in[0,N]} \bigl\|\alpha^k(t)\bigr\|_2<\infty$ for every $N\in\mathbb{N};$ 
\item\label{LimsupBound} 	$\limsup_{(k,t)\to(\infty,\infty)} \bigl\|\alpha^k(t)\big\|_2<\infty$. 
\end{enumerate}
Then the sequence $(\alpha^k)_{k\in\mathbb{N}}$ is uniformly bounded, \emph{i.e.} $\sup_{k\in\mathbb N}  \|\alpha\|_\infty<\infty$.
\end{lemma}
\begin{proof}
Assume that the sequence is unbounded, \emph{i.e.} for every $M\in\mathbb N$ there exists $k(M)\in\mathbb{N}$ such that $\|\alpha^{k(M)}\|_\infty>M$. 
We can extract a subsequence of $\bigl(k(M)\bigr)_{M\in\mathbb N}$, called $\bigl(k(M_n)\bigr)_{n\in\mathbb N}$,  such that $k(M_l)< k(M_m)$ whenever $l< m$ and $M_n\to\infty$ as $n\to\infty$. 
This can be done as follows.
For $n=1$ we define $k(M_1):=k(1)$ and for $i\in\mathbb N\setminus\{1\}$ we set
\begin{equation*}
	k(M_{i+1}):=\min\Bigl\{k(M)\in\mathbb N: \|\alpha^{k(M)}\|_\infty> \Bigl\lceil \| \alpha^{k(M_i)}\|_\infty\Bigr\rceil\Bigr\},
\end{equation*}
where $\lceil x\rceil$ denotes the least integer greater than or equal to $x.$
Since we have assumed that for every $k\in\mathbb N$ the value $\big\|\alpha^k(\infty)\|_2$ is finite, we have in particular (using that $\alpha^k\in\mathbb{D}(\mathbb{R}_+;\mathbb{R}^p)$) that $\|\alpha^k\|_\infty<\infty.$
Therefore
\[
\Bigl\{k(M)\in\mathbb N: \|\alpha^{k(M)}\|_\infty> \Bigl\lceil \| \alpha^{k(M_i)}\|_\infty\Bigr\rceil\Bigr\}\neq \emptyset,
\]
 and its minimum exists as $\mathbb N$ is a well-ordered set under the usual order.
In view of these comments, the subsequence $\bigl(k(M_n)\bigr)_{n\in\mathbb N}$ is well-defined and, according to our initial assumption, it has to be unbounded.

By definition of $\alpha^{k(M)}$, there exists $t_{k(M)}\in\mathbb{R}_+$ such that $\|\alpha^{k(M)}(t_{k(M)})\|_\infty>M$, a property which holds in particular for the subsequence indexed by $(k(M_n))_{n\in\mathbb N}.$
Let us distinguish, now, two cases for the sequence $(t_{k(M_n)})_{n\in\mathbb N}$.
\begin{enumerate}[label=\raisebox{0.5pt}{\text{\tiny$\blacksquare$}},itemindent=0.cm, leftmargin=*]
	\item If $(t_{k(M_n)})_{n\in\mathbb N}$ is bounded, then \ref{LocBound} leads to a contradiction, for $N:=\sup\{t_{k(M_n)},{n\in\mathbb N}\}$.
	\item If $(t_{k(M_n)})_{n\in\mathbb N}$ is unbounded, then there exists a subsequence $(t_{k(M_{n_l})})_{l\in\mathbb N}$ such that
$t_{k(M_{n_l})}\xrightarrow{\hspace{1em}}\infty$ as $l\rightarrow\infty$. 
Therefore, we have also
\begin{equation*}
	\lim_{l\to\infty} \big\|\alpha^{k(M_{n_l})}(t_{k(M_{n_l})})\big\|_\infty\ge \lim_{l\to\infty} M_{n_l}=\infty.
\end{equation*}
But, then \ref{LimsupBound} leads to a contradiction.
\end{enumerate}

Now, we have only to verify that the set $\{\|\alpha^k\|_\infty:k\in\{1,\dots,k(1)\}\}$ is bounded in order to conclude that the sequence $(\alpha^k)_{k\in\mathbb N}$ is $\|\cdot\|_\infty$-bounded. 
But the above is clear since it is a finite set of finite numbers; use that $\|\alpha^k(\infty)\|_2<\infty$ and that $\alpha^k\in\mathbb{D}(\mathbb{R}_+;\mathbb{R}^p).$
\end{proof}

%% file: Appendix_SBSDE_1.tex

\subsection{Some helpful lemmata for the proof of Theorem \ref{BSDERobMainTheorem}.}

\begin{lemma}\label{lem:Construction_DCirc}
There exists a countable set $D^{\circ,\ell\times m}\subset C_{c}(\mathbb{R}_+;\mathbb{R}^{\ell\times m})$ such that for every $Z\in\mathbb{H}^{\infty,\circ}$ and for every $\varepsilon >0$, there exists $\widetilde{Z}^{\varepsilon}\in D^{\circ,\ell \times m}$ with the property 
\begin{align*}
\int_{(0,\infty)} \textup{Tr}\bigg[ (Z_s- \widetilde{Z}^{\varepsilon}_s) \frac{\d \langle X \rangle_s}{\d C_s} (Z_s - \widetilde{Z}^{\varepsilon}_s)^{\top}\bigg] \d C_s<\varepsilon,\; \mathbb{P}\text{\rm--a.s.}
\end{align*}
\end{lemma}

\begin{proof}
Since $Z\in \mathbb{H}^{\infty,\circ}$, then 
\begin{align*}
\int_{(0,\infty)} \textup{Tr}\bigg[ Z_s \frac{\d \langle X \rangle_s}{\d C_s} Z_s^{\top}\bigg] \d C_s<\infty,\; \mathbb{P}\text{\rm--a.s.}
\end{align*}
Let $\widetilde{\Omega}\subset \Omega$ be the set of $\omega$ for which the aforementioned property holds and $\langle X\rangle$ is finite.
Then, for $\overline{Z}^{n}:=Z \mathds{1}_{\{\Vert Z\Vert_{2}\le n\}}$ holds 
\begin{align*}
\int_{(0,\infty)} \textup{Tr}\bigg[ (Z_s- \overline{Z}^{n}_s) \frac{\d \langle X \rangle_s}{\d C_s} (Z_s - \overline{Z}^{n}_s)^{\top}\bigg] \d C_s\xrightarrow[n\to\infty]{} 0,\; \forall \omega \in \widetilde{\Omega}.
\end{align*}
Therefore, for every $\varepsilon >0$, there exists $n_0$ (depending on $\omega$) such that 
\begin{align}\label{ineq:epsilonhalf-1}
\int_{(0,\infty)} \textup{Tr}\bigg[ (Z_s- \overline{Z}^{n}_s) \frac{\d \langle X \rangle_s}{\d C_s} (Z_s - \overline{Z}^{n}_s)^{\top}\bigg] \d C_s<\frac{\varepsilon}{2},\; \forall n\ge n_0, \forall \omega \in \widetilde{\Omega}.
\end{align}
The fact that $\overline{Z}^{n}$ is bounded, allows us to write 
\begin{align*}
\int_{(0,\infty)} \textup{Tr}\bigg[\, \overline{Z}^{n}_s \frac{\d \langle X \rangle_s}{\d C_s} \overline{Z}_{s}^{n\top}\bigg] \d C_s 
=\sum_{k=1}^{\ell} \sum_{i,j=1}^{m} 
\int_{(0,\infty)} \overline{Z}_{s}^{n,ki}\overline{Z}_{s}^{n,kj}  \d \langle X^{\infty,\circ,i},X^{\infty,\circ,j}\rangle_s ,
\; \forall \omega \in \widetilde{\Omega},
\end{align*}
where we have denoted by $\overline{Z}_{s}^{n,ki}$ the $ki$-element of the matrix $\overline{Z}^{n}_s$, for $k\in \{1,\dots,\ell\}$ and $i\in\{1,\dots,m\}$. 
Kunita--Watanabe's inequality (in conjunction with Young's Inequality) writes for all $k\in\{1,\dots,\ell\}$, $(i,j)\in\{1,\dots,m\}^2$, and $\omega \in \widetilde{\Omega}$
\begin{align}\label{ineq:epsilonhalf-2}
\int_{(0,\infty)} \overline{Z}_{s}^{n,ki}\overline{Z}_{s}^{n,kj}  \d \textup{Var}\big[\langle X^{\infty,\circ,i},X^{\infty,\circ,j}\rangle\big]_s
\leq  
\frac{1}{2}\int_{(0,\infty)} \vert \overline{Z}_{s}^{n,ki}\vert^{2}   \d \langle X^{\infty,\circ,i}\rangle_s
+\frac{1}{2}\int_{(0,\infty)} \vert \overline{Z}_{s}^{n,kj}\vert^{2}   \d \langle X^{\infty,\circ,j}\rangle_s.
\end{align}
Now we use the fact that $\langle X^{\infty,\circ,i}\rangle$ is a finite, Borel measure on $\mathbb{R}_+$.
Therefore, the set 
\begin{align*}
\bigg\{f:\mathbb{R}_+\to\mathbb{R}: f=\sum_{i=1}^{p}\alpha_i\mathds{1}_{I_i}, \text{ with } \alpha_i\in\mathbb{Q},\; I_i \text{ bounded interval with rational endpoints } \forall i\in\{1,\dots,p \}\bigg\},
\end{align*}
is countable and dense in $\mathbb{L}^{2}(\mathbb{R}_+,\mathcal{B}(\mathbb{R}_+), \langle X^{\infty,\circ,i}\rangle)$, for every $i\in\{1,\dots,m\}$. 
The regularity of the measure $\langle X^{\infty,\circ,i}\rangle$ (as finite, Borel defined on a locally compact space for every $\omega\in\widetilde{\Omega}$) and Tietze's extension theorem allow\footnote{The local compactness allows for the following property: for every $K$ compact and $U$ open, there exists open $V$ such that its closure $\overline{V}$ is compact and $ K\subset V\subset \overline{V}\subset U$. Then, Tietze's extension theorem guarantees the existence of $f$ in $C_{c}(\mathbb{R}_+;\mathbb{R})$ such that $\mathds{1}_K\le f \le \mathds{1}_U$ with $\textup{supp}(f)$ compact.} us to choose a countable subset of $C_{c}(\mathbb{R}_+;\mathbb{R}^{\ell\times m})$, denoted by $D^{\circ}$, which is dense in $\mathbb{L}^{2}(\mathbb{R}_+,\mathcal{B}(\mathbb{R}_+), \langle X^{\infty,\circ,i}\rangle)$ for every $\omega\in\widetilde{\Omega}$.
We can easily pass to the required set $D^{\circ,\ell\times m}$ by observing now that for every $k\in\{1,\dots,\ell\}$ and for every $i\in\{1,\dots, m\}$, for every $\omega\in\widetilde{\Omega}$, we can choose $\widetilde{Z}^{\varepsilon,ki}$\footnote{We suppress the dependence on $\omega$.} such that 
\begin{align*}
\int_{(0,\infty)} \big\vert\overline{Z}^{n_0,ki}_s- \widetilde{Z}^{\varepsilon,ki}_s\big\vert^{2} \d \langle X^{\infty,\circ,i}\rangle_s<\frac{\varepsilon}{2\ell m^2},\; \forall \omega\in\widetilde{\Omega},
\end{align*}
which in conjunction with \eqref{ineq:epsilonhalf-1} and \eqref{ineq:epsilonhalf-2} results in 
\begin{align*}
\int_{(0,\infty)} \textup{Tr}\bigg[ (Z_s- \widetilde{Z}^{\varepsilon}_s) \frac{\d \langle X \rangle_s}{\d C_s} (Z_s - \widetilde{Z}^{\varepsilon}_s)^{\top}\bigg] \d C_s<\varepsilon,\; \forall \omega \in \widetilde{\Omega}.
\end{align*}
\end{proof}

\begin{lemma}\label{lem:Construction_DNatural}
Let $\nu$ be a measure on $(E,\mathcal{B}(E))$ such that $\int_{\mathbb{R}_+\times \mathbb{R}^{n}} \Vert x\Vert_2^{2} \nu(\d t,\d x)<\infty$ and $\nu(E_0)=0$. 
Then, there exists a countable set $D^{\natural}\subset C_{c|\widetilde{E}}(E;\mathbb{R}^{\ell})$ such that $D^{\natural}$ is dense in 
$\mathbb{L}^{2}(E,\mathcal{B}(E),\nu)$. 
\end{lemma}

\begin{proof}
Since $\nu(E_0)=0$, it is sufficient to prove that there exists a countable subset of $C_c(\widetilde{E};\mathbb{R}^{\ell})$ which is dense in 
$\mathbb{L}^{2}(\widetilde{E},\mathcal{B}(\widetilde{E}),\nu|_{\mathcal{B}(\widetilde{E})})$.
The integrability assumption implies that $\nu|_{\mathcal{B}(\widetilde{E})}$ is finite on every compact subset of $\widetilde{E}$.
In view of the above property, $\nu|_{\mathcal{B}(\widetilde{E})}$ is a Radon measure\footnote{In other words, it is finite on compact subsets of $\widetilde{E}$, inner-regular on all open sets and outer-regular on every Borel set.} as a locally finite Borel measure. 
The reader may observe that it is crucial that $\widetilde{E}$ is open, so that we exclude compacts subsets of $E$ which intersect $E_0$.

\medskip

Now that we have verified that $\nu|_{\mathcal{B}(\widetilde{E})}$ is Radon, we can use the well-known fact that the family of the step functions with rational $\mathbb{R}^\ell$-values and supported on rectangles whose vertices' endpoints are rational is $\mathbb{L}^2$-dense in the collection of simple functions, which in turn is dense in $\mathbb{L}^{2}\big(\widetilde{E},\mathcal{B}(\widetilde{E}),\nu|_{\mathcal{B}(\widetilde{E})}\big)$.
We will denote by $\mathcal{Z}$ the aforementioned family of step functions.
Hence, $\mathcal{Z}$ is a countable $\mathbb{L}^2$-dense set.
Let us define
\begin{gather*}
\mathcal{Z}_c := \big\{ f\in\mathcal{Z}: \text{supp}(f)\text{ is compact and the rectangles in the representation are closed} \big\},
\shortintertext{as well as}
\mathcal{Z}_o := \big\{ f\in\mathcal{Z}: \text{the rectangles in the representation are open} \big\}.
\end{gather*}
For every $f_c\in \mathcal{Z}_c$ and for every $f_o\in\mathcal{Z}_o$ such that\footnote{We denote by $f|_A$ the restriction of the function $f$ on the set $A$.}
$f_c|_{\text{supp}(f_c)}=f_o|_{\text{supp}(f_c)}$, there exists by Tietze's extension theorem $f\in C_{c}(\widetilde{E};\mathbb{R}^\ell)$ such that $f_c\le f \le f_o$.
The collection of all $f\in C_{c}(\widetilde{E};\mathbb{R}^\ell)$ with the above property is countable and $\mathbb{L}^2$-dense; we will denote it by $D^{\natural}$.
Without  destroying the properties of $D^{\natural}$, we can assume that the zero function, denoted by $0$, lies in $D^\natural$. 
Now, we can isometrically identify $C_{c}(\widetilde{E};\mathbb{R}^\ell)$ and $C_{c|\widetilde{E}}(E;\mathbb{R}^\ell)$ and conclude by using the fact that $\nu(E_0)=0$.
\end{proof}

\begin{lemma}\label{lem:Uniform-Picard-Constamts}
Under Conditions \ref{BgeneratorUI} and \ref{Bintegrator}.\ref{BintegratorPhiBounds}, there exists 
$k_{\star,0}\in\N$, such that
\begin{align*}
\sup_{k\ge k_{\star,0}} M^{\Phi^k}_{\star}(\hat{\beta}) <\frac{1}{4},
\end{align*}
with 
$\mathbb{R}_+\times \mathbb{R}_+ \ni (\beta,\Phi)\longmapsto M^{\Phi}_{\star}(\beta)$, defined as in {\rm\cite[Lemma 3.4]{papapantoleon2016existence}}. 
\end{lemma}

\begin{proof}
For the convenience of the reader, we restate the notation of the aforementioned lemma.
For $\beta,\Phi>0$ and $\mathcal{C}_\beta:=\{(\gamma,\delta) \in(0,\beta]^2, \gamma<\delta\}$ we define
\begin{gather*}
\Pi_\star^\Phi(\gamma,\delta):=\frac8\gamma+\frac{9}{\delta}+{9\delta}\frac{\textup{e}^{(\delta-\gamma)\Phi}}{\gamma(\delta-\gamma)},
 \text{ with } 
M^\Phi_\star(\beta):=\inf_{(\gamma,\delta)\in\mathcal C_{\beta}}\Pi^\Phi_\star(\gamma,\delta)
.
\end{gather*}
The infimum 
 is attained at $(\bar{\gamma}_\star^\Phi(\beta),\beta)$; the exact values are given in \cite[Lemma 3.4]{papapantoleon2016existence}. 
Hence, using Assumption \ref{Bintegrator}.\ref{BintegratorPhiBounds} and Remark \ref{rem:BSDE-solution}.\ref{foot:StandardData}, we can assume that
there exists 
$\Phi^{k_{\star,0}}$ such that 
	$$M_\star^{\Phi^{k_{\star,0}}}(\hat{\beta})=\Pi_\star^{\Phi^{k_{\star,0}}}(\overline{\gamma}_\star^{\Phi^{k_{\star,0}}},\hat{\beta})<\frac{1}{4}.$$
We observe now that, when we fix $\gamma, \delta\in (0,\beta]^2$, then the function
\begin{equation}\label{funct:decreas}
 \mathbb R_+\ni \Phi\mapsto\Pi^{\Phi}_{\star}(\gamma,\delta)\in\mathbb R_+,
\end{equation}
is decreasing and continuous under the usual topology of $\mathbb R$.
Consequently, by \cite[Lemma 3.4]{papapantoleon2016existence} and the fact that the function in \eqref{funct:decreas} is decreasing, we obtain 
\begin{align*}
 	M_\star^{\Phi^{k}}(\hat{\beta})	
 		= 	\Pi_\star^{\Phi^{k}}(\overline{\gamma}_\star^{\Phi^{k}},\hat{\beta})
 		\le \Pi_\star^{\Phi^{k}}(\overline{\gamma}_\star^{\Phi^{k_{\star,0}}},\hat{\beta}) 
 		\le \Pi_\star^{\Phi^{k_{\star,0}}}(\overline{\gamma}_\star^{\Phi^{k_{\star,0}}},\hat{\beta}) 
 		= M_\star^{\Phi^{k_{\star,0}}}(\hat{\beta})<\frac{1}{4},\text{ for }k\ge k_{\star,0}.
	\numberthis\label{XiStar}
\end{align*}
\end{proof}

\begin{lemma}\label{cor:DoobApplication}
Let $(\zeta^k)_{k\in\overline{\mathbb N}}$ be a sequence of $\mathbb R^p$-valued random variables such that $\bigl(\bigl\|\zeta^k\bigr\|_2^2\bigr)_{k\in\overline{\mathbb N}}$ is uniformly integrable.
Then the sequence $\Bigl(\sup_{t\in\mathbb R_+} \bigl\| \mathbb E\bigl[\zeta^k \big|\mathcal G^k_t\bigr]\bigr\|_2^2\Bigr)$ is uniformly integrable.
\end{lemma}
\begin{proof}
Let $\Phi$ be a moderate Young function\footnote{See \cite[Appendix 2]{papapantoleon2019stability}.} for which the sequence $\bigl(\bigl\|\zeta^k\bigr\|_2^2\bigr)_{k\in\overline{\mathbb N}}$ satisfies the de La Vallée Poussin--Meyer criterion (see \citeauthor*{meyer1978sur} \cite[Lemme]{meyer1978sur}), \emph{i.e.}
\begin{equation}\label{Bound:dlVPM-seq-xik}
	\sup_{k\in\overline{\mathbb N}} \mathbb E\Bigl[ \Phi\bigl( \bigl\|\zeta^k\bigr\|_2^2\bigr) \Bigr]<\infty.
\end{equation}
Then, by \cite[Proposition A.9]{papapantoleon2019stability}, $\Psi:=\Phi\circ \textup{quad}$ is a moderate Young function.
Using the fact that $\Phi$ is increasing, we can write \eqref{Bound:dlVPM-seq-xik} as
\begin{equation}\label{Bound:dlVPM-seq-xik-1}
M:=\sup_{k\in\overline{\mathbb N}}\mathbb E\Bigl[ \Psi\bigl(\|\zeta^k\|_2\bigr)\Bigr]<\infty.
\end{equation}
The latter form will be more convenient for later use. 
Before we proceed to prove the claim of the corollary, we provide some helpful results.
In order to ease notation, let us denote the Orlicz norm of $\Psi$ by $\big\| \cdot \big\|_{\Psi}$, which we simplify for the case $\big\| \|\zeta^k\|_2\big\|_{\Psi}$ by $\Vert \zeta^k\Vert_{2,\Psi}$ for $k\in\overline{\mathbb N}$; see \citeauthor*{rao2002applications} \cite[Subsection III.3.1 Theorem 3, p. 54]{rao2002applications}.
Observe that $\Vert \zeta^k\Vert_{2,\Psi}<\infty$, for every $k\in\overline{\mathbb N}$ because of \eqref{Bound:dlVPM-seq-xik-1}.
We head to prove that $\sup_{k\in\overline{\mathbb N}}\Vert \zeta^k\Vert_{2,\Psi} <\infty.$
To this end, observe that

\medskip
\hspace{0.3em}\raisebox{1pt}{\tiny$\blacksquare$} if $M\le 1$, then $\sup_{k\in\overline{\mathbb N}}\zeta^k\le1$ by the definition of the Orlicz norm;

\medskip
\hspace{0.3em}\raisebox{1pt}{\tiny$\blacksquare$} if $M> 1$, then using the convexity of $\Psi$ and the fact that $\Psi(0)=0$ we obtain
\begin{gather*}
\mathbb E\bigg[ \Psi\bigg(\frac{\|\zeta^k\|_2}{M} \bigg)\bigg] 
\le 
\frac{\mathbb E\bigl[ \Psi\bigl(\|\zeta^k\|_2 \bigr)\bigr] }{M}\overset{\eqref{Bound:dlVPM-seq-xik-1}}{\le} 1, \text{ which implies } \Vert \zeta^k\Vert_{2,\Psi}\le M.
\shortintertext{and consequently}
	\sup_{k\in\overline{\mathbb N}} \Vert \zeta^k\Vert_{2,\Psi}\le 1\vee  M\overset{\eqref{Bound:dlVPM-seq-xik-1}}{<}\infty.
\numberthis\label{Bound:Xik}
\end{gather*}

We proceed now to prove the uniform integrability of $\bigl(\sup_{t\in\mathbb R_+} \bigl\|\mathbb E\big[ \zeta^k\big| \mathcal G^k_t\bigr] \bigr\|_2^2\bigr)_{k\in\overline{\mathbb N}}$.
It suffices to prove that 
\[
\sup_{k\in\overline{\mathbb N}}\mathbb E\bigg[ \Phi\bigg(\sup_{t\in\mathbb R_+}\bigl\| \mathbb E[\zeta^k | \mathcal G_t^k]\bigr\|_2^2\bigg) \bigg]<\infty,
\]
or equivalently
that $\sup_{k\in\overline{\mathbb N}}\mathbb E\Bigl[ \Psi\bigr(\sqrt{2}S^k\bigr) \Bigr]<\infty,$
where we have defined $S^k:=\sup_{t\in\mathbb R_+}\big\|\mathbb E[\zeta^k|\mathcal G^k_t]\bigr\|_2$.
Using the properties of moderate Young functions, \emph{e.g.} see \citeauthor{long1993martingale} \cite[Theorem 3.1.1, p. 82]{long1993martingale}, we can obtain 
\begin{equation}\label{Bound:Exp-Sup-OrlNorm}
\mathbb E\Bigl[ \Psi\bigr(\sqrt{2}S^k\bigr) \Bigr]
\le 1\vee \widebar{c}_{_\Psi}^{\Vert\sqrt{2} S^k\Vert_\Psi},
\end{equation}
where $\widebar{c}_{_{\Psi}}$ depends only on $\Psi$\footnote{Actually, for a Young function $\Theta$ it is defined $\widebar{c}_{_\Theta}:=\sup_{x>0}\frac{x \theta(x)}{\Theta(x)}$ for $\theta(x)$ the right derivative of $\Theta$. $\widebar{c}_{_\Theta}$ is finite if and only if $\Theta$ is moderate; see \cite[Theorem 3.1.1 (c), p. 82]{long1993martingale}. Analogously, we define $\underline{c}_{_\Theta}:=\inf_{x>0}\frac{x \theta(x)}{\Theta(x)}$. }.
Recall that $\Psi=\Phi\circ\textup{quad}$ is a moderate Young function. 
Hence, for $\Psi^{\star}$ the conjugate Young function of $\Psi$, we have by \cite[Proposition A.9]{papapantoleon2019stability} that $\Psi^\star$ is also moderate.
Therefore, $\widebar{c}_{_\Psi},$ $\widebar{c}_{_{\Psi^\star}}\in(1,\infty)$; we used that $\widebar{c}_{_{\Psi^\star}}=\underline{c}_{_\Psi}$.
Consequently, we can conclude the required property if we prove that $\sup_{k\in\overline{\mathbb N}}\|S^k\|_{\Psi}<\infty,$ where the $\sqrt{2}$ can be taken out, since $\|\cdot\|_\Psi$ is a norm.
To this end we will use Doob's $\mathbb L^{\Upsilon}$-inequality.
We can obtain now, by standard properties of norms and the fact that 
\[
\Vert x\Vert_2=\bigg(\sum_{i=1}^{\ell}|x_i|^{2}\bigg)\frac{1}{2}\le \sum_{i=1}^{\ell}|x_i|=\Vert x\Vert_1,
\]
and consequently $S^k\le \sup_{t\in\mathbb R_+}\bigl\|\mathbb E[\zeta^k|\mathcal G^k_t]\bigr\|_1$, the following inequalities
\begin{align*}
\big\| S^k\big\|_{\Psi} 
	&\le \Bigl\| \sup_{t\in\mathbb R_+}\bigl\|\mathbb E[\zeta^k|\mathcal G^k_t]\bigr\|_1\Bigr\|_{\Psi}
	\le \sum_{j=1}^p \Bigl\| \sup_{t\in\mathbb R_+}\bigl|\mathbb E[\zeta^{k,j}|\mathcal G^k_t]\bigr|\Bigr\|_{\Psi}
	\overset{\text{Doob's In.}}{\le}
		2\widebar{c}_{_{\Psi^\star}} \sum_{j=1}^p \bigl\| \zeta^{k,i}\bigr\|_{\Psi}\\
	&\le 2p\widebar{c}_{_{\Psi^\star}} \bigl\| \|\zeta^k\|_2 \bigr\|_{\Psi}
	 =2p\widebar{c}_{_{\Psi^\star}} \Vert \zeta^k\Vert_{2,\Psi},
\end{align*}
where we applied Doob's inequality for moderate Young functions; see \citeauthor*{dellacherie1982probabilities} \cite[Paragraph VI.103, p. 169]{dellacherie1982probabilities}.
The above inequality and \eqref{Bound:Xik} imply the desired
\begin{equation*}
	\sup_{k\in\overline{\mathbb N}}	\big\| S^k\big\|_{\Psi} 
		\le 2p\widebar{c}_{_{\Psi^\star}} \sup_{k\in\overline{\mathbb N}} \Vert \zeta^k\Vert_{2,\Psi}<\infty,
\end{equation*}
which, in view of \eqref{Bound:Exp-Sup-OrlNorm}, implies also the finiteness of $\sup_{k\in\overline{\mathbb N}}\mathbb E\bigl[ \Psi(\sqrt 2 S^k)\bigr]$.
Therefore, the sequence 
\begin{align*}
\bigg(\sup_{t\in\mathbb R_+} \bigl\| \mathbb E\bigl[\zeta^k \big|\mathcal G^k_t\bigr]\bigr\|_2^2\bigg)_{k\in\overline{\mathbb{N}}},
\end{align*}
is uniformly integrable since it satisfies the de La Vallée Poussin theorem for the Young function $\Phi$.
\end{proof}